\documentclass[11pt]{amsart}
\usepackage{amscd, amsfonts, amsmath, amssymb, amsthm, enumerate, epsfig, flexisym, float, graphicx, hhline, hyperref, pdfpages, subcaption, tabularray, tikz, url}
\usepackage[autostyle]{csquotes}
\usepackage[all,cmtip]{xy}

\makeatletter
\@namedef{subjclassname@2020}{\textup{2020} Mathematics Subject Classification}
\makeatother

\newtheorem{theorem}{Theorem}[section]

\newtheorem{conjecture}[theorem]{Conjecture}
\theoremstyle{definition}

\newtheorem{definition}[theorem]{Definition}

\newtheorem{remark}[theorem]{Remark}
\numberwithin{equation}{subsection}
\newtheorem*{ack}{Acknowledgement}

\newcommand{\Conj}{\operatorname{Conj}}
\newcommand{\Core}{\operatorname{Core}}

\DeclareFontFamily{U}{mathb}{\hyphenchar\font45}
\DeclareFontShape{U}{mathb}{m}{n}{<5> <6> <7> <8> <9> <10> <10.95> <12> <14.4> <17.28> <20.74> <24.88> mathb10}{}
\DeclareSymbolFont{mathb}{U}{mathb}{m}{n}
\DeclareMathSymbol{\blackdiamond}{2}{mathb}{"0C}

\usetikzlibrary{cd, arrows, knots, patterns, positioning, shapes.callouts, shapes.geometric, shapes.symbols}
\tikzset{>=stealth', rect1/.style={rectangle, draw=black, text width=38mm, minimum height=17mm, text centered},
rect2/.style={rectangle, draw=black, text width=55mm, minimum height=28mm, text centered},
arrow/.style={->}}
\begin{document}

\title[Biorderability of knot quandles]{Biorderability of knot quandles of knots up to eight crossings}
\author[V. Gupta]{Vaishnavi Gupta}
\author[H. Raundal]{Hitesh Raundal}
\address{Department of Mathematics\\ Faculty of Science\\ University of Allahabad\\ Colnelganj\\ George Town\\ Prayagraj\\ Uttar Pradesh 211002\\ India}
\email{vaishnavigupta7779@gmail.com}
\email{hiteshrndl@gmail.com}


\begin{abstract}
The paper investigates biorderability of knot quandles of prime knots up to eight crossings. We prove that knot quandles of knots $6_3$, $8_7$, $8_8$, $8_{10}$ and $8_{16}$ can not be biorderable. However, we see that knot quandles of knots $4_1$, $6_1$, $6_2$, $7_6$, $7_7$, $8_1$, $8_2$, $8_3$, $8_4$, $8_5$, $8_6$, $8_9$, $8_{11}$, $8_{12}$, $8_{13}$, $8_{14}$, $8_{17}$, $8_{18}$, $8_{20}$ and $8_{21}$ could be biorderable. We also give linear orders on the generating set of the knot quandle of a knot (among these knots) that could be extendable to biorders on the quandle.
\end{abstract}

\maketitle

\section{Introduction}\label{introduction}

Quandles are algebraic systems with a binary operation that encodes the three Reidemeister moves of planar diagrams of knots in the 3-space. These objects have shown appearance in a wide spectrum of mathematics including knot theory \cite{j1,j2,m}, group theory, mapping class groups \cite{Zablow2003,Zablow1999}, set-theoretic solutions of the quantum Yang-Baxter equation \cite{e}, Riemannian symmetric spaces \cite{Loos} and Hopf algebras \cite{ag}, to name a few. Though quandles already appeared under different forms in the literature, study of these objects gained true momentum after the fundamental works of Joyce \cite{j1,j2} and Matveev \cite{m}, who showed that knot quandles are complete invariants of knots up to orientation of the ambient space. Although quandles are strong invariants of knots, the isomorphism problem for them is hard. This has motivated search for newer properties, constructions and invariants of quandles themselves.

It is known that existence of a linear order on a group has profound implications on its structure. For instance, a left-orderable group cannot have torsion and a bi-orderable group cannot have even generalized torsion. Concerning groups arising in topology, the literature shows that many such groups are left-orderable. Braid groups are the most relevant examples of left-orderable groups which are not bi-orderable \cite{d}. On the other hand, pure braid groups are known to be bi-orderable \cite{fr}. Rourke and Wiest \cite{rw} extended this result by showing that mapping class groups of all Riemann surfaces with non-empty boundary are left-orderable. In general these groups are not bi-orderable. Orderability of 3-manifold groups has been investigated extensively where left-orderability is a rather common property. Concerning knot complements, it is known that fundamental groups of knot complements are left-orderable \cite{brw}, whereas fundamental groups of not all knot complements are bi-orderable \cite{pr}. For more on the literature, we refer the reader to the recent monograph \cite{cr} by Clay and Rolfsen which explores orderability of groups motivated by topology, like fundamental groups of surfaces and 3-manifolds, braid and mapping class groups, groups of homeomorphisms, etc.

The notion of orderability can be defined for magmas just as for groups. Since quandles represent interesting examples of non-associative magmas and knot quandles are deeply related to knot groups, it seems natural to explore orderability of quandles. Orderability of magmas including conjugation quandles of bi-orderable groups has been analysed in the work of Dabkowska et al \cite{DDHPV}. A work of Bardakov et al \cite{bps} investigated orderability of quandles for studying zero-divisors in quandle rings. In a recent work \cite{rss}, it is proved that knot quandles of certain positive (or negative) knots are not bi-orderable and almost all torus knot quandles are not right-orderable. Interconnections between orderability of quandles and that of their enveloping groups also has been explored. The results prove that orderability of knot quandles behave quite differently than that of corresponding knot groups. For instance, the knot quandle of the trefoil knot is neither left nor right orderable, whereas knot quandles of many fibered prime knots are right-orderable. Further, in \cite{drs}, it is proved that free involutory quandles are left orderable and that certain generalized Alexander quandles are bi-orderable.

In this paper, we investigate biorderability of knot quandles of prime knots up to eight crossings. It is known that knot quandles of knots $3_1$, $5_1$, $5_2$, $7_1$, $7_2$, $7_3$, $7_4$, $7_5$, $8_{15}$ and $8_{19}$ are not biorderable (see Corollary 6.7, Theorem 7.2 and Corollary 7.4 in \cite{rss}). We explicitly prove that knot quandles of knots $6_3$, $8_7$, $8_8$, $8_{10}$ and $8_{16}$ can not be biorderable (see tables \ref{tbl6_3}, \ref{tbl8_7}, \ref{tbl8_8}, \ref{tbl8_10} and \ref{tbl8_16}). However, we see that knot quandles of knots $4_1$, $6_1$, $6_2$, $7_6$, $7_7$, $8_1$, $8_2$, $8_3$, $8_4$, $8_5$, $8_6$, $8_9$, $8_{11}$, $8_{12}$, $8_{13}$, $8_{14}$, $8_{17}$, $8_{18}$, $8_{20}$ and $8_{21}$ could be biorderable (see tables \ref{tbl4_1}, \ref{tbl6_1}, \ref{tbl6_2}, \ref{tbl7_6}, \ref{tbl7_7}, \ref{tbl8_1}, \ref{tbl8_2}, \ref{tbl8_3}, \ref{tbl8_4}, \ref{tbl8_5}, \ref{tbl8_6}, \ref{tbl8_9}, \ref{tbl8_11}, \ref{tbl8_12}, \ref{tbl8_13}, \ref{tbl8_14}, \ref{tbl8_17}, \ref{tbl8_18}, \ref{tbl8_20} and \ref{tbl8_21}).  We also give linear orders on the generating set of the knot quandle of a knot (among these knots) that could be extendable to biorders on the quandle.

The paper is organized as follows. Section \ref{definitions-known-results} recalls some basic definitions, examples and results that are required in subsequent sections. In sections \ref{knots_up_to_five_crossings}, \ref{six_crossing_knots}, \ref{seven_crossing_knots} and \ref{eight_crossing_knots}, we shall explore biorderability of knot quandles of prime knots up to eight crossings by investigating all possibilities. For each knot, we write the data in the form of a table.

\medskip

\section{Preliminaries}\label{definitions-known-results}

We begin by defining the main object of our interest.

\begin{definition}
A {\it quandle} is a non-empty set $Q$ together with a binary operation $*$ satisfying the following axioms:
\begin{enumerate}
\item[{\textbf Q1}] $x*x=x$\, for all $x\in Q$.
\item[{\textbf Q2}] For each $x, y \in Q$, there exists a unique $z \in Q$ such that $x=z*y$.
\item[{\textbf Q3}] $(x*y)*z=(x*z)*(y*z)$\, for all $x,y,z\in Q$.
\end{enumerate}
\end{definition}

The axiom {\textbf Q2} is equivalent to bijectivity of the right multiplication by each element of $Q$. This gives a dual binary operation $*^{-1}$ on $Q$ defined as $x*^{-1}y=z$ if $x=z*y$. Thus, the axiom {\textbf Q2} is equivalent to saying that
\begin{equation*}
(x*y)*^{-1}y=x\qquad\textrm{and}\qquad\left(x*^{-1}y\right)*y=x
\end{equation*}
for all $x,y\in Q$, and hence it allows us cancellation from right.
\par

Topologically, the three quandle axioms correspond to the three Reidemeister moves of planar diagrams of knots in the 3-space, which was observed independently in the foundational works of Joyce \cite{j1,j2} and Matveev \cite{m}.
\medskip

Following are examples of quandles.
\begin{itemize}
\item The binary operation $x*y=y^{-1}xy$ of conjugation in a group $G$ turns it into the quandle $\Conj(G)$ called the \textit{conjugation quandle} of $G$.
\item If $G$ is a group, then the binary operation $x*y=yx^{-1}y$ turns $G$ into the quandle $\Core(G)$ called the \textit{core quandle} of $G$.
\item If $K$ is a knot in the 3-space, then Joyce \cite{j1,j2} and Matveev \cite{m} associated a quandle $Q(K)$ to $K$ called the \textit{knot quandle} of $K$. We fix a diagram $D$ of $K$ and label its arcs. Then the knot quandle $Q(K)$ is generated by labelings of arcs of $D$ with a defining relation at each crossing in $D$ given as shown in Figure \ref{fig1}.
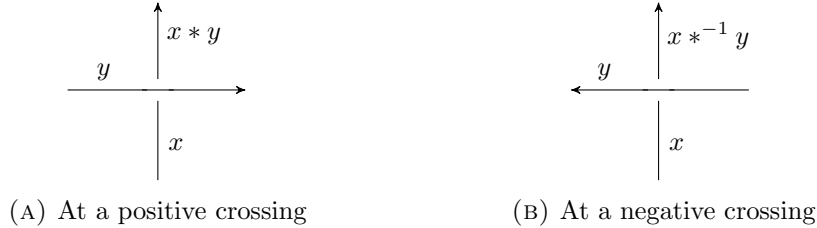
\begin{figure}[H]
\begin{subfigure}{0.4\textwidth}
\centering
\begin{tikzpicture}[scale=0.6]
\node at (0.4,-1.2) {{\small $x$}};
\node at (-1.2,0.4) {{\small $y$}};
\node at (0.8,1.2) {{\small $x*y$}};
\begin{knot}[clip width=6, clip radius=4pt]
\strand[->] (-2,0)--(2,0);
\strand[->] (0,-2)--(0,2);
\end{knot}
\end{tikzpicture}
\caption{At a positive crossing}
\end{subfigure}
\begin{subfigure}{0.4\textwidth}
\centering
\begin{tikzpicture}[scale=0.6]
\node at (1.1,1.2) {{\small $x*^{-1}y$}};
\node at (-1.2,0.4) {{\small $y$}};
\node at (0.4,-1.2) {{\small $x$}};
\begin{knot}[clip width=6, clip radius=4pt]
\strand[->] (2,0)--(-2,0);
\strand[->] (0,-2)--(0,2);
\end{knot}
\end{tikzpicture}
\caption{At a negative crossing}
\end{subfigure}
\caption{Relations}
\label{fig1}
\end{figure}
\end{itemize}

A \textit{homomorphism of quandles} $P$ and $Q$ is a map $\phi:P\to Q$ with $\phi(x*y)=\phi(x)*\phi(y)$ for all $x,y\in P$. By cancellation in $P$ and $Q$, we obtain $\phi(x*^{-1}y)=\phi(x)*^{-1}\phi(y)$ for all $x,y\in P$.
\medskip

Let us define the following term that we shall use throughout the paper.

\begin{definition}
A knot $K$ is said to be {\it positive} if there exists a diagram $D$ of $K$ such that all its crossings are positive.
\end{definition}

A negative knot is defined similarly.
\medskip

We state the Generalized Kauffman–Harary (GKH) Conjecture below (see \cite[Conjecture 3.1]{aps}).

\begin{conjecture}
Let $K$ be a prime alternating knot, and $D$ be an alternating diagram of $K$ without nugatory crossings (i.e. $D$ is a minimal diagram). Then different arcs of $D$ represent different elements of $H_1(X^{(2)}_K,\mathbb{Z})$, where $X^{(2)}_K$ denotes the double cover of $S^3$ branched along $K$.
\end{conjecture}

The authors in \cite{aps} proved the following (see \cite[Theorem 4.2]{aps}).

\begin{theorem}\label{GKH_theorem}
The GKH conjecture holds for all alternating Montesinos knots.
\end{theorem}

\begin{remark}\label{GKH_remark_1}
Let $K$ be a non-trivial Montesinos knot that is prime and alternating. Consider an alternating diagram $D$ of $K$ without a nugatory crossing (i.e. $D$ is a minimal diagram). Let $x_0, x_1,\ldots, x_n$ be generators of the knot quandle $Q(K)$ corresponding to arcs in $D$. By Theorem \ref{GKH_theorem}, there is a quandle homomorphism $\phi:Q(K)\to H_1(X^{(2)}_K,\mathbb{Z})$ such that the generators $x_0, x_1,\ldots, x_n$ map to different elements of $H_1(X^{(2)}_K,\mathbb{Z})$. Hence, the elements $x_0, x_1,\ldots, x_n$ in $Q(K)$ are all pairwise distinct.
\end{remark}

\begin{remark}\label{GKH_remark_2}
Let $K$ be a prime non-trivial knot of eight or fewer crossings, excluding the knots $8_{12}$, $8_{14}$, $8_{16}$, $8_{17}$, $8_{18}$, $8_{19}$, $8_{20}$ and $8_{21}$. By \cite[Theorem 18]{dm}, the knot $K$ is a pretzel and hence a Montesinos one. Also, the knot $K$ is alternating. Thus, Remark \ref{GKH_remark_1} holds true in this case.
\end{remark}

We now define orderability of quandles.

\begin{definition}
A quandle $Q$ is said to be \textit{left-orderable} if there is a (strict) linear order $<$ on $Q$ such that $x<y$ implies $z*x<z*y$ for all $x,y,z\in Q$. Similarly, a quandle $Q$ is \textit{right-orderable} if there is a linear order $<^\prime$ on $Q$ such that $x<^\prime y$ implies $x*z<^\prime y*z$ for all $x,y,z\in Q$. A quandle is \textit{bi-orderable} if it has a linear order with respect to which it is both left and right ordered.
\end{definition}

For example, a trivial quandle (a quandle $Q$ such that $x*y=x$ for all $x,y\in Q$) with more than one element is right-orderable but not left-orderable. Notice the contrast to groups where left-orderability implies right-orderability and vice-versa.

In \cite{rss}, left, right and bi-oderability is discussed for many classes of knots and links. For example, biorderability of positive and negative links, and right orderability of torus links is investigated in the same paper. Further, in \cite{drs}, orderability of Dehn quandles of groups is discussed.

\medskip

\section{Knots up to Five Crossings}\label{knots_up_to_five_crossings}

By \cite[Corollary 7.4]{rss}, the knot quandle of the trefoil knot (i.e. of the knot $3_1$) is neither left nor right orderable, and thus it is not biorderable. Further, by \cite[Corollary 6.7]{rss}, knot quandles of knots $5_1$ and $5_2$ are not biorderable.

\begin{figure}[H]
	\centering
	\includegraphics[scale=0.3]{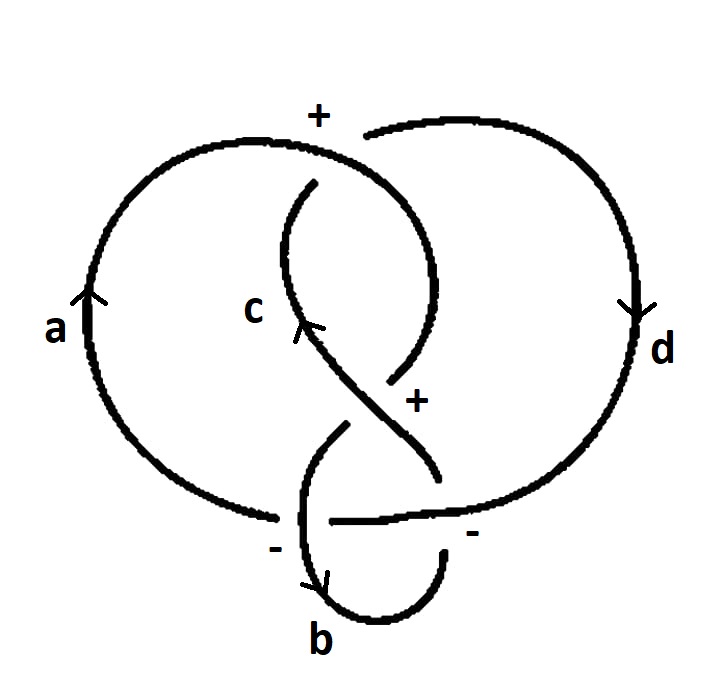}
	\caption{Knot $4_1$}
	\label{fig4_1}
\end{figure}

Since all the roots of the Alexander polynomial of the figure eight knot (i.e. of the knot $4_1$, which is a fibred knot) are real and positive, its knot quandle is right orderable \cite[Corollary 6.4]{rss}. Suppose there is a biordering $<$ on $Q(4_1)$. Let $a_1\blackdiamond a_2\blackdiamond\cdots \blackdiamond a_k$ denote either $a_1< a_2<\cdots< a_k$ or $a_1> a_2>\cdots > a_k$. Now, consider the diagram of the knot $4_1$ as in Figure \ref{fig4_1}. The quandle $Q(4_1)$ is generated by labelings of arcs of the diagram, and the relations are given by (see \cite[Proposition 3.9]{rss} for the relations on the right):
\begin{align}
	c*a&=d&c\blackdiamond d\blackdiamond a\label{4_1_1} \\
	a*c&=b&a\blackdiamond b\blackdiamond c\label{4_1_2} \\
	a*b&=d&a\blackdiamond d\blackdiamond b\label{4_1_3} \\
	c*d&=b&c\blackdiamond b\blackdiamond d\label{4_1_4}
\end{align}

Applying the relations above ((\ref{4_1_1}) to (\ref{4_1_4})), we have the table as follows.

\begin{table}[H]
	\caption{Orders that could be extendable to biorders on $Q(4_1)$}
	\label{tbl4_1}
	\begin{center}
	\(\begin{array}{|c|c|l|l|l|}\hline
	\text{Sr.} & \text{Label} & \text{Expression} & \text{Obtained by} & \text{{\color{red} Contradiction to} $\backslash$}\\
	\text{No.} &  & &  & \text{\color{blue} Compatible with}\\[2mm]\hline
	1 & \text{A} & a\blackdiamond d\blackdiamond b & \ref{4_1_3} & \\
	2 & \text{B} & a\blackdiamond d\blackdiamond b\blackdiamond c & \text{A},\;\ref{4_1_2} & {\color{blue} \ref{4_1_1},\;\ref{4_1_4}}\\

 \hline
\end{array}\)
\end{center}
\end{table}

Linear order on the generating set of $Q(4_1)$ as in \text{B} (in Table \ref{tbl4_1}) could be extendable to a biorder on $Q(4_1)$.

\medskip

\section{Six Crossing Knots}\label{six_crossing_knots}

In this section, we consider prime knots of six crossings with their diagrams as in the Rolfsen knot table (see \cite[Appendix A]{l}). The labelings of the arcs of a diagram are labeled as $a$, $b$, $c$, $d$, $e$ and $f$ (see figures \ref{fig6_1}, \ref{fig6_2} and \ref{fig6_3}) which are the generators for the knot quandle of the corresponding six crossing prime knot. Since all six crossing prime knots are alternating and Montesinos (see \cite[Theorem 18]{dm}), by Remark \ref{GKH_remark_2}, the generators $a$, $b$, $c$, $d$, $e$ and $f$ are all pairwise distinct.

\subsection{Knot $6_1$}\label{6_1}

\begin{figure}[H]
\centering
\includegraphics[scale=0.3]{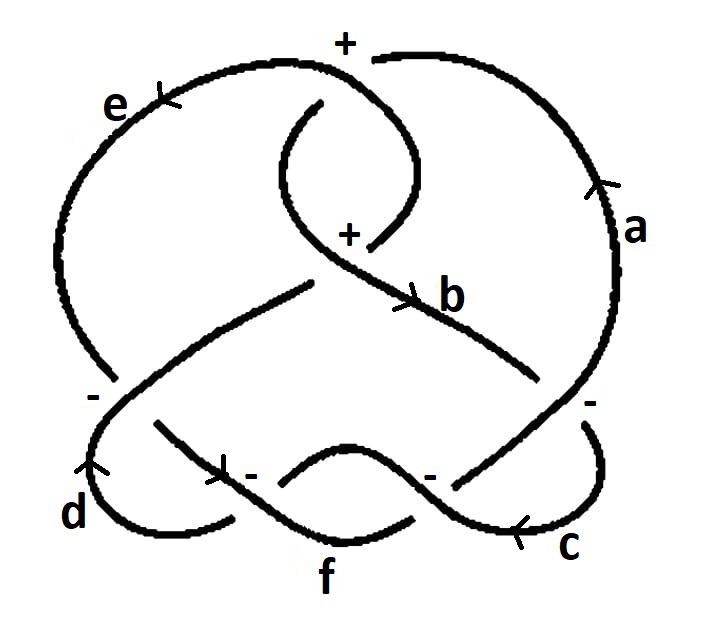}
\caption{Knot $6_1$}
\label{fig6_1}
\end{figure}

Consider the diagram of the knot $6_1$ as in Figure \ref{fig6_1}. The quandle $Q(6_1)$ is generated by labelings of arcs of the diagram, and the relations are given by

\begin{align}
a*e&=b&a\blackdiamond b\blackdiamond e\label{6_1_1} \\
d*b&=e&d\blackdiamond e\blackdiamond b\label{6_1_2} \\
f*d&=e&f\blackdiamond e\blackdiamond d\label{6_1_3} \\
c*a&=b&c\blackdiamond b\blackdiamond a\label{6_1_4} \\
d*f&=c&d\blackdiamond c\blackdiamond f\label{6_1_5} \\
a*c&=f&a\blackdiamond f\blackdiamond c \label{6_1_6}
\end{align}

Applying the relations above ((\ref{6_1_1}) to (\ref{6_1_6})), we have the table as follows.

\begin{longtblr}[caption = {Orders that could be extendable to biorders on $Q(6_1)$}, label = {tbl6_1}]
{colspec = {|c|c|l|l|l|}, rowhead = 2}
\hline
\text{Sr.} & \text{Label} & \text{Expression} & \text{Obtained by} & \text{{\color{red} Contradiction to} $\backslash$}\\
\text{No.} &  &  &  & \text{\color{blue} Compatible with}\\[2mm]\hline

1 & A & $a \blackdiamond b \blackdiamond e$ & \ref{6_1_1}  & \\
2 & B & $a \blackdiamond b \blackdiamond e \blackdiamond d$ & A,\;\ref{6_1_2} & \\
3 & \text{$C_1$} & $f \blackdiamond a \blackdiamond b \blackdiamond e \blackdiamond d$ & B,\;\ref{6_1_3} & \\
4 & \text{$C_2$} & $a \blackdiamond f \blackdiamond b \blackdiamond e \blackdiamond d$ & B,\;\ref{6_1_3} & \\
5 & \text{$C_3$} & $a \blackdiamond b \blackdiamond f \blackdiamond e \blackdiamond d$ & B,\;\ref{6_1_3} & \\

6 & \text{$D_1$} & $f \blackdiamond a \blackdiamond b \blackdiamond c \blackdiamond e \blackdiamond d$ & \text{$C_1$},\;\ref{6_1_4} & {\color{blue} \ref{6_1_5},}\;{\color{red}\ref{6_1_6}}\\
7 & \text{$D_2$} & $f \blackdiamond a \blackdiamond b \blackdiamond e \blackdiamond c \blackdiamond d$ & \text{$C_1$},\;\ref{6_1_4} & {\color{blue} \ref{6_1_5},}\;{\color{red}\ref{6_1_6}}\\
8 & \text{$D_3$} & $f \blackdiamond a \blackdiamond b \blackdiamond e \blackdiamond d \blackdiamond c$ & \text{$C_1$},\;\ref{6_1_4} & {\color{red} \ref{6_1_5},\;\ref{6_1_6}}\\
9 & \text{$E_1$} & $a \blackdiamond f \blackdiamond b \blackdiamond c \blackdiamond e \blackdiamond d$ & \text{$C_2$},\;\ref{6_1_4} & {\color{blue} \ref{6_1_5},\;\ref{6_1_6}} \\
10 & \text{$E_2$} & $a \blackdiamond f \blackdiamond b \blackdiamond e \blackdiamond c \blackdiamond d$ & \text{$C_2$},\;\ref{6_1_4} & {\color{blue} \ref{6_1_5},\;\ref{6_1_6}} \\
11 & \text{$E_3$} & $a \blackdiamond f \blackdiamond b \blackdiamond e \blackdiamond d \blackdiamond c$ & \text{$C_2$},\;\ref{6_1_4} & {\color{red} \ref{6_1_5},}\;{\color{blue}\ref{6_1_6}} \\
12 & \text{$F_1$} & $a \blackdiamond b \blackdiamond c \blackdiamond f \blackdiamond e \blackdiamond d$ & \text{$C_3$},\;\ref{6_1_4} & {\color{red} \ref{6_1_5},\;\ref{6_1_6}} \\
13 & \text{$F_2$} & $a \blackdiamond b \blackdiamond f \blackdiamond c \blackdiamond e \blackdiamond d$ & \text{$C_3$},\;\ref{6_1_4} & {\color{blue} \ref{6_1_5},\;\ref{6_1_6}} \\
14 & \text{$F_3$} & $a \blackdiamond b \blackdiamond f \blackdiamond e \blackdiamond c \blackdiamond d$ & \text{$C_3$},\;\ref{6_1_4} & {\color{blue} \ref{6_1_5},\;\ref{6_1_6}} \\
15 & \text{$F_4$} & $a \blackdiamond b \blackdiamond f \blackdiamond e \blackdiamond d \blackdiamond c$ & \text{$C_3$},\;\ref{6_1_4} & {\color{red} \ref{6_1_5},}\;{\color{blue}\ref{6_1_6}} \\
\hline
\end{longtblr}

In Table \ref{tbl6_1}, except \text{$E_1$}, \text{$E_2$}, \text{$F_2$} and \text{$F_3$}, all other possibilities contradict either (\ref{6_1_5}) or (\ref{6_1_6}), or both. Linear orders on the generating set of $Q(6_1)$ as in \text{$E_1$}, \text{$E_2$}, \text{$F_2$} and \text{$F_3$} could be extendable to biorders on $Q(6_1)$.

\subsection{Knot $6_2$}\label{6_2}

\begin{figure}[H]
\centering
\includegraphics[scale=0.3]{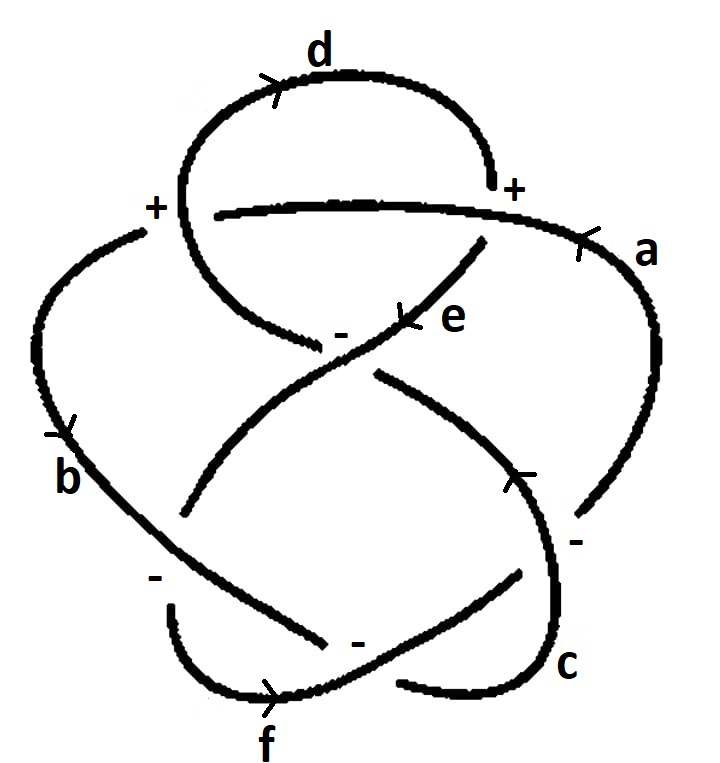}
\caption{Knot $6_2$}
\label{fig6_2}
\end{figure}

Consider the diagram of the knot $6_2$ as in Figure \ref{fig6_2}. The quandle $Q(6_2)$ is generated by labelings of arcs of the diagram, and the relations are given by

\begin{align}
a*d&=b&a\blackdiamond b\blackdiamond d\label{6_2_1}\\
f*b&=e&f\blackdiamond e\blackdiamond b\label{6_2_2}\\
c*f&=b&c\blackdiamond b\blackdiamond f\label{6_2_3}\\
a*c&=f&a\blackdiamond f\blackdiamond c\label{6_2_4}\\
d*e&=c&d\blackdiamond c\blackdiamond e\label{6_2_5}\\
d*a&=e&d\blackdiamond e\blackdiamond a\label{6_2_6}
\end{align}

We have the table as follows.

\begin{longtblr}[caption = {Orders that could be extendable to biorders on $Q(6_2)$}, label = {tbl6_2}]
{colspec = {|c|c|l|l|l|}, rowhead = 2}
	\hline
    \text{Sr.} & \text{Label} & \text{Expression} & \text{Obtained by} & \text{{\color{red} Contradiction to} $\backslash$}\\
	\text{No.} &  &  &  & \text{\color{blue} Compatible with}\\[2mm]\hline
	1 & A & $a \blackdiamond b \blackdiamond d$ & \ref{6_2_1} & \\
	2 & \text{$B_1$} & $a \blackdiamond e \blackdiamond b \blackdiamond d$ & A,\;\ref{6_2_6}& \\
	3 & \text{$B_2$} & $a \blackdiamond b \blackdiamond e \blackdiamond d$ & A,\;\ref{6_2_6}& \\
	4 & \text{$C_1$} & $a \blackdiamond f \blackdiamond e \blackdiamond b \blackdiamond d$ & \text{$B_1$},\;\ref{6_2_2} & \\
	5 & \text{$C_2$} & $f \blackdiamond a \blackdiamond e \blackdiamond b \blackdiamond d$ & \text{$B_1$},\;\ref{6_2_2} & \\
	6 & D & $a \blackdiamond b \blackdiamond e \blackdiamond c \blackdiamond d$ & \text{$B_2$},\;\ref{6_2_5} & \\
	7 & \text{$E_1$} & $a \blackdiamond f \blackdiamond e \blackdiamond b \blackdiamond c \blackdiamond d$ & \text{$C_1$},\;\ref{6_2_3} & {\color{blue} \ref{6_2_4},\;\ref{6_2_5}} \\
	8 & \text{$E_2$} & $a \blackdiamond f \blackdiamond e \blackdiamond b \blackdiamond d \blackdiamond c$ & \text{$C_1$},\;\ref{6_2_3} & {\color{blue} \ref{6_2_4},}\;{\color{red} \ref{6_2_5}} \\
	9 & \text{$F_1$} & $f \blackdiamond a \blackdiamond e \blackdiamond b \blackdiamond c \blackdiamond d$ & \text{$C_2$},\;\ref{6_2_3} & {\color{red} \ref{6_2_4},}\;{\color{blue}\ref{6_2_5}} \\
	10 & \text{$F_2$} & $f \blackdiamond a \blackdiamond e \blackdiamond b \blackdiamond d \blackdiamond c$ &\text{$C_2$},\;\ref{6_2_3} & {\color{red} \ref{6_2_4},\;\ref{6_2_5}} \\
	11 & \text{$G_1$} & $a \blackdiamond b \blackdiamond e \blackdiamond f \blackdiamond c \blackdiamond d$ & \text{$D$},\;\ref{6_2_2} & {\color{red} \ref{6_2_3},}\;{\color{blue} \ref{6_2_4}}\\
	12 & \text{$G_2$} & $a \blackdiamond b \blackdiamond e \blackdiamond c \blackdiamond f \blackdiamond d$ & \text{$D$},\;\ref{6_2_2} & {\color{red} \ref{6_2_3},\;\ref{6_2_4}} \\
	13 & \text{$G_3$} & $a \blackdiamond b \blackdiamond e \blackdiamond c \blackdiamond d \blackdiamond f$ & \text{$D$},\;\ref{6_2_2} & {\color{red} \ref{6_2_3},\;\ref{6_2_4}} \\
    \hline
	\end{longtblr}

Linear order on the generating set of $Q(6_2)$ as in $E_1$ (in Table \ref{tbl6_2}) could be extendable to a biorder on $Q(6_2)$.

\subsection{Knot $6_3$}\label{6_3}

\begin{figure}[H]
\includegraphics[scale=0.3]{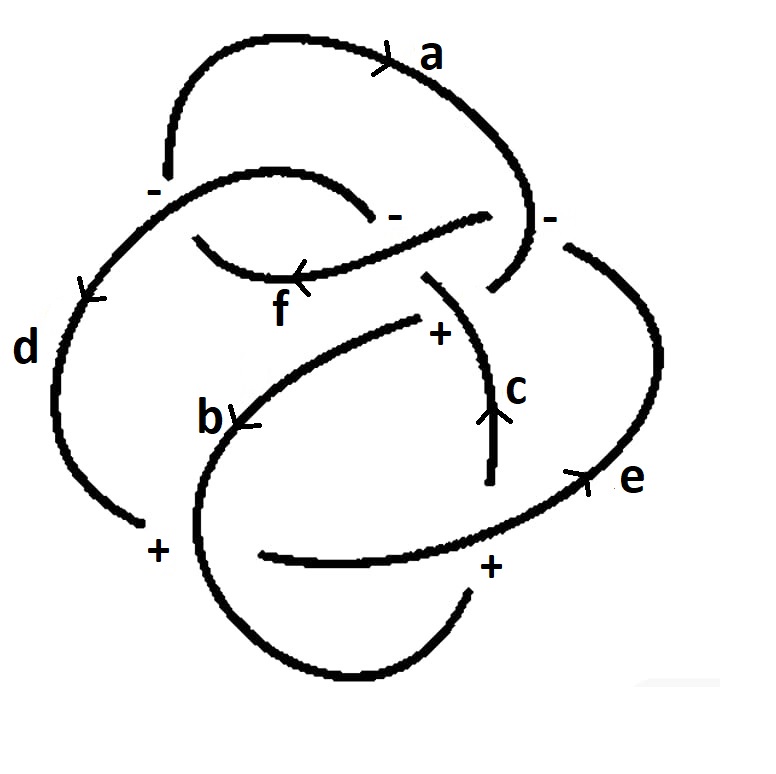}
\caption{Knot $6_3$}
\label{fig6_3}
\end{figure}

Consider the diagram of the knot $6_3$ as in Figure \ref{fig6_3}. The quandle $Q(6_3)$ is generated by labelings of arcs of the diagram, and the relations are given by

\begin{align}
a*d&=f&a\blackdiamond f\blackdiamond d\label{6_3_1}\\
f*a&=e&f\blackdiamond e\blackdiamond a\label{6_3_2}\\
a*c&=b&a\blackdiamond b\blackdiamond c\label{6_3_3}\\
d*b&=e&d\blackdiamond e\blackdiamond b\label{6_3_4}\\
b*e&=c&b\blackdiamond c\blackdiamond e\label{6_3_5}\\
d*f&=c&d\blackdiamond c\blackdiamond f\label{6_3_6}
\end{align}

We have the table as follows.

\begin{longtblr}[caption = {Orders that could be extendable to biorders on $Q(6_3)$}, label = {tbl6_3}]
{colspec = {|c|c|l|l|l|}, rowhead = 2}
		     \hline
            \text{Sr.} & \text{Label} & \text{Expression} & \text{Obtained by} & \text{{\color{red} Contradiction to} $\backslash$}\\
			\text{No.} &  &  &  & \text{\color{blue} Compatible with}\\[2mm]\hline
			1 & \text{$A$} & $a \blackdiamond f \blackdiamond d$ & \ref{6_3_1} & \\
			2 & \text{$B$} & $a \blackdiamond e \blackdiamond f \blackdiamond c \blackdiamond d$ & A,\;\ref{6_3_2},\;\ref{6_3_6} & \\
			3 & \text{$C_1$} & $a \blackdiamond b \blackdiamond e \blackdiamond f \blackdiamond c \blackdiamond d$ & \text{$B$},\;\ref{6_3_4} & {\color{blue} \ref{6_3_3},}\;{\color{red} \ref{6_3_5}} \\
			4 & \text{$C_2$} & $b \blackdiamond a \blackdiamond e \blackdiamond f \blackdiamond c \blackdiamond d$ & \text{$B$},\;\ref{6_3_4} & {\color{red} \ref{6_3_3},}\;{\color{red}\ref{6_3_5}}\\
			\hline
\end{longtblr}

Since the possibilities $C_1$ and $C_2$ in Table \ref{tbl6_3} contradict either (\ref{6_3_3}) or (\ref{6_3_5}), or both, we can say that $Q(6_3)$ can not be biorderable.

\medskip

\section{Seven Crossing Knots}\label{seven_crossing_knots}

The knots $7_1$, $7_2$, $7_3$, $7_4$ and $7_5$ are prime, alternating and positive (or negative). They are all pretzel and hence Montesinos one, see \cite[Theorem 18]{dm}. By \cite[Corollary 6.7]{rss}, the quandles of these knots are not biorderable. In fact, by \cite[Theorem 7.2]{rss}, the quandle of $7_1$ is not right-orderable.

Now, we consider knots $7_6$ and $7_7$ with their diagrams as in the Rolfsen knot table (see \cite[Appendix A]{l}). The labelings of the arcs of a diagram are labeled as $a$, $b$, $c$, $d$, $e$, $f$ and $g$ (see figures \ref{fig7_6} and \ref{fig7_7}) which are the generators for the knot quandle of the corresponding knot. Since knots $7_6$ and $7_7$ are alternating and Montesinos (see \cite[Theorem 18]{dm}), by Remark \ref{GKH_remark_2}, the generators $a$, $b$, $c$, $d$, $e$, $f$ and $g$ are all pairwise distinct.

\subsection{Knot $7_6$}\label{7_6}

\begin{figure}[H]
\centering
\includegraphics[scale=0.3]{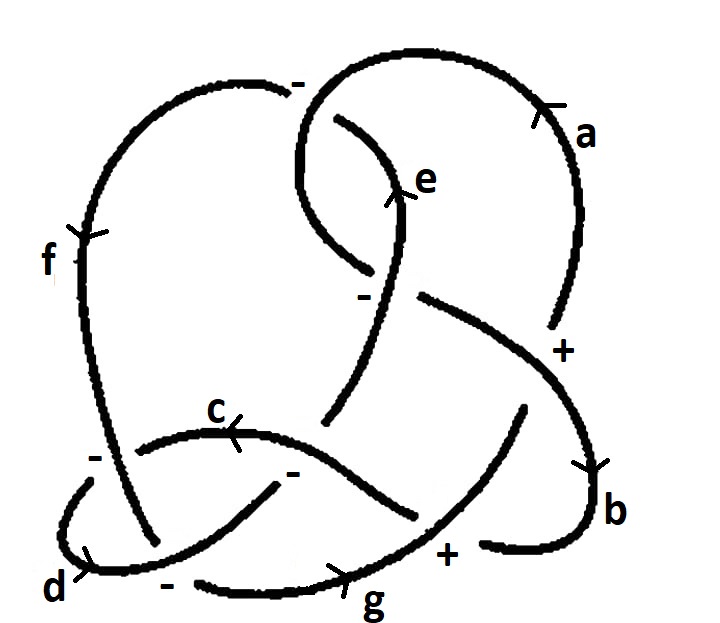}
\caption{Knot $7_6$}
\label{fig7_6}
\end{figure}

 Consider the diagram of the knot $7_6$ as in Figure \ref{fig7_6}. The quandle $Q(7_6)$ is generated by labelings of arcs of the diagram, and the relations are given by

\begin{align}
b*e&=a&b\blackdiamond a\blackdiamond e \label{7_6_1} \\
g*b&=a&g\blackdiamond a\blackdiamond b \label{7_6_2} \\
b*g&=c&b\blackdiamond c\blackdiamond g \label{7_6_3} \\
e*c&=d&e\blackdiamond d\blackdiamond c \label{7_6_4} \\
d*f&=c&d\blackdiamond c\blackdiamond f \label{7_6_5} \\
g*d&=f&g\blackdiamond f\blackdiamond d \label{7_6_6} \\
f*a&=e&f\blackdiamond e\blackdiamond a \label{7_6_7}
\end{align}

We have the table as follows.

\begin{longtblr}[caption = {Orders that could be extendable to biorders on $Q(7_6)$}, label = {tbl7_6}]
{colspec = {|c|c|l|l|l|}, rowhead = 2}
                    \hline
    				\text{Sr.} & \text{Label} & \text{Expression} & \text{Obtained by} & \text{{\color{red} Contradiction to} $\backslash$}\\
    				\text{No.} &  &  &  & \text{\color{blue} Compatible with}\\[2mm]
    				\hline
    1 & \text{A} & $b\blackdiamond a\blackdiamond e\blackdiamond f$ & \ref{7_6_1},\;\ref{7_6_7} & \\
    2 & \text{$B_1$} &  $b\blackdiamond a\blackdiamond g\blackdiamond e\blackdiamond f\blackdiamond d$ & \text{A},\;\ref{7_6_2},\;\ref{7_6_6} & \\
    3 & \text{$B_2$} &  $b\blackdiamond a\blackdiamond e\blackdiamond g\blackdiamond f\blackdiamond d$ & \text{A},\;\ref{7_6_2},\;\ref{7_6_6} & \\
    4 & \text{$B_3$} &  $b\blackdiamond a\blackdiamond e\blackdiamond f\blackdiamond g$ & \text{A},\;\ref{7_6_2} & \\
    5 & \text{$C_1$} & $b\blackdiamond c\blackdiamond a\blackdiamond g\blackdiamond e\blackdiamond f\blackdiamond d$ & \text{$B_1$},\;\ref{7_6_3} &{\color{red} \ref{7_6_4},}\;{\color{red} \ref{7_6_5}} \\
    6 &  \text{$C_2$} & $b\blackdiamond a\blackdiamond c\blackdiamond g\blackdiamond e\blackdiamond f\blackdiamond d$ & \text{$B_1$},\;\ref{7_6_3} &{\color{red} \ref{7_6_4},}\;{\color{red} \ref{7_6_5}} \\
    7 &  \text{$D_1$} & $b\blackdiamond c\blackdiamond a\blackdiamond e\blackdiamond g\blackdiamond f\blackdiamond d$ & \text{$B_2$},\;\ref{7_6_3} &{\color{red} \ref{7_6_4},}\;{\color{red} \ref{7_6_5}} \\
    8 & \text{$D_2$} & $b\blackdiamond a\blackdiamond c\blackdiamond e\blackdiamond g\blackdiamond f\blackdiamond d$ & \text{$B_2$},\;\ref{7_6_3} &{\color{red} \ref{7_6_4},}\;{\color{red} \ref{7_6_5}} \\
    9 & \text{$D_3$} & $b\blackdiamond a\blackdiamond e\blackdiamond c\blackdiamond g\blackdiamond f\blackdiamond d$ & \text{$B_2$},\;\ref{7_6_3} &{\color{red} \ref{7_6_4},}\;{\color{red} \ref{7_6_5}} \\
    10 & \text{$E_1$} & $b\blackdiamond c\blackdiamond a\blackdiamond e\blackdiamond f\blackdiamond g$ & \text{$B_3$},\;\ref{7_6_3} & \\
    11 & \text{$E_2$} & $b\blackdiamond a\blackdiamond c\blackdiamond d\blackdiamond e\blackdiamond f\blackdiamond g$ & \text{$B_3$},\;\ref{7_6_3},\;\ref{7_6_4} & {\color{red} \ref{7_6_5},}\;{\color{blue} \ref{7_6_6}} \\
    12 & \text{$E_3$} & $b\blackdiamond a\blackdiamond e\blackdiamond d\blackdiamond c\blackdiamond f\blackdiamond g$ & \text{$B_3$},\;\ref{7_6_3},\;\ref{7_6_4} & {\color{blue} \ref{7_6_5},}\;{\color{blue} \ref{7_6_6}} \\
    13 & \text{$E_4$} & $b\blackdiamond a\blackdiamond e\blackdiamond f\blackdiamond c\blackdiamond g$ & \text{$B_3$},\;\ref{7_6_3} & \\
    14 &  \text{$F_1$} & $b\blackdiamond d\blackdiamond c\blackdiamond a\blackdiamond e\blackdiamond f\blackdiamond g$ & \text{$E_1$},\;\ref{7_6_5} & {\color{red} \ref{7_6_4},}\;{\color{blue} \ref{7_6_6}} \\
    15 &  \text{$F_2$} & $d\blackdiamond b\blackdiamond c\blackdiamond a\blackdiamond e\blackdiamond f\blackdiamond g$ & \text{$E_1$},\;\ref{7_6_5} & {\color{red} \ref{7_6_4},}\;{\color{blue} \ref{7_6_6}} \\
    16 &  \text{$G_1$} & $b\blackdiamond a\blackdiamond e\blackdiamond d\blackdiamond f\blackdiamond c\blackdiamond g$ & \text{$E_4$},\;\ref{7_6_4} &{\color{red} \ref{7_6_5},}\;{\color{blue} \ref{7_6_6}} \\
    17 &  \text{$G_2$} & $b\blackdiamond a\blackdiamond e\blackdiamond f\blackdiamond d\blackdiamond c\blackdiamond g$ & \text{$E_4$},\;\ref{7_6_4} &{\color{red} \ref{7_6_5},}\;{\color{red} \ref{7_6_6}} \\

    \hline
\end{longtblr}

The linear order on the generating set of $Q(7_6)$ as in \text{$E_3$} (in Table \ref{tbl7_6}) could be extendable to a biorder on $Q(7_6)$.

\subsection{Knot $7_7$}\label{7_7}

\begin{figure}[H]
    	\centering
    	\includegraphics[scale=0.34]{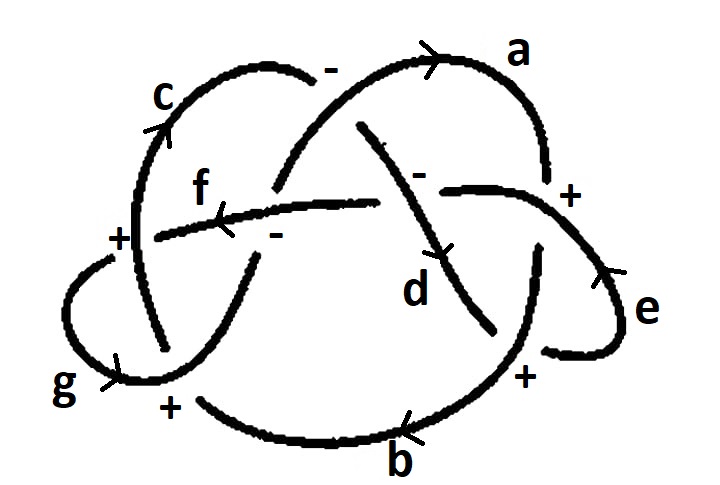}
    	\caption{Knot $7_7$}
        \label{fig7_7}
    \end{figure}

Consider the diagram of the knot $7_7$ as in Figure \ref{fig7_7}. The quandle $Q(7_7)$ is generated by labelings of arcs of the diagram, and the relations are given by

    \begin{align}
    	a*e&=b&a\blackdiamond b\blackdiamond e \label{7_7_1} \\
    	d*b&=e&d\blackdiamond e\blackdiamond b \label{7_7_2} \\
    	b*g&=c&b\blackdiamond c\blackdiamond g \label{7_7_3} \\
    	f*c&=g&f\blackdiamond g\blackdiamond c \label{7_7_4} \\
    	d*a&=c&d\blackdiamond c\blackdiamond a \label{7_7_5} \\
    	a*f&=g&a\blackdiamond g\blackdiamond f \label{7_7_6} \\
    	f*d&=e&f\blackdiamond e\blackdiamond d \label{7_7_7}
    \end{align}

We have the table as follows.

\begin{longtblr}[caption = {Orders that could be extendable to biorders on $Q(7_7)$}, label = {tbl7_7}]
{colspec = {|c|c|l|l|l|}, rowhead = 2}
	 \hline
   	\text{Sr.} & \text{Label} & \text{Expression} & \text{Obtained by} & \text{{\color{red} Contradiction to} $\backslash$}\\
 	   \text{No.} &  &  &  & \text{\color{blue} Compatible with}\\[2mm]
    \hline
    1 & \text{A} & $b\blackdiamond c\blackdiamond g\blackdiamond f$ & \ref{7_7_3},\;\ref{7_7_4} & \\
    2 & \text{$B_1$} & $d\blackdiamond e\blackdiamond b\blackdiamond c\blackdiamond a\blackdiamond g\blackdiamond f$ & \text{A},\;\ref{7_7_6},\;\ref{7_7_1},\;\ref{7_7_7} &{\color{blue} \ref{7_7_2},}\;{\color{blue} \ref{7_7_5}} \\
    3 &  \text{$B_2$} & $d\blackdiamond e\blackdiamond b\blackdiamond a\blackdiamond c\blackdiamond g\blackdiamond f$ & \text{A},\;\ref{7_7_6},\;\ref{7_7_1},\;\ref{7_7_7} &{\color{blue} \ref{7_7_2},}\;{\color{red} \ref{7_7_5}} \\
    4 &  \text{$B_3$} & $a\blackdiamond b\blackdiamond c\blackdiamond g\blackdiamond f$ & \text{A},\;\ref{7_7_6} & \\
   5 & \text{$C_1$} & $a\blackdiamond b\blackdiamond c\blackdiamond d\blackdiamond g\blackdiamond f$ & \text{$B_3$},\;\ref{7_7_5} & \\
   6 & \text{$C_2$} & $a\blackdiamond b\blackdiamond c\blackdiamond g\blackdiamond d\blackdiamond e\blackdiamond f$ & \text{$B_3$},\;\ref{7_7_5},\;\ref{7_7_7} & {\color{blue} \ref{7_7_1},}\;{\color{red} \ref{7_7_2}}\\
   7 &  \text{$C_3$} & $a\blackdiamond b\blackdiamond c\blackdiamond g\blackdiamond f\blackdiamond e\blackdiamond d$ & \text{$B_3$},\;\ref{7_7_5},\;\ref{7_7_7} & {\color{blue} \ref{7_7_1},}\;{\color{blue} \ref{7_7_2}}\\
   8 &   \text{$D_1$} & $a\blackdiamond b\blackdiamond c\blackdiamond d\blackdiamond e\blackdiamond g\blackdiamond f$ & \text{$C_1$},\;\ref{7_7_7} & {\color{blue} \ref{7_7_1},}\;{\color{red} \ref{7_7_2}}\\
   9 &  \text{$D_2$} & $a\blackdiamond b\blackdiamond c\blackdiamond d\blackdiamond g\blackdiamond e\blackdiamond f$ & \text{$C_1$},\;\ref{7_7_7} & {\color{blue} \ref{7_7_1},}\;{\color{red} \ref{7_7_2}}\\

   	\hline
   	\end{longtblr}

 Linear orders on the generating set of $Q(7_7)$ as in \text{$B_1$} and \text{$C_3$} (in Table \ref{tbl7_7}) could be extendable to biorders on $Q(7_7)$.

\medskip

\section{Eight Crossing Knots}\label{eight_crossing_knots}

In this section, we consider prime knots of eight crossings with their diagrams as in the Rolfsen knot table (see \cite[Appendix A]{l}). The labelings of the arcs of a diagram are labeled as $a$, $b$, $c$, $d$, $e$, $f$, $g$ and $h$ (see figures \ref{fig8_1} to \ref{fig8_21}) which are the generators for the knot quandle of the corresponding eight crossing prime knot. Since, except $8_{12}$, $8_{14}$, $8_{16}$, $8_{17}$, $8_{18}$, $8_{19}$, $8_{20}$ and $8_{21}$, all other eight crossing prime knots are alternating and Montesinos (see \cite[Theorem 18]{dm}), by Remark \ref{GKH_remark_2}, the generators $a$, $b$, $c$, $d$, $e$, $f$, $g$ and $h$ are all pairwise distinct.

\subsection{Knot $8_1$}\label{8_1}

\begin{figure}[H]
    	\centering
    	\includegraphics[scale=0.3]{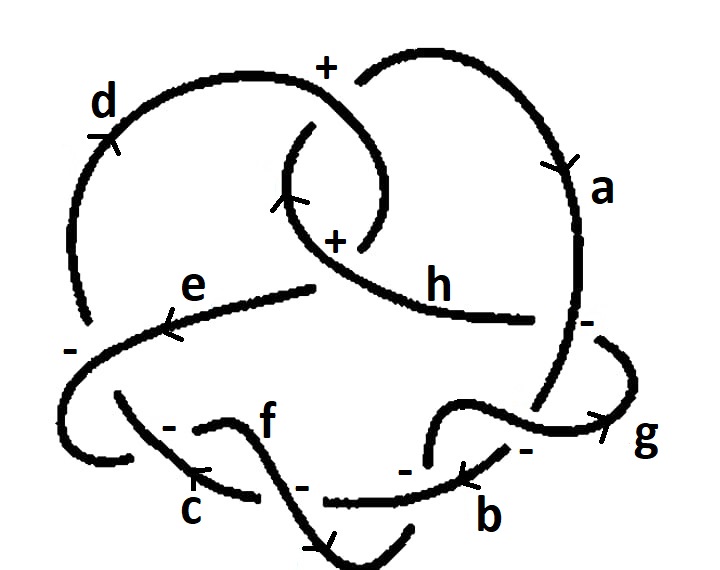}
    	\caption{Knot $8_1$}
        \label{fig8_1}
    \end{figure}

 Consider the diagram of the knot $8_1$ as in Figure \ref{fig8_1}. The quandle $Q(8_1)$ is generated by labelings of arcs of the diagram, and the relations are given by
     \begin{align}
    	b*g&=a&b\blackdiamond a\blackdiamond g  \label{8_1_1}\\
    	g*b&=f&g\blackdiamond f\blackdiamond b \label{8_1_2}\\
    	c*f&=b&c\blackdiamond b\blackdiamond f \label{8_1_3}\\
    	f*c&=e&f\blackdiamond e\blackdiamond c \label{8_1_4}\\
    	d*e&=c&d\blackdiamond c\blackdiamond e  \label{8_1_5}\\
    	h*d&=a&h\blackdiamond a\blackdiamond d \label{8_1_6}\\
    	d*h&=e&d\blackdiamond e\blackdiamond h \label{8_1_7}\\
    	h*a&=g&h\blackdiamond g\blackdiamond a \label{8_1_8}
    \end{align}

We have the table as follows.

\begin{longtblr}[caption = {Orders that could be extendable to biorders on $Q(8_1)$}, label = {tbl8_1}]
{colspec = {|c|c|l|l|l|}, rowhead = 2}
                \hline
       			\text{Sr.} & \text{Label} & \text{Expression} & \text{Obtained by} & \text{{\color{red} Contradiction to} $\backslash$}\\
       			\text{No.} &  & &  & \text{\color{blue} Compatible with}\\[2mm]
       			\hline
       			1 & \text{A} & $b\blackdiamond a\blackdiamond g$ & \ref{8_1_1} & \\
       			2 & \text{$B_1$} & $c\blackdiamond b\blackdiamond f\blackdiamond a\blackdiamond g\blackdiamond h$ & A,\; \ref{8_1_2},\; \ref{8_1_3},\;\ref{8_1_8} & \\
       			3 & \text{$B_2$} & $c\blackdiamond b\blackdiamond a\blackdiamond f\blackdiamond g\blackdiamond h$ & A,\;\ref{8_1_2},\; \ref{8_1_3},\;\ref{8_1_8} & \\
       			4 & \text{$C_1$} & $d\blackdiamond c\blackdiamond e\blackdiamond b\blackdiamond f\blackdiamond a\blackdiamond g\blackdiamond h$ & \text{$B_1$},\; \ref{8_1_4},\;\ref{8_1_5} & {\color{blue} \ref{8_1_6},\;\ref{8_1_7}} \\
       			5 & \text{$C_2$} & $d\blackdiamond c\blackdiamond b\blackdiamond e\blackdiamond f\blackdiamond a\blackdiamond g\blackdiamond h$ & \text{$B_1$},\; \ref{8_1_4},\;\ref{8_1_5} & {\color{blue} \ref{8_1_6},\;\ref{8_1_7}} \\
       			6 & \text{$D_1$} & $d\blackdiamond c\blackdiamond e\blackdiamond b\blackdiamond a\blackdiamond f\blackdiamond g\blackdiamond h$ & \text{$B_2$},\; \ref{8_1_4},\;\ref{8_1_5} & {\color{blue} \ref{8_1_6},\;\ref{8_1_7}} \\
       			7 & \text{$D_2$} & $d\blackdiamond c\blackdiamond b\blackdiamond a\blackdiamond e\blackdiamond f\blackdiamond g\blackdiamond h$ & \text{$B_2$},\; \ref{8_1_4},\;\ref{8_1_5} & {\color{blue} \ref{8_1_6},\;\ref{8_1_7}} \\
       			8 & \text{$D_3$} & $d\blackdiamond c\blackdiamond b\blackdiamond e\blackdiamond a\blackdiamond f\blackdiamond g\blackdiamond h$ & \text{$B_2$},\; \ref{8_1_4},\;\ref{8_1_5} & {\color{blue} \ref{8_1_6},\;\ref{8_1_7}} \\

       			\hline
       	\end{longtblr}

In Table \ref{tbl8_1}, linear orderings \text{$C_1$}, \text{$C_2$}, \text{$D_1$}, \text{$D_2$} and \text{$D_3$} are compatible with all quandle relations. These orders could be extendable to biorders on $Q(8_1)$.

\subsection{Knot $8_2$}\label{8_2}

\begin{figure}[H]
       	 	\centering
       	 	\includegraphics[scale=0.3]{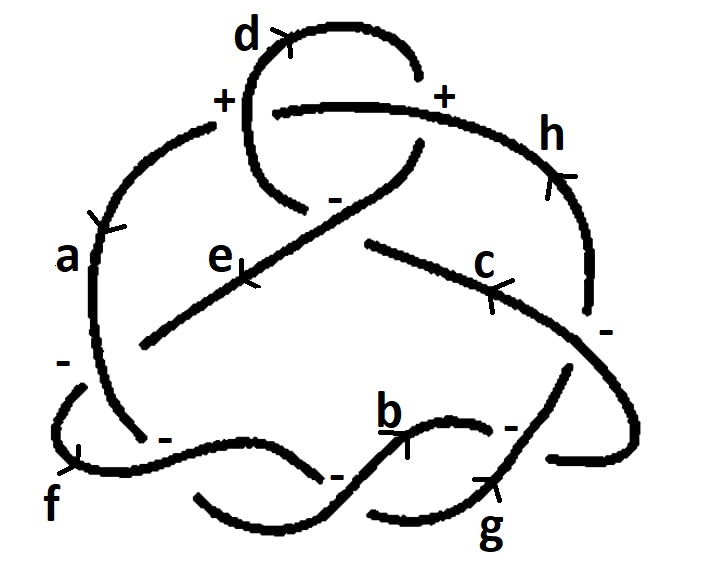}
       	 	\caption{Knot $8_2$}
            \label{fig8_2}
       	 \end{figure}

Consider the diagram of the knot $8_2$ as in Figure \ref{fig8_2}. The quandle $Q(8_2)$ is generated by labelings of arcs of the diagram, and the relations are given by

       	 \begin{align}
       	 	h*d&=a&h\blackdiamond a\blackdiamond d \label{8_2_1}\\
       	 	d*h&=e&d\blackdiamond e\blackdiamond h \label{8_2_2}\\
       	 	d*e&=c&d\blackdiamond c\blackdiamond e \label{8_2_3}\\
       	 	f*a&=e&f\blackdiamond e\blackdiamond a \label{8_2_4}\\
       	 	b*f&=a&b\blackdiamond a\blackdiamond f \label{8_2_5}\\
       	 	g*b&=f&g\blackdiamond f\blackdiamond b \label{8_2_6}\\
       	 	c*g&=b&c\blackdiamond b\blackdiamond g \label{8_2_7}\\
       	 	h*c&=g&h\blackdiamond g\blackdiamond c \label{8_2_8}
       	 \end{align}

We have the table as follows.

\begin{longtblr}[caption = {Orders that could be extendable to biorders on $Q(8_2)$}, label = {tbl8_2}]
{colspec = {|c|c|l|l|l|}, rowhead = 2}
       	    \hline
       		\text{Sr.} & \text{Label} & \text{Expression} & \text{Obtained by} & \text{{\color{red} Contradiction to} $\backslash$}\\
    		\text{No.} &  & &  & \text{\color{blue} Compatible with}\\[2mm]
   			\hline
   		1 & \text{A} & $c\blackdiamond b\blackdiamond g\blackdiamond h$ & \ref{8_2_7},\; \ref{8_2_8} & \\
        2 & \text{B} & $c\blackdiamond b\blackdiamond f\blackdiamond g\blackdiamond h$ & \text{A},\; \ref{8_2_6} & \\
        3 & \text{C} & $c\blackdiamond b\blackdiamond a\blackdiamond f\blackdiamond g\blackdiamond h$ & \text{B},\; \ref{8_2_5} & \\
        4 & \text{D} & $c\blackdiamond b\blackdiamond a\blackdiamond e\blackdiamond f\blackdiamond g\blackdiamond h$ & \text{C},\; \ref{8_2_4} & \\
        5 & \text{E} & $d\blackdiamond c\blackdiamond b\blackdiamond a\blackdiamond e\blackdiamond f\blackdiamond g\blackdiamond h$ & \text{D},\; \ref{8_2_3} & {\color{blue} \ref{8_2_1}},\;{\color{blue} \ref{8_2_2} } \\
   			\hline
   	\end{longtblr}

In Table \ref{tbl8_2}, linear ordering \text{E} is compatible with all quandle relations. This order could be extendable to a biorder on $Q(8_2)$.

   	\subsection{Knot $8_3$}\label{8_3}

\begin{figure}[H]
   		\centering
   		\includegraphics[scale=0.3]{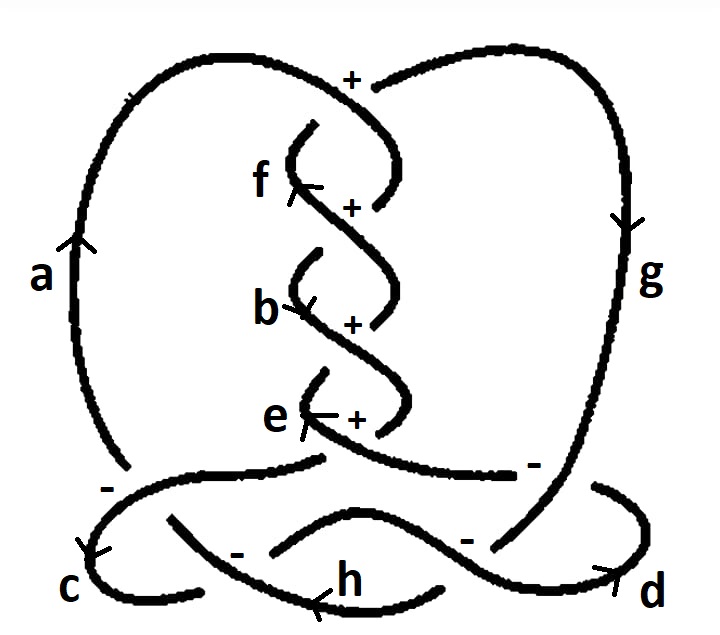}
   		\caption{Knot $8_3$}
        \label{fig8_3}
   	\end{figure}

   	Consider the diagram of the knot $8_3$ as in Figure \ref{fig8_3}. The quandle $Q(8_3)$ is generated by labelings of arcs of the diagram, and the relations are given by

   	\begin{align}
   	f*a&=g&f\blackdiamond g\blackdiamond a \label{8_3_1}\\
   	a*f&=b&a\blackdiamond b\blackdiamond f \label{8_3_2}\\
   	e*b&=f&e\blackdiamond f\blackdiamond b \label{8_3_3}\\
   	b*e&=c&b\blackdiamond c\blackdiamond e \label{8_3_4}\\
   	a*c&=h&a\blackdiamond h\blackdiamond c \label{8_3_5}\\
   	d*h&=c&d\blackdiamond c\blackdiamond h \label{8_3_6}\\
   	h*d&=g&h\blackdiamond g\blackdiamond d \label{8_3_7}\\
   	e*g&=d&e\blackdiamond d\blackdiamond g \label{8_3_8}
   	\end{align}

   	We have the table as follows.

 \begin{longtblr}[caption = {Orders that could be extendable to biorders on $Q(8_3)$}, label = {tbl8_3}]
{colspec = {|c|c|l|l|l|}, rowhead = 2}
                \hline
   				\text{Sr.} & \text{Label} & \text{Expression} & \text{Obtained by} & \text{{\color{red} Contradiction to} $\backslash$}\\
   				\text{No.} &  &  &  & \text{\color{blue} Compatible with}\\[2mm]
   				\hline
   	1 & \text{A} & $a\blackdiamond h\blackdiamond c\blackdiamond d$ & \ref{8_3_5},\;\ref{8_3_6} & \\
   	2 & \text{$B_1$} & $a\blackdiamond h\blackdiamond g\blackdiamond c\blackdiamond d\blackdiamond e$ &\text{A},\;\ref{8_3_7},\;\ref{8_3_8} & \\
   	3 & \text{$B_2$} & $a\blackdiamond h\blackdiamond c\blackdiamond g\blackdiamond d\blackdiamond e$ &\text{A},\;\ref{8_3_7},\;\ref{8_3_8} & \\
  	4 & \text{$C_1$} & $a\blackdiamond h\blackdiamond g\blackdiamond f\blackdiamond c\blackdiamond d\blackdiamond e$ &\text{$B_1$},\;\ref{8_3_1} & \\

   	5 & \text{$C_2$} & $a\blackdiamond h\blackdiamond g\blackdiamond c\blackdiamond f\blackdiamond d\blackdiamond e$ &\text{$B_1$},\;\ref{8_3_1} & \\
   	6 & \text{$C_3$} & $a\blackdiamond h\blackdiamond g\blackdiamond c\blackdiamond d\blackdiamond f\blackdiamond e$ &\text{$B_1$},\;\ref{8_3_1} & \\
   	7 & \text{$C_4$} & $a\blackdiamond h\blackdiamond g\blackdiamond c\blackdiamond d\blackdiamond e\blackdiamond f\blackdiamond b$ &\text{$B_1$},\;\ref{8_3_1},\;\ref{8_3_3} & {\color{red} \ref{8_3_2},\;\ref{8_3_4}} \\
   	8 & \text{$D_1$} & $a\blackdiamond h\blackdiamond c\blackdiamond g\blackdiamond f\blackdiamond d\blackdiamond e$ &\text{$B_2$},\;\ref{8_3_1} & \\
   	9 & \text{$D_2$} & $a\blackdiamond h\blackdiamond c\blackdiamond g\blackdiamond d\blackdiamond f\blackdiamond e$ &\text{$B_2$},\;\ref{8_3_1} & \\
   10 & \text{$D_3$} & $a\blackdiamond h\blackdiamond c\blackdiamond g\blackdiamond d\blackdiamond e\blackdiamond f\blackdiamond b$ &\text{$B_2$},\;\ref{8_3_1},\;\ref{8_3_3} &{\color{red} \ref{8_3_2},\;\ref{8_3_4}} \\
   11 & \text{$E_1$} & $a\blackdiamond b\blackdiamond h\blackdiamond g\blackdiamond f\blackdiamond c\blackdiamond d\blackdiamond e$ &\text{$C_1$},\;\ref{8_3_2} &{\color{blue} \ref{8_3_3},\;\ref{8_3_4}} \\
   12 & \text{$E_2$} & $a\blackdiamond h\blackdiamond b\blackdiamond g\blackdiamond f\blackdiamond c\blackdiamond d\blackdiamond e$ &\text{$C_1$},\;\ref{8_3_2} &{\color{blue} \ref{8_3_3},\;\ref{8_3_4}} \\
   13 & \text{$E_3$} & $a\blackdiamond h\blackdiamond g\blackdiamond b\blackdiamond f\blackdiamond c\blackdiamond d\blackdiamond e$ &\text{$C_1$},\;\ref{8_3_2} &{\color{blue} \ref{8_3_3},\;\ref{8_3_4}} \\
   14 & \text{$F_1$} & $a\blackdiamond h\blackdiamond g\blackdiamond b\blackdiamond c\blackdiamond f\blackdiamond d\blackdiamond e$ &\text{$C_2$},\;\ref{8_3_4} &{\color{blue} \ref{8_3_2},\;\ref{8_3_3}} \\
   15 & \text{$F_2$} & $a\blackdiamond h\blackdiamond b\blackdiamond g\blackdiamond c\blackdiamond f\blackdiamond d\blackdiamond e$ &\text{$C_2$},\;\ref{8_3_4} &{\color{blue} \ref{8_3_2},\;\ref{8_3_3}} \\
   16 & \text{$F_3$} & $a\blackdiamond b\blackdiamond h\blackdiamond g\blackdiamond c\blackdiamond f\blackdiamond d\blackdiamond e$ &\text{$C_2$},\;\ref{8_3_4} &{\color{blue} \ref{8_3_2},\;\ref{8_3_3}} \\
   17 & \text{$F_4$} & $b\blackdiamond a\blackdiamond h\blackdiamond g\blackdiamond c\blackdiamond f\blackdiamond d\blackdiamond e$ &\text{$C_2$},\;\ref{8_3_4} &{\color{red} \ref{8_3_2},}\;{\color{blue} \ref{8_3_3}} \\
   18 & \text{$G_1$} & $a\blackdiamond h\blackdiamond g\blackdiamond b\blackdiamond c\blackdiamond d\blackdiamond f\blackdiamond e$ &\text{$C_3$},\;\ref{8_3_4} &{\color{blue} \ref{8_3_2},\;\ref{8_3_3}} \\
   19 & \text{$G_2$} & $a\blackdiamond h\blackdiamond b\blackdiamond g\blackdiamond c\blackdiamond d\blackdiamond f\blackdiamond e$ &\text{$C_3$},\;\ref{8_3_4} &{\color{blue} \ref{8_3_2},\;\ref{8_3_3}} \\
   20 & \text{$G_3$} & $a\blackdiamond b\blackdiamond h\blackdiamond g\blackdiamond c\blackdiamond d\blackdiamond f\blackdiamond e$ &\text{$C_3$},\;\ref{8_3_4} &{\color{blue} \ref{8_3_2},\;\ref{8_3_3}} \\
   21 & \text{$G_4$} & $b\blackdiamond a\blackdiamond h\blackdiamond g\blackdiamond c\blackdiamond d\blackdiamond f\blackdiamond e$ &\text{$C_3$},\;\ref{8_3_4} &{\color{red} \ref{8_3_2},}\;{\color{blue} \ref{8_3_3}} \\
   22 & \text{$H_1$} & $a\blackdiamond h\blackdiamond b\blackdiamond c\blackdiamond g\blackdiamond f\blackdiamond d\blackdiamond e$ &\text{$D_1$},\;\ref{8_3_4} &{\color{blue} \ref{8_3_2},\;\ref{8_3_3}} \\
   23 & \text{$H_2$} & $a\blackdiamond b\blackdiamond h\blackdiamond c\blackdiamond g\blackdiamond f\blackdiamond d\blackdiamond e$ &\text{$D_1$},\;\ref{8_3_4} &{\color{blue} \ref{8_3_2},\;\ref{8_3_3}} \\
   24 & \text{$H_3$} & $b\blackdiamond a\blackdiamond h\blackdiamond c\blackdiamond g\blackdiamond f\blackdiamond d\blackdiamond e$ &\text{$D_1$},\;\ref{8_3_4} &{\color{red} \ref{8_3_2},}\;{\color{blue} \ref{8_3_3}} \\
   25 & \text{$I_1$} & $a\blackdiamond h\blackdiamond b\blackdiamond c\blackdiamond g\blackdiamond d\blackdiamond f\blackdiamond e$ &\text{$D_2$},\;\ref{8_3_4} &{\color{blue} \ref{8_3_2},\;\ref{8_3_3}} \\
   26 & \text{$I_2$} & $a\blackdiamond b\blackdiamond h\blackdiamond c\blackdiamond g\blackdiamond d\blackdiamond f\blackdiamond e$ &\text{$D_2$},\;\ref{8_3_4} &{\color{blue} \ref{8_3_2},\;\ref{8_3_3}} \\
   27 & \text{$I_3$} & $b\blackdiamond a\blackdiamond h\blackdiamond c\blackdiamond g\blackdiamond d\blackdiamond f\blackdiamond e$ &\text{$D_2$},\;\ref{8_3_4} &{\color{red} \ref{8_3_2},}\;{\color{blue} \ref{8_3_3}} \\
    \hline
   	\end{longtblr}

Linear orders on the generating set of $Q(8_3)$ as in \text{$E_1$}, \text{$E_2$}, \text{$E_3$}, \text{$F_1$}, \text{$F_2$}, \text{$F_3$}, \text{$G_1$}, \text{$G_2$}, \text{$G_3$}, \text{$H_1$}, \text{$H_2$}, \text{$I_1$} and \text{$I_2$} in Table \ref{tbl8_3} could be extendable to biorders on $Q(8_3)$.

\subsection{Knot $8_4$}\label{8_4}

\begin{figure}[H]
	\centering
	\includegraphics[scale=0.3]{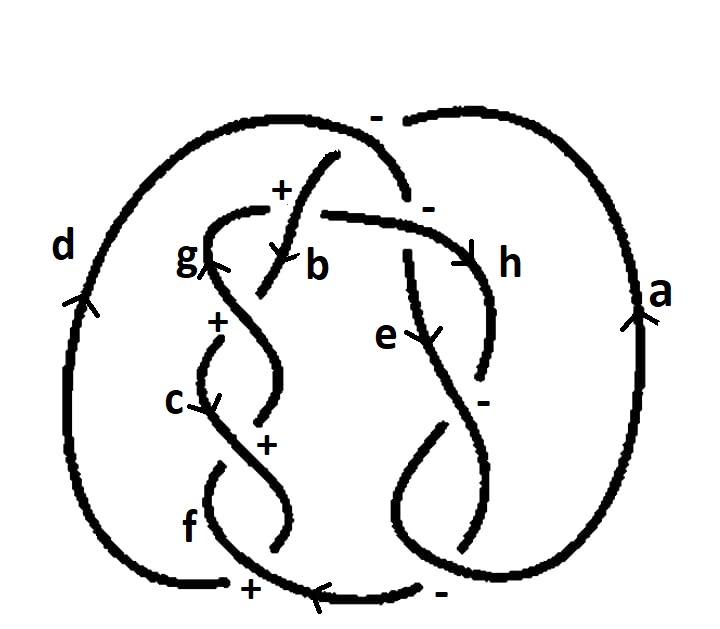}
	\caption{Knot $8_4$}
    \label{fig8_4}
\end{figure}

  Consider the diagram of the knot $8_4$ as in Figure \ref{fig8_4}. The quandle $Q(8_4)$ is generated by labelings of arcs of the diagram, and the relations are given by

  \begin{align}
  	b*d&=a&b\blackdiamond a\blackdiamond d \label{8_4_1}\\
  	g*b&=h&g\blackdiamond h\blackdiamond b \label{8_4_2}\\
  	b*g&=c&b\blackdiamond c\blackdiamond g \label{8_4_3}\\
  	f*c&=g&f\blackdiamond g\blackdiamond c \label{8_4_4}\\
  	c*f&=d&c\blackdiamond d\blackdiamond f \label{8_4_5}\\
  	f*a&=e&f\blackdiamond e\blackdiamond a \label{8_4_6}\\
  	a*e&=h&a\blackdiamond h\blackdiamond e \label{8_4_7}\\
  	e*h&=d&e\blackdiamond d\blackdiamond h \label{8_4_8}
  \end{align}

We have the table as follows.

\begin{longtblr}[caption = {Orders that could be extendable to biorders on $Q(8_4)$}, label = {tbl8_4}]
{colspec = {|c|c|l|l|l|}, rowhead = 2}
            \hline
  			\text{Sr.} & \text{Label} & \text{Expression} & \text{Obtained by} & \text{{\color{red} Contradiction to} $\backslash$}\\
  			\text{No.} &  & &  & \text{\color{blue} Compatible with}\\[2mm]
  			\hline
  			1&\text{A} & $f\blackdiamond e\blackdiamond d\blackdiamond h\blackdiamond a\blackdiamond b$ &\ref{8_4_6},\;\ref{8_4_7},\;\ref{8_4_8},\;\ref{8_4_1}& \\
  			2&\text{$B_1$} & $f\blackdiamond e\blackdiamond d\blackdiamond c\blackdiamond h\blackdiamond a\blackdiamond b$ &\text{A},\;\ref{8_4_5}& \\
  			3&\text{$B_2$} & $f\blackdiamond e\blackdiamond d\blackdiamond h\blackdiamond c\blackdiamond a\blackdiamond b$ &\text{A},\;\ref{8_4_5}& \\
  			4&\text{$B_3$} & $f\blackdiamond e\blackdiamond d\blackdiamond h\blackdiamond a\blackdiamond c\blackdiamond b$ &\text{A},\;\ref{8_4_5}& \\
  			5&\text{$B_4$} & $f\blackdiamond e\blackdiamond d\blackdiamond h\blackdiamond a\blackdiamond b\blackdiamond c\blackdiamond g$ &\text{A},\;\ref{8_4_5},\;\ref{8_4_3}&  {\color{red} \ref{8_4_2},\;\color{red} \ref{8_4_4}} \\
  			6&\text{$C_1$} & $f\blackdiamond e\blackdiamond d\blackdiamond g\blackdiamond h\blackdiamond a\blackdiamond c\blackdiamond b$ &\text{$B_3$},\;\ref{8_4_2}&  {\color{blue} \ref{8_4_3},\;\color{blue} \ref{8_4_4}} \\
  			7&\text{$C_2$} & $f\blackdiamond e\blackdiamond g\blackdiamond d\blackdiamond h\blackdiamond a\blackdiamond c\blackdiamond b$ &\text{$B_3$},\;\ref{8_4_2}&  {\color{blue} \ref{8_4_3},\;\color{blue} \ref{8_4_4}} \\
  			8&\text{$C_3$} & $f\blackdiamond g\blackdiamond e\blackdiamond d\blackdiamond h\blackdiamond a\blackdiamond c\blackdiamond b$ &\text{$B_3$},\;\ref{8_4_2}&  {\color{blue} \ref{8_4_3},\;\color{blue} \ref{8_4_4}} \\
  			9&\text{$C_4$} & $g\blackdiamond f\blackdiamond e\blackdiamond d\blackdiamond h\blackdiamond a\blackdiamond c\blackdiamond b$ &\text{$B_3$},\;\ref{8_4_2}&  {\color{blue} \ref{8_4_3},\;\color{red} \ref{8_4_4}} \\
  			10&\text{$D_1$} & $f\blackdiamond e\blackdiamond d\blackdiamond g\blackdiamond h\blackdiamond c\blackdiamond a\blackdiamond b$ &\text{$B_2$},\;\ref{8_4_2}&  {\color{blue} \ref{8_4_3},\;\color{blue} \ref{8_4_4}} \\
  			11&\text{$D_2$} & $f\blackdiamond e\blackdiamond g\blackdiamond d\blackdiamond h\blackdiamond c\blackdiamond a\blackdiamond b$ &\text{$B_2$},\;\ref{8_4_2}&  {\color{blue} \ref{8_4_3},\;\color{blue} \ref{8_4_4}} \\
  	    	12&\text{$D_3$} & $f\blackdiamond g\blackdiamond e\blackdiamond d\blackdiamond h\blackdiamond c\blackdiamond a\blackdiamond b$ &\text{$B_2$},\;\ref{8_4_2}&  {\color{blue} \ref{8_4_3},\;\color{blue} \ref{8_4_4}} \\
  	   		13&\text{$D_4$} & $g\blackdiamond f\blackdiamond e\blackdiamond d\blackdiamond h\blackdiamond c\blackdiamond a\blackdiamond b$ &\text{$B_2$},\;\ref{8_4_2}&  {\color{blue} \ref{8_4_3},\;\color{red} \ref{8_4_4}} \\
  	   		14&\text{$E_1$} & $f\blackdiamond g\blackdiamond e\blackdiamond d\blackdiamond c\blackdiamond h\blackdiamond a\blackdiamond b$ &\text{$B_1$},\;\ref{8_4_4}&  {\color{blue} \ref{8_4_2},\;\color{blue} \ref{8_4_3}}\\
  	   		15&\text{$E_2$} & $f\blackdiamond e\blackdiamond g\blackdiamond d\blackdiamond c\blackdiamond h\blackdiamond a\blackdiamond b$ &\text{$B_1$},\;\ref{8_4_4}&  {\color{blue} \ref{8_4_2},\;\color{blue} \ref{8_4_3}}\\
  			16&\text{$E_3$} & $f\blackdiamond e\blackdiamond d\blackdiamond g\blackdiamond c\blackdiamond h\blackdiamond a\blackdiamond b$ &\text{$B_1$},\;\ref{8_4_4}&  {\color{blue} \ref{8_4_2},\;\color{blue} \ref{8_4_3}}\\

  				\hline
  	\end{longtblr}

Linear orders on the generating set of $Q(8_4)$ as in \text{$C_1$}, \text{$C_2$}, \text{$C_3$}, \text{$D_1$},  \text{$D_2$},  \text{$D_3$},  \text{$E_1$}, \text{$E_2$} and \text{$E_3$} in Table \ref{tbl8_4} could be extendable to biorders on $Q(8_4)$.

\subsection{Knot $8_5$}\label{8_5}

\begin{figure}[H]
 	\centering
 	\includegraphics[scale=0.3]{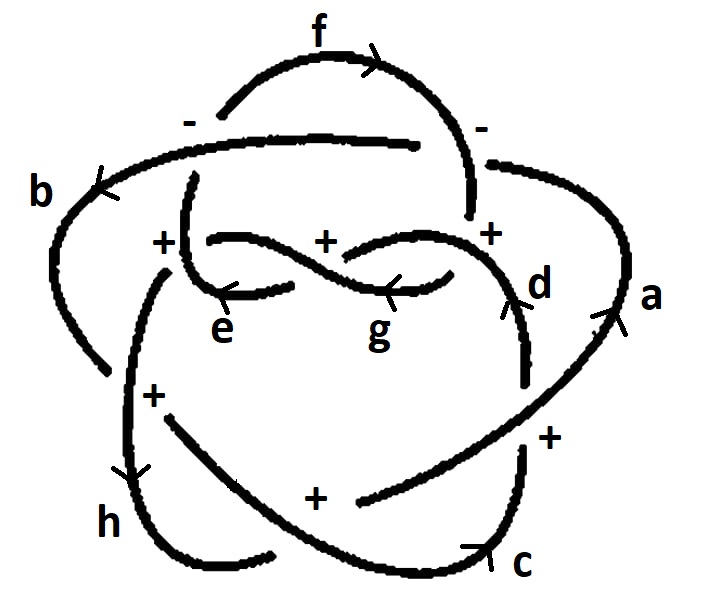}
 	\caption{Knot $8_5$}
    \label{fig8_5}
 \end{figure}

 Consider the diagram of the knot $8_5$ as in Figure \ref{fig8_5}. The quandle $Q(8_5)$ is generated by labelings of arcs of the diagram, and the relations are given by
   \begin{align}
 	f*b&=e&f\blackdiamond e\blackdiamond b \label{8_5_1}\\
 	b*f&=a&b\blackdiamond a\blackdiamond f \label{8_5_2}\\
 	f*d&=g&f\blackdiamond g\blackdiamond d \label{8_5_3}\\
 	d*g&=e&d\blackdiamond e\blackdiamond g \label{8_5_4}\\
 	g*e&=h&g\blackdiamond h\blackdiamond e \label{8_5_5}\\
 	b*h&=c&b\blackdiamond c\blackdiamond h \label{8_5_6}\\
 	h*c&=a&h\blackdiamond a\blackdiamond c \label{8_5_7}\\
 	c*a&=d&c\blackdiamond d\blackdiamond a \label{8_5_8}
 \end{align}

 We have the table as follows.

\begin{longtblr}[caption = {Orders that could be extendable to biorders on $Q(8_5)$}, label = {tbl8_5}]
{colspec = {|c|c|l|l|l|}, rowhead = 2}
            \hline
 			\text{Sr.} & \text{Label} & \text{Expression} & \text{Obtained by} & \text{{\color{red} Contradiction to} $\backslash$}\\
 			\text{No.} &  & &  & \text{\color{blue} Compatible with}\\[2mm]
 			\hline
 	1& \text{A}& $h\blackdiamond a\blackdiamond d\blackdiamond c\blackdiamond b$ &\ref{8_5_7},\;\ref{8_5_8},\;\ref{8_5_6}& \\
 	2& \text{$B_1$}& $h\blackdiamond f\blackdiamond a\blackdiamond d\blackdiamond c\blackdiamond b$ &\text{A},\;\ref{8_5_2}& \\
 	3& \text{$B_2$}& $f\blackdiamond h\blackdiamond a\blackdiamond d\blackdiamond c\blackdiamond b$ &\text{A},\;\ref{8_5_2}& \\
 	4& \text{$C_1$}& $e\blackdiamond h\blackdiamond f\blackdiamond g\blackdiamond a\blackdiamond d\blackdiamond c\blackdiamond b$ &\text{$B_1$},\;\ref{8_5_3},\;\ref{8_5_5}&	{\color{red} \ref{8_5_1},}\;{\color{red} \ref{8_5_4}} \\
    5& \text{$C_2$}& $h\blackdiamond f\blackdiamond a\blackdiamond g\blackdiamond e\blackdiamond d\blackdiamond c\blackdiamond b$    &\text{$B_1$},\;\ref{8_5_3},\;\ref{8_5_4}&	{\color{blue} \ref{8_5_1},}\;{\color{red} \ref{8_5_5}} \\
 	6& \text{$D_1$}& $f\blackdiamond g\blackdiamond h\blackdiamond a\blackdiamond d\blackdiamond c\blackdiamond b$ &\text{$B_2$},\;\ref{8_5_3}& \\
 	7& \text{$D_2$}& $f\blackdiamond h\blackdiamond g\blackdiamond a\blackdiamond d\blackdiamond c\blackdiamond b$ &\text{$B_2$},\;\ref{8_5_3}& \\
    8& \text{$D_3$}& $f\blackdiamond h\blackdiamond a\blackdiamond g\blackdiamond e\blackdiamond d\blackdiamond c\blackdiamond b$     &\text{$B_2$},\;\ref{8_5_3},\;\ref{8_5_4}&	{\color{blue} \ref{8_5_1},}\;{\color{red} \ref{8_5_5}} \\
    9& \text{$E_1$}& $f\blackdiamond g\blackdiamond e\blackdiamond h\blackdiamond a\blackdiamond d\blackdiamond c\blackdiamond b$     &\text{$D_1$},\;\ref{8_5_4}&	{\color{blue} \ref{8_5_1},}\;{\color{red} \ref{8_5_5}} \\
    10& \text{$E_2$}& $f\blackdiamond g\blackdiamond h\blackdiamond e\blackdiamond a\blackdiamond d\blackdiamond c\blackdiamond b$     &\text{$D_1$},\;\ref{8_5_4}&	{\color{blue} \ref{8_5_1},}\;{\color{blue} \ref{8_5_5}} \\
    11& \text{$E_3$}& $f\blackdiamond g\blackdiamond h\blackdiamond a\blackdiamond e\blackdiamond d\blackdiamond c\blackdiamond b$    &\text{$D_1$},\;\ref{8_5_4}&	{\color{blue} \ref{8_5_1},}\;{\color{blue} \ref{8_5_5}} \\
    12& \text{$F_1$}& $f\blackdiamond h\blackdiamond g\blackdiamond e\blackdiamond a\blackdiamond d\blackdiamond c\blackdiamond b$    &\text{$D_2$},\;\ref{8_5_4}&	{\color{blue} \ref{8_5_1},}\;{\color{red} \ref{8_5_5}} \\
    13& \text{$F_2$}& $f\blackdiamond h\blackdiamond g\blackdiamond a\blackdiamond e\blackdiamond d\blackdiamond c\blackdiamond b$   &\text{$D_2$},\;\ref{8_5_4}&	{\color{blue} \ref{8_5_1},}\;{\color{red} \ref{8_5_5}} \\

		\hline
       \end{longtblr}

Linear orders on the generating set of $Q(8_5)$ as in \text{$E_2$} and \text{$E_3$} (in Table \ref{tbl8_5}) could be extendable to biorders on $Q(8_5)$.

\subsection{Knot $8_6$}\label{8_6}

\begin{figure}[H]
	\centering
	\includegraphics[scale=0.3]{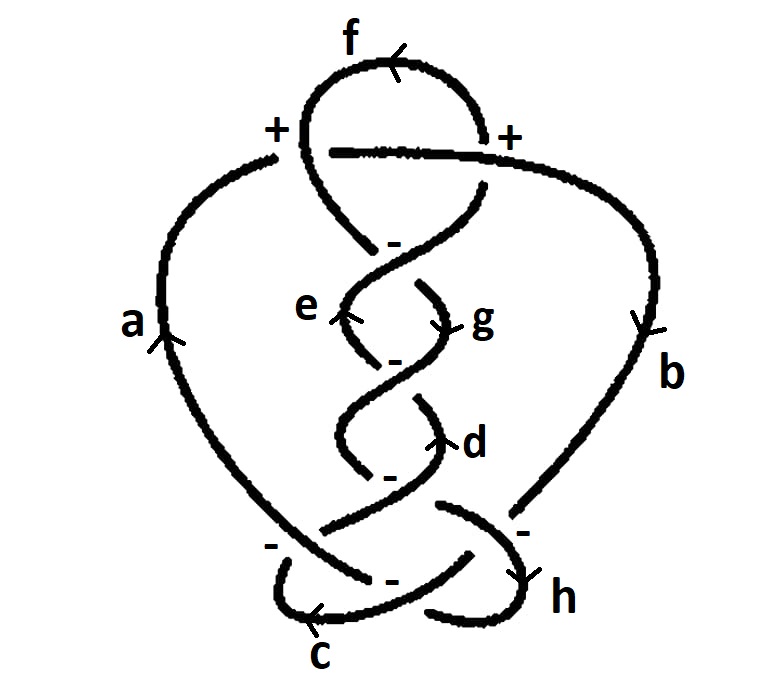}
	\caption{Knot $8_6$}
    \label{fig8_6}
\end{figure}

Consider the diagram of the knot $8_6$ as in Figure \ref{fig8_6}. The quandle $Q(8_6)$ is generated by labelings of arcs of the diagram, and the relations are given by
\begin{align}
	a*f&=b&a\blackdiamond b\blackdiamond f \label{8_6_1}\\
	e*b&=f&e\blackdiamond f\blackdiamond b \label{8_6_2}\\
	g*e&=f&g\blackdiamond f\blackdiamond e \label{8_6_3}\\
	e*g&=d&e\blackdiamond d\blackdiamond g \label{8_6_4}\\
	h*d&=g&h\blackdiamond g\blackdiamond d \label{8_6_5}\\
	d*a&=c&d\blackdiamond c\blackdiamond a \label{8_6_6}\\
	a*c&=h&a\blackdiamond h\blackdiamond c \label{8_6_7}\\
	c*h&=b&c\blackdiamond b\blackdiamond h \label{8_6_8}
\end{align}

We have the table as follows.

\begin{longtblr}[caption = {Orders that could be extendable to biorders on $Q(8_6)$}, label = {tbl8_6}]
{colspec = {|c|c|l|l|l|}, rowhead = 2}
            \hline
			\text{Sr.} & \text{Label} & \text{Expression} & \text{Obtained by} & \text{{\color{red} Contradiction to} $\backslash$}\\
			\text{No.} &  & &  & \text{\color{blue} Compatible with}\\[2mm]
			\hline
	1& \text{A}& $d\blackdiamond c\blackdiamond b\blackdiamond h\blackdiamond a$ &\ref{8_6_8},\;\ref{8_6_7},\;\ref{8_6_6}& \\
	2& \text{$B_1$}& $e\blackdiamond d\blackdiamond g\blackdiamond c\blackdiamond b\blackdiamond h\blackdiamond a$ & \text{A},\;\ref{8_6_5},\;\ref{8_6_4}&\\
	3& \text{$B_2$}& $e\blackdiamond d\blackdiamond c\blackdiamond g\blackdiamond b\blackdiamond h\blackdiamond a$ & \text{A},\;\ref{8_6_5},\;\ref{8_6_4}&\\
	4& \text{$B_3$}& $e\blackdiamond d\blackdiamond c\blackdiamond b\blackdiamond g\blackdiamond h\blackdiamond a$ & \text{A},\;\ref{8_6_5},\;\ref{8_6_4}&\\
	5&\text{$C_1$}& $e\blackdiamond f\blackdiamond d\blackdiamond g\blackdiamond c\blackdiamond b\blackdiamond h\blackdiamond a$ & \text{$B_1$},\;\ref{8_6_3}& {\color{blue} \ref{8_6_1},\;\ref{8_6_2}}\\
	6&\text{$C_2$}& $e\blackdiamond d\blackdiamond f\blackdiamond g\blackdiamond c\blackdiamond b\blackdiamond h\blackdiamond a$ & \text{$B_1$},\;\ref{8_6_3}& {\color{blue} \ref{8_6_1},\;\ref{8_6_2}}\\
	7&\text{$D_1$}& $e\blackdiamond f\blackdiamond d\blackdiamond c\blackdiamond g\blackdiamond b\blackdiamond h\blackdiamond a$ & \text{$B_2$},\;\ref{8_6_3}& {\color{blue} \ref{8_6_1},\;\ref{8_6_2}}\\
	8&\text{$D_2$}& $e\blackdiamond d\blackdiamond f\blackdiamond c\blackdiamond g\blackdiamond b\blackdiamond h\blackdiamond a$ & \text{$B_2$},\;\ref{8_6_3}& {\color{blue} \ref{8_6_1},\;\ref{8_6_2}}\\
	9&\text{$D_3$}& $e\blackdiamond d\blackdiamond c\blackdiamond f\blackdiamond g\blackdiamond b\blackdiamond h\blackdiamond a$ & \text{$B_2$},\;\ref{8_6_3}& {\color{blue} \ref{8_6_1},\;\ref{8_6_2}}\\
	10&\text{$E_1$}& $e\blackdiamond f\blackdiamond d\blackdiamond c\blackdiamond b\blackdiamond g\blackdiamond h\blackdiamond a$ & \text{$B_3$},\;\ref{8_6_2}& {\color{blue} \ref{8_6_1},\;\ref{8_6_3}}\\
	11&\text{$E_2$}& $e\blackdiamond d\blackdiamond f\blackdiamond c\blackdiamond b\blackdiamond g\blackdiamond h\blackdiamond a$ & \text{$B_3$},\;\ref{8_6_2}& {\color{blue} \ref{8_6_1},\;\ref{8_6_3}}\\
	12&\text{$E_3$}& $e\blackdiamond d\blackdiamond c\blackdiamond f\blackdiamond b\blackdiamond g\blackdiamond h\blackdiamond a$ & \text{$B_3$},\;\ref{8_6_2}& {\color{blue} \ref{8_6_1},\;\ref{8_6_3}}\\

	\hline
	\end{longtblr}

In Table \ref{tbl8_6}, linear orderings \text{$C_1$}, \text{$C_2$}, \text{$D_1$}, \text{$D_2$}, \text{$D_3$}, \text{$E_1$}, \text{$E_2$} and \text{$E_3$} are compatible with all quandle relations. These orders could be extendable to biorders on $Q(8_6)$.

\subsection{Knot $8_7$}\label{8_7}

\begin{figure}[H]
	\centering
	\includegraphics[scale=0.3]{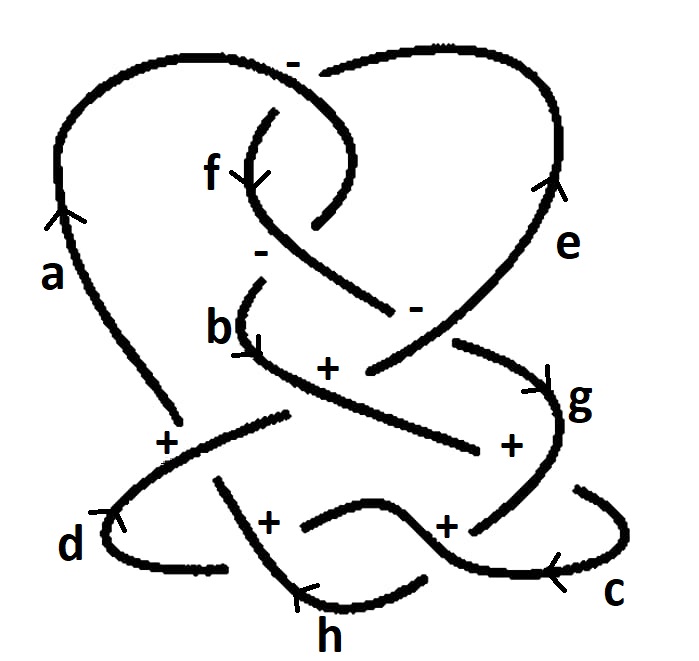}
	\caption{Knot $8_7$}
     \label{fig8_7}
\end{figure}

Consider the diagram of the knot $8_7$ as in Figure \ref{fig8_7}. The quandle $Q(8_7)$ is generated by labelings of arcs of the diagram, and the relations are given by
\begin{align}
	f*a&=e&f\blackdiamond e\blackdiamond a \label{8_7_1}\\
	b*f&=a&b\blackdiamond a\blackdiamond f \label{8_7_2}\\
	g*e&=f&g\blackdiamond f\blackdiamond e \label{8_7_3}\\
	d*b&=e&d\blackdiamond e\blackdiamond b \label{8_7_4}\\
	h*d&=a&h\blackdiamond a\blackdiamond d \label{8_7_5}\\
	c*h&=d&c\blackdiamond d\blackdiamond h \label{8_7_6}\\
	g*c&=h&g\blackdiamond h\blackdiamond c \label{8_7_7}\\
	b*g&=c&b\blackdiamond c\blackdiamond g \label{8_7_8}
\end{align}

We have the table as follows.

\begin{longtblr}[caption = {Orders that could be extendable to biorders on $Q(8_7)$}, label = {tbl8_7}]
{colspec = {|c|c|l|l|l|}, rowhead = 2}
            \hline
			\text{Sr.} & \text{Label} & \text{Expression} & \text{Obtained by} & \text{{\color{red} Contradiction to} $\backslash$}\\
			\text{No.} &  &  &  & \text{\color{blue} Compatible with}\\[2mm]
			\hline
			1& \text{A} & $b\blackdiamond a\blackdiamond e\blackdiamond f\blackdiamond g$ & \ref{8_7_1},\;\ref{8_7_2},\;\ref{8_7_3}& \\
			2& \text{$B_1$} & $b\blackdiamond a\blackdiamond e\blackdiamond d\blackdiamond f\blackdiamond g$ & \text{A},\;\ref{8_7_4}& \\
			3& \text{$B_2$} & $b\blackdiamond a\blackdiamond e\blackdiamond f\blackdiamond d\blackdiamond g$ & \text{A},\;\ref{8_7_4}& \\
	    	4& \text{$B_3$} & $b\blackdiamond a\blackdiamond e\blackdiamond f\blackdiamond g\blackdiamond d$ & \text{A},\;\ref{8_7_4}& \\
	    	5& \text{$C_1$} & $b\blackdiamond h\blackdiamond a\blackdiamond e\blackdiamond d\blackdiamond f\blackdiamond g$ &\text{$B_1$},\;\ref{8_7_5}& \\
	    	6& \text{$C_2$} & $c\blackdiamond h\blackdiamond b\blackdiamond a\blackdiamond e\blackdiamond d\blackdiamond f\blackdiamond g$ &\text{$B_1$},\;\ref{8_7_5},\;\ref{8_7_7}&{\color{red} \ref{8_7_6},\;\ref{8_7_8}} \\
	   		7& \text{$D_1$} & $b\blackdiamond c\blackdiamond h\blackdiamond a\blackdiamond e\blackdiamond d\blackdiamond f\blackdiamond g$ &\text{$C_1$},\;\ref{8_7_7}&{\color{red} \ref{8_7_6},}\;{\color{blue} \ref{8_7_8}} \\
	   		8& \text{$D_2$} & $c\blackdiamond b\blackdiamond h\blackdiamond a\blackdiamond e\blackdiamond d\blackdiamond f\blackdiamond g$ &\text{$C_1$},\;\ref{8_7_7}&{\color{red} \ref{8_7_6},}\;{\color{red} \ref{8_7_8}} \\
	    	9& \text{$E_1$} & $b\blackdiamond h\blackdiamond a\blackdiamond e\blackdiamond f\blackdiamond d\blackdiamond g$ &\text{$B_2$},\;\ref{8_7_5}& \\
	    	10& \text{$E_2$} & $c\blackdiamond h\blackdiamond b\blackdiamond a\blackdiamond e\blackdiamond f\blackdiamond d\blackdiamond g$ &\text{$B_2$},\;\ref{8_7_5},\;\ref{8_7_7}&{\color{red} \ref{8_7_6},}\;{\color{red} \ref{8_7_8}} \\
	    	11& \text{$F_1$} & $b\blackdiamond c\blackdiamond h\blackdiamond a\blackdiamond e\blackdiamond f\blackdiamond d\blackdiamond g$ &\text{$E_1$},\;\ref{8_7_7}&{\color{red} \ref{8_7_6},}\;{\color{blue} \ref{8_7_8}} \\
	    	12& \text{$F_2$} & $c\blackdiamond b\blackdiamond h\blackdiamond a\blackdiamond e\blackdiamond f\blackdiamond d\blackdiamond g$ &\text{$E_1$},\;\ref{8_7_7}&{\color{red} \ref{8_7_6},}\;{\color{red} \ref{8_7_8}} \\
	    	13& \text{$G_1$} & $b\blackdiamond h\blackdiamond a\blackdiamond e\blackdiamond f\blackdiamond g\blackdiamond d$ &\text{$B_3$},\;\ref{8_7_5}& \\
	    	14& \text{$G_2$} & $c\blackdiamond h\blackdiamond b\blackdiamond a\blackdiamond e\blackdiamond f\blackdiamond g\blackdiamond d$ &\text{$B_3$},\;\ref{8_7_5},\;\ref{8_7_7}&{\color{red} \ref{8_7_6},\;\ref{8_7_8}}\\
	   		15& \text{$H_1$} & $b\blackdiamond c\blackdiamond h\blackdiamond a\blackdiamond e\blackdiamond f\blackdiamond g\blackdiamond d$ &\text{$G_1$},\;\ref{8_7_7}&{\color{red} \ref{8_7_6},}\;{\color{blue} \ref{8_7_8}}\\
	   		16& \text{$H_2$} & $c\blackdiamond b\blackdiamond h\blackdiamond a\blackdiamond e\blackdiamond f\blackdiamond g\blackdiamond d$ &\text{$G_1$},\;\ref{8_7_7}&{\color{red} \ref{8_7_6},}\;{\color{red} \ref{8_7_8}}\\
\hline
\end{longtblr}

In Table \ref{tbl8_7}, every possible linear order (i.e. orders \text{$C_2$}, \text{$D_1$}, \text{$D_2$}, \text{$E_2$}, \text{$F_1$}, \text{$F_2$}, \text{$G_2$}, \text{$H_1$} and \text{$H_2$}) on the generating set contradict either (\ref{8_7_6}) or (\ref{8_7_8}), or both. Thus, we conclude that $Q(8_7)$ can not be biorderable.

\subsection{Knot $8_8$}\label{8_8}

\begin{figure}[H]
		\centering
		\includegraphics[scale=0.3]{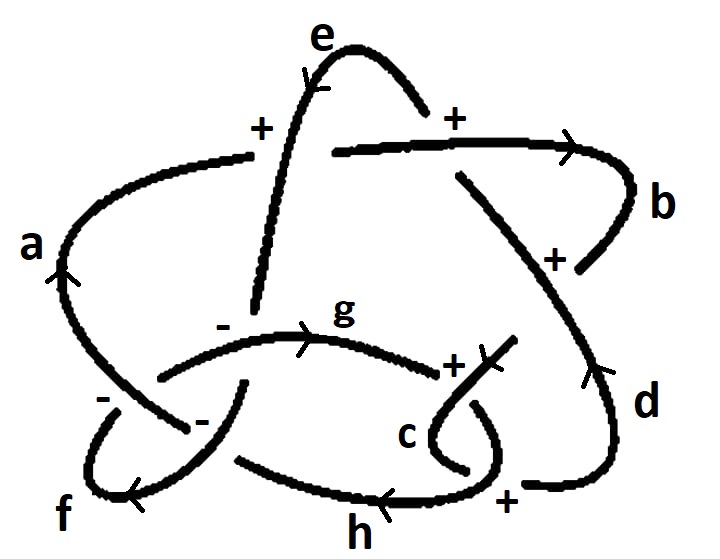}
		\caption{Knot $8_8$}
        \label{fig8_8}
	\end{figure}

Consider the diagram of the knot $8_8$ as in Figure \ref{fig8_8}. The quandle $Q(8_8)$ is generated by labelings of arcs of the diagram, and the relations are given by

	\begin{align}
		a*e&=b&a\blackdiamond b\blackdiamond e \label{8_8_1}\\
		d*b&=e&d\blackdiamond e\blackdiamond b \label{8_8_2}\\
		b*d&=c&b\blackdiamond c\blackdiamond d \label{8_8_3}\\
		g*c&=h&g\blackdiamond h\blackdiamond c \label{8_8_4}\\
		c*h&=d&c\blackdiamond d\blackdiamond h \label{8_8_5}\\
		a*f&=h&a\blackdiamond h\blackdiamond f \label{8_8_6}\\
		g*a&=f&g\blackdiamond f\blackdiamond a \label{8_8_7}\\
		f*g&=e&f\blackdiamond e\blackdiamond g \label{8_8_8}
	\end{align}

We have the table as follows.

\begin{longtblr}[caption = {Orders that could be extendable to biorders on $Q(8_8)$}, label = {tbl8_8}]
{colspec = {|c|c|l|l|l|}, rowhead = 2}
            \hline
			\text{Sr.} & \text{Label} & \text{Expression} & \text{Obtained by} & \text{{\color{red} Contradiction to} $\backslash$}\\
			\text{No.} &  & &  & \text{\color{blue} Compatible with}\\[2mm]
			\hline
	1& \text{A}& $b\blackdiamond c\blackdiamond d\blackdiamond h\blackdiamond g$ & \ref{8_8_3},\;\ref{8_8_5},\;\ref{8_8_4}& \\
	2& \text{$B_1$}& $a\blackdiamond b\blackdiamond e\blackdiamond c\blackdiamond d\blackdiamond h\blackdiamond g$ & \text{A},\;\ref{8_8_2},\;\ref{8_8_1}& \\
	3& \text{$B_2$}& $a\blackdiamond b\blackdiamond c\blackdiamond e\blackdiamond d\blackdiamond h\blackdiamond g$ & \text{A},\;\ref{8_8_2},\;\ref{8_8_1}& \\
	4& \text{$C_1$}& $a\blackdiamond b\blackdiamond e\blackdiamond c\blackdiamond d\blackdiamond h\blackdiamond f\blackdiamond g$ & \text{$B_1$},\;\ref{8_8_6}& {\color{blue} \ref{8_8_7},}\;{\color{red} \ref{8_8_8}} \\
	5& \text{$C_2$}& $a\blackdiamond b\blackdiamond e\blackdiamond c\blackdiamond d\blackdiamond h\blackdiamond g\blackdiamond f$ & \text{$B_1$},\;\ref{8_8_6}& {\color{red} \ref{8_8_7},}\;{\color{red} \ref{8_8_8}} \\
	6& \text{$D_1$}& $a\blackdiamond b\blackdiamond c\blackdiamond e\blackdiamond d\blackdiamond h\blackdiamond f\blackdiamond g$ & \text{$B_2$},\;\ref{8_8_6}& {\color{blue} \ref{8_8_7},}\;{\color{red} \ref{8_8_8}} \\
	7& \text{$D_2$}& $a\blackdiamond b\blackdiamond c\blackdiamond e\blackdiamond d\blackdiamond h\blackdiamond g\blackdiamond f$ & \text{$B_2$},\;\ref{8_8_6}& {\color{red} \ref{8_8_7},}\;{\color{red} \ref{8_8_8}} \\

	\hline
\end{longtblr}

In Table \ref{tbl8_8}, linear orders \text{$C_1$}, \text{$C_2$}, \text{$D_1$} and \text{$D_2$} contradict either (\ref{8_8_7}) or (\ref{8_8_8}), or both. Therefore, we conclude that $Q(8_8)$ can not be biorderable.

\subsection{Knot $8_9$}\label{8_9}

\begin{figure}[H]
	\centering
	\includegraphics[scale=0.3]{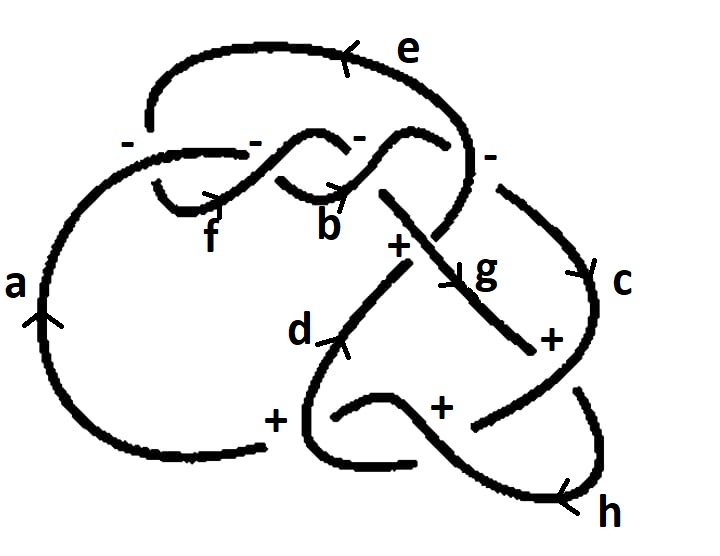}
	\caption{Knot $8_9$}
    \label{fig8_9}
\end{figure}

Consider the diagram of the knot $8_9$ as in Figure \ref{fig8_9}. The quandle $Q(8_9)$ is generated by labelings of arcs of the diagram, and the relations are given by

\begin{align}
	f*a&=e&f\blackdiamond e\blackdiamond a \label{8_9_1}\\
	b*f&=a&b\blackdiamond a\blackdiamond f \label{8_9_2}\\
	g*b&=f&g\blackdiamond f\blackdiamond b \label{8_9_3}\\
	c*e&=b&c\blackdiamond b\blackdiamond e \label{8_9_4}\\
	d*g&=e&d\blackdiamond e\blackdiamond g \label{8_9_5}\\
	g*c&=h&g\blackdiamond h\blackdiamond c \label{8_9_6}\\
	c*h&=d&c\blackdiamond d\blackdiamond h \label{8_9_7}\\
	h*d&=a&h\blackdiamond a\blackdiamond d \label{8_9_8}
\end{align}

We have the table as follows.

\begin{longtblr}[caption = {Orders that could be extendable to biorders on $Q(8_9)$}, label = {tbl8_9}]
{colspec = {|c|c|l|l|l|}, rowhead = 2}
            \hline
			\text{Sr.} & \text{Label} & \text{Expression} & \text{Obtained by} & \text{{\color{red} Contradiction to} $\backslash$}\\
			\text{No.} &  & &  & \text{\color{blue} Compatible with}\\[2mm]
			\hline
1& \text{A}& $g\blackdiamond f\blackdiamond e\blackdiamond a\blackdiamond b\blackdiamond c$ &\ref{8_9_1},\;\ref{8_9_2},\;\ref{8_9_3},\;\ref{8_9_4}& \\
2& \text{$B_1$} & $g\blackdiamond f\blackdiamond e\blackdiamond d\blackdiamond a\blackdiamond b\blackdiamond c$ & \text{A},\;\ref{8_9_5}& \\
3& \text{$B_2$} & $g\blackdiamond f\blackdiamond e\blackdiamond a\blackdiamond d\blackdiamond b\blackdiamond c$ & \text{A},\;\ref{8_9_5}& \\
4& \text{$B_3$} & $g\blackdiamond f\blackdiamond e\blackdiamond a\blackdiamond b\blackdiamond d\blackdiamond c$ & \text{A},\;\ref{8_9_5}& \\
5& \text{$B_4$} & $g\blackdiamond f\blackdiamond e\blackdiamond a\blackdiamond b\blackdiamond c\blackdiamond d\blackdiamond h$ & \text{A},\;\ref{8_9_5},\;\ref{8_9_7}&{\color{red} \ref{8_9_6},}\;{\color{red} \ref{8_9_8}} \\
6& \text{$C_1$} & $g\blackdiamond f\blackdiamond e\blackdiamond d\blackdiamond a\blackdiamond h\blackdiamond b\blackdiamond c$ & \text{$B_1$},\;\ref{8_9_8}&{\color{blue} \ref{8_9_6},}\;{\color{red} \ref{8_9_7}} \\
7& \text{$C_2$} & $g\blackdiamond f\blackdiamond e\blackdiamond d\blackdiamond a\blackdiamond b\blackdiamond h\blackdiamond c$ & \text{$B_1$},\;\ref{8_9_8}&{\color{blue} \ref{8_9_6},}\;{\color{red} \ref{8_9_7}} \\
8& \text{$C_3$} & $g\blackdiamond f\blackdiamond e\blackdiamond d\blackdiamond a\blackdiamond b\blackdiamond c\blackdiamond h$ & \text{$B_1$},\;\ref{8_9_8}&{\color{red} \ref{8_9_6},}\;{\color{red} \ref{8_9_7}} \\
9& \text{$D_1$} & $g\blackdiamond f\blackdiamond e\blackdiamond h\blackdiamond a\blackdiamond d\blackdiamond b\blackdiamond c$ & \text{$B_2$},\;\ref{8_9_8}&{\color{blue} \ref{8_9_6},}\;{\color{blue} \ref{8_9_7}} \\
10& \text{$D_2$} & $g\blackdiamond f\blackdiamond h\blackdiamond e\blackdiamond a\blackdiamond d\blackdiamond b\blackdiamond c$ & \text{$B_2$},\;\ref{8_9_8}&{\color{blue} \ref{8_9_6},}\;{\color{blue} \ref{8_9_7}} \\
11& \text{$D_3$} & $g\blackdiamond h\blackdiamond f\blackdiamond e\blackdiamond a\blackdiamond d\blackdiamond b\blackdiamond c$ & \text{$B_2$},\;\ref{8_9_8}&{\color{blue} \ref{8_9_6},}\;{\color{blue} \ref{8_9_7}} \\
12& \text{$D_4$} & $h\blackdiamond g\blackdiamond f\blackdiamond e\blackdiamond a\blackdiamond d\blackdiamond b\blackdiamond c$ & \text{$B_2$},\;\ref{8_9_8}&{\color{red} \ref{8_9_6},}\;{\color{blue} \ref{8_9_7}} \\
13& \text{$E_1$} & $g\blackdiamond f\blackdiamond e\blackdiamond h\blackdiamond a\blackdiamond b\blackdiamond d\blackdiamond c$ & \text{$B_3$},\;\ref{8_9_8}&{\color{blue} \ref{8_9_6},}\;{\color{blue} \ref{8_9_7}} \\
14& \text{$E_2$} & $g\blackdiamond f\blackdiamond h\blackdiamond e\blackdiamond a\blackdiamond b\blackdiamond d\blackdiamond c$ & \text{$B_3$},\;\ref{8_9_8}&{\color{blue} \ref{8_9_6},}\;{\color{blue} \ref{8_9_7}} \\
15& \text{$E_3$} & $g\blackdiamond h\blackdiamond f\blackdiamond e\blackdiamond a\blackdiamond b\blackdiamond d\blackdiamond c$ & \text{$B_3$},\;\ref{8_9_8}&{\color{blue} \ref{8_9_6},}\;{\color{blue} \ref{8_9_7}} \\
16& \text{$E_4$} & $h\blackdiamond g\blackdiamond f\blackdiamond e\blackdiamond a\blackdiamond b\blackdiamond d\blackdiamond c$ & \text{$B_3$},\;\ref{8_9_8}&{\color{red} \ref{8_9_6},}\;{\color{blue} \ref{8_9_7}} \\

		\hline
	\end{longtblr}

Linear orders on the generating set of $Q(8_9)$ as in  \text{$D_1$}, \text{$D_2$}, \text{$D_3$}, \text{$E_1$}, \text{$E_2$} and \text{$E_3$} in Table \ref{tbl8_9} could be extendable to biorders on $Q(8_9)$.

\subsection{Knot $8_{10}$}\label{8_10}

\begin{figure}[H]
	\centering
	\includegraphics[scale=0.3]{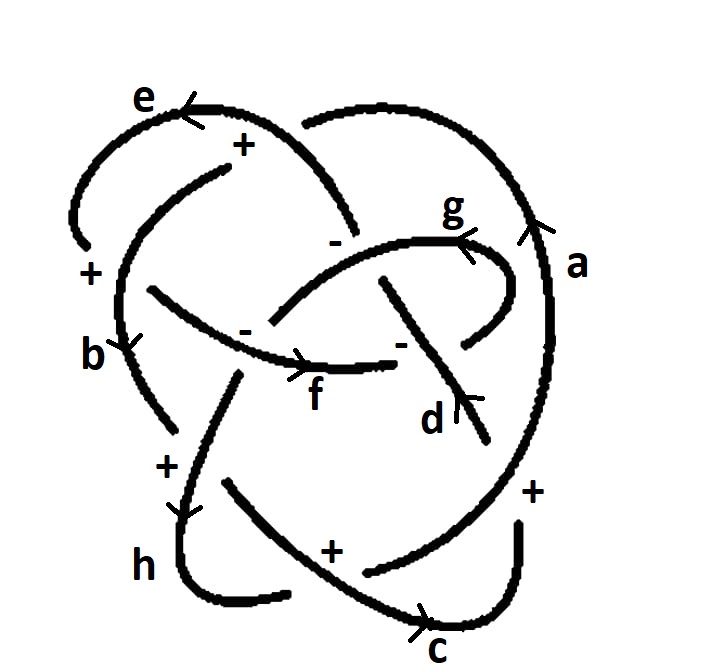}
	\caption{Knot $8_{10}$}
    \label{fig8_10}
\end{figure}

Consider the diagram of the knot $8_{10}$ as in Figure \ref{fig8_10}. The quandle $Q(8_{10})$ is generated by labelings of arcs of the diagram, and the relations are given by

\begin{align}
	a*e&=b&a\blackdiamond b\blackdiamond e \label{8_10_1}\\
	e*b&=f&e\blackdiamond f\blackdiamond b \label{8_10_2}\\
	b*h&=c&b\blackdiamond c\blackdiamond h \label{8_10_3}\\
	h*c&=a&h\blackdiamond a\blackdiamond c \label{8_10_4}\\
	c*a&=d&c\blackdiamond d\blackdiamond a \label{8_10_5}\\
	h*f&=g&h\blackdiamond g\blackdiamond f \label{8_10_6}\\
	e*g&=d&e\blackdiamond d\blackdiamond g \label{8_10_7}\\
	g*d&=f&g\blackdiamond f\blackdiamond d \label{8_10_8}
\end{align}

We have the table as follows.

\begin{longtblr}[caption = {Orders that could be extendable to biorders on $Q(8_{10})$}, label = {tbl8_10}]
{colspec = {|c|c|l|l|l|}, rowhead = 2}
            \hline
			\text{Sr.} & \text{Label} & \text{Expression} & \text{Obtained by} & \text{{\color{red} Contradiction to} $\backslash$}\\
			\text{No.} &  & &  & \text{\color{blue} Compatible with}\\[2mm]
			\hline
	1&\text{A} & $e\blackdiamond b\blackdiamond c\blackdiamond a\blackdiamond h$ & \ref{8_10_3},\;\ref{8_10_4},\;\ref{8_10_1} & \\
	2&\text{B} & $e\blackdiamond f\blackdiamond b\blackdiamond c\blackdiamond d\blackdiamond a\blackdiamond h$ & \text{A},\;\ref{8_10_2},\;\ref{8_10_5} & \\
	3&\text{$C_1$} & $e\blackdiamond g\blackdiamond f\blackdiamond b\blackdiamond c\blackdiamond d\blackdiamond a\blackdiamond h$ & \text{B},\;\ref{8_10_8} &{\color{red} \ref{8_10_6},\;\ref{8_10_7}} \\
	4&\text{$C_2$} & $g\blackdiamond e\blackdiamond f\blackdiamond b\blackdiamond c\blackdiamond d\blackdiamond a\blackdiamond h$ & \text{B},\;\ref{8_10_8} &{\color{red} \ref{8_10_6},\;\ref{8_10_7}} \\
    \hline
\end{longtblr}

In Table \ref{tbl8_10}, linear orders \text{$C_1$} and \text{$C_2$} on the generating set contradict (\ref{8_10_6}) and (\ref{8_10_7}). Therefore, we conclude that $Q(8_{10})$ can not be biorderable.

\subsection{Knot $8_{11}$}\label{8_11}

\begin{figure}[H]
	\centering
	\includegraphics[scale=0.3]{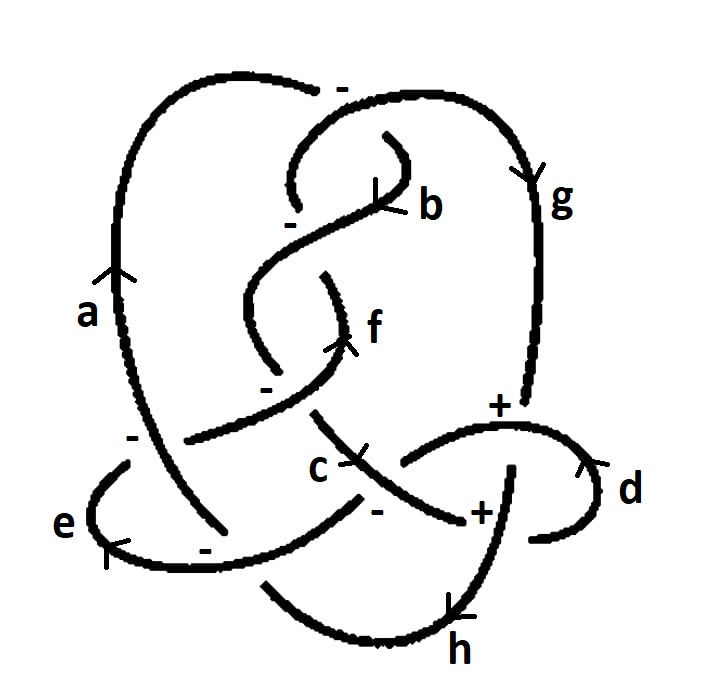}
	\caption{Knot $8_{11}$}
    \label{fig8_11}
\end{figure}

Consider the diagram of the knot $8_{11}$ as in Figure \ref{fig8_11}. The quandle $Q(8_{11})$ is generated by labelings of arcs of the diagram, and the relations are given by

\begin{align}
	b*g&=a&b\blackdiamond a\blackdiamond g \label{8_11_1}\\
	g*b&=f&g\blackdiamond f\blackdiamond b \label{8_11_2}\\
	c*f&=b&c\blackdiamond b\blackdiamond f \label{8_11_3}\\
	e*c&=d&e\blackdiamond d\blackdiamond c \label{8_11_4}\\
	a*e&=h&a\blackdiamond h\blackdiamond e \label{8_11_5}\\
	f*a&=e&f\blackdiamond e\blackdiamond a \label{8_11_6}\\
	g*d&=h&g\blackdiamond h\blackdiamond d \label{8_11_7}\\
	c*h&=d&c\blackdiamond d\blackdiamond h \label{8_11_8}
\end{align}

We have the table as follows.

\begin{longtblr}[caption = {Orders that could be extendable to biorders on $Q(8_{11})$}, label = {tbl8_11}]
{colspec = {|c|c|l|l|l|}, rowhead = 2}
            \hline
			\text{Sr.} & \text{Label} & \text{Expression} & \text{Obtained by} & \text{{\color{red} Contradiction to} $\backslash$}\\
			\text{No.} &  &  &  & \text{\color{blue} Compatible with}\\[2mm]
			\hline
1 & \text{A}& $c\blackdiamond d\blackdiamond h\blackdiamond g$ & \ref{8_11_7},\;\ref{8_11_8}& \\
2 & \text{$B_1$}& $c\blackdiamond d\blackdiamond e\blackdiamond h\blackdiamond g$ & \text{A},\;\ref{8_11_4}& \\
3 & \text{$B_2$}& $c\blackdiamond d\blackdiamond h\blackdiamond e\blackdiamond g$ & \text{A},\;\ref{8_11_4}& \\
4 & \text{$B_3$}& $c\blackdiamond d\blackdiamond h\blackdiamond g\blackdiamond e$ & \text{A},\;\ref{8_11_4}& \\
5 & \text{$C_1$}& $c\blackdiamond d\blackdiamond e\blackdiamond h\blackdiamond a\blackdiamond g$ & \text{$B_1$},\;\ref{8_11_5}& \\
6 & \text{$C_2$}& $c\blackdiamond d\blackdiamond e\blackdiamond h\blackdiamond g\blackdiamond a\blackdiamond b\blackdiamond f$ & \text{$B_1$},\;\ref{8_11_5},\;\ref{8_11_1},\;\ref{8_11_3}& {\color{red} \ref{8_11_2},\;\ref{8_11_6}} \\
7 & \text{$D_1$}& $c\blackdiamond d\blackdiamond f\blackdiamond e\blackdiamond h\blackdiamond a\blackdiamond g$ & \text{$C_1$},\;\ref{8_11_6}& \\
8 & \text{$D_2$}& $c\blackdiamond b\blackdiamond f\blackdiamond d\blackdiamond e\blackdiamond h\blackdiamond a\blackdiamond g$ & \text{$C_1$},\;\ref{8_11_6},\;\ref{8_11_3}& {\color{blue} \ref{8_11_1},\;\ref{8_11_2}} \\
9 & \text{$D_3$}& $f\blackdiamond b\blackdiamond c\blackdiamond d\blackdiamond e\blackdiamond h\blackdiamond a\blackdiamond g$ & \text{$C_1$},\;\ref{8_11_6},\;\ref{8_11_3}& {\color{blue} \ref{8_11_1},}\;{\color{red} \ref{8_11_2}} \\
10 & \text{$E_1$}& $c\blackdiamond b\blackdiamond d\blackdiamond f\blackdiamond e\blackdiamond h\blackdiamond a\blackdiamond g$ & \text{$D_1$},\;\ref{8_11_3}& {\color{blue} \ref{8_11_1},\;\ref{8_11_2}} \\
11 & \text{$E_2$}& $c\blackdiamond d\blackdiamond b\blackdiamond f\blackdiamond e\blackdiamond h\blackdiamond a\blackdiamond g$ & \text{$D_1$},\;\ref{8_11_3}& {\color{blue} \ref{8_11_1},\;\ref{8_11_2}} \\
12 & \text{$F_1$}& $c\blackdiamond d\blackdiamond a\blackdiamond h\blackdiamond e\blackdiamond g$ & \text{$B_2$},\;\ref{8_11_5}& \\
13 & \text{$F_2$}& $c\blackdiamond a\blackdiamond d\blackdiamond h\blackdiamond e\blackdiamond g$ & \text{$B_2$},\;\ref{8_11_5}& \\
14 & \text{$F_3$}& $f\blackdiamond b\blackdiamond a\blackdiamond c\blackdiamond d\blackdiamond h\blackdiamond e\blackdiamond g$ & \text{$B_2$},\;\ref{8_11_5},\;\ref{8_11_1},\;\ref{8_11_3}& {\color{red} \ref{8_11_2},\;\ref{8_11_6}} \\
15 & \text{$G_1$}& $c\blackdiamond d\blackdiamond a\blackdiamond h\blackdiamond e\blackdiamond f\blackdiamond g$ & \text{$F_1$},\;\ref{8_11_6}& \\
16 & \text{$G_2$}& $c\blackdiamond d\blackdiamond a\blackdiamond h\blackdiamond e\blackdiamond g\blackdiamond f\blackdiamond b$ & \text{$F_1$},\;\ref{8_11_6},\;\ref{8_11_2}& {\color{red} \ref{8_11_1},\;\ref{8_11_3}} \\
17 & \text{$H_1$}& $c\blackdiamond d\blackdiamond b\blackdiamond a\blackdiamond h\blackdiamond e\blackdiamond f\blackdiamond g$ & \text{$G_1$},\;\ref{8_11_1}& {\color{blue} \ref{8_11_2},\;\ref{8_11_3}} \\
18 & \text{$H_2$}& $c\blackdiamond b\blackdiamond d\blackdiamond a\blackdiamond h\blackdiamond e\blackdiamond f\blackdiamond g$ & \text{$G_1$},\;\ref{8_11_1}& {\color{blue} \ref{8_11_2},\;\ref{8_11_3}} \\
19 & \text{$H_3$}& $b\blackdiamond c\blackdiamond d\blackdiamond a\blackdiamond h\blackdiamond e\blackdiamond f\blackdiamond g$ & \text{$G_1$},\;\ref{8_11_1}& {\color{blue} \ref{8_11_2},}\;{\color{red} \ref{8_11_3}} \\
20 & \text{$I_1$}& $c\blackdiamond b\blackdiamond a\blackdiamond d\blackdiamond h\blackdiamond e\blackdiamond g$ & \text{$F_2$},\;\ref{8_11_1}& \\
21 & \text{$I_2$}& $f\blackdiamond b\blackdiamond c\blackdiamond a\blackdiamond d\blackdiamond h\blackdiamond e\blackdiamond g$ & \text{$F_2$},\;\ref{8_11_1},\;\ref{8_11_3}& {\color{red} \ref{8_11_2},\;\ref{8_11_6}} \\
22 & \text{$J_1$}& $c\blackdiamond b\blackdiamond a\blackdiamond d\blackdiamond h\blackdiamond e\blackdiamond f\blackdiamond g$ & \text{$I_1$},\;\ref{8_11_6}& {\color{blue} \ref{8_11_2},\;\ref{8_11_3}} \\
23 & \text{$J_2$}& $c\blackdiamond b\blackdiamond a\blackdiamond d\blackdiamond h\blackdiamond e\blackdiamond g\blackdiamond f$ & \text{$I_1$},\;\ref{8_11_6}& {\color{red} \ref{8_11_2},\;}{\color{blue} \ref{8_11_3}} \\
24 & \text{$K_1$}& $c\blackdiamond d\blackdiamond a\blackdiamond h\blackdiamond g\blackdiamond e\blackdiamond f\blackdiamond b$ & \text{$B_3$},\;\ref{8_11_5},\;\ref{8_11_6},\;\ref{8_11_2}& {\color{red} \ref{8_11_1},\;}{\color{red} \ref{8_11_3}} \\
25 & \text{$K_2$}& $c\blackdiamond a\blackdiamond d\blackdiamond h\blackdiamond g\blackdiamond e\blackdiamond f\blackdiamond b$ & \text{$B_3$},\;\ref{8_11_5},\;\ref{8_11_6},\;\ref{8_11_2}& {\color{red} \ref{8_11_1},\;}{\color{red} \ref{8_11_3}} \\
26 & \text{$K_3$}& $a\blackdiamond c\blackdiamond d\blackdiamond h\blackdiamond g\blackdiamond e\blackdiamond f\blackdiamond b$ & \text{$B_3$},\;\ref{8_11_5},\;\ref{8_11_6},\;\ref{8_11_2}& {\color{red} \ref{8_11_1},\;}{\color{red} \ref{8_11_3}} \\

\hline
\end{longtblr}

Linear orders on the generating set of $Q(8_{11})$ as in \text{$D_2$}, \text{$E_1$}, \text{$E_2$}, \text{$H_1$}, \text{$H_2$} and \text{$J_1$} in Table \ref{tbl8_11} could be extendable to biorders on $Q(8_{11})$.

\subsection{Knot $8_{12}$}\label{8_12}

\begin{figure}[H]
	\centering
	\includegraphics[scale=0.3]{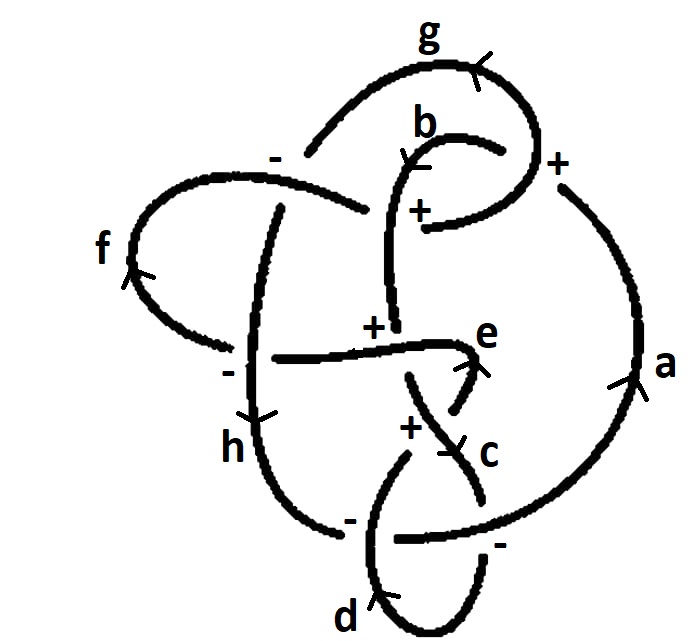}
	\caption{Knot $8_{12}$}
    \label{fig8_12}
\end{figure}

Consider the diagram of the knot $8_{12}$ as in Figure \ref{fig8_12}. The quandle $Q(8_{12})$ is generated by labelings of arcs of the diagram, and the relations are given by

\begin{align}
	a*g&=b&a\blackdiamond b\blackdiamond g \label{8_12_1}\\
	f*b&=g&f\blackdiamond g\blackdiamond b \label{8_12_2}\\
	b*e&=c&b\blackdiamond c\blackdiamond e \label{8_12_3}\\
	d*c&=e&d\blackdiamond e\blackdiamond c \label{8_12_4}\\
	d*a&=c&d\blackdiamond c\blackdiamond a \label{8_12_5}\\
	a*d&=h&a\blackdiamond h\blackdiamond d \label{8_12_6}\\
	f*h&=e&f\blackdiamond e\blackdiamond h \label{8_12_7}\\
	h*f&=g&h\blackdiamond g\blackdiamond f \label{8_12_8}
\end{align}

The determinant of the knot $8_{12}$ is $29$. By \cite[Proposition 4.4]{ms}, the generators $a$, $b$, $c$, $d$, $e$, $f$, $g$ and $h$ are all pairwise distinct. Applying the relations above ((\ref{8_12_1}) to (\ref{8_12_8})), we have the table as follows.

\begin{longtblr}[caption = {Orders that could be extendable to biorders on $Q(8_{12})$}, label = {tbl8_12}]
{colspec = {|c|c|l|l|l|}, rowhead = 2}
            \hline
			\text{Sr.} & \text{Label} & \text{Expression} & \text{Obtained by} & \text{{\color{red} Contradiction to} $\backslash$}\\
			\text{No.} &  &  &  & \text{\color{blue} Compatible with}\\[2mm]
			\hline
1 & \text{$A$} & $a\blackdiamond c\blackdiamond d$ & \ref{8_12_5} & \\
2 & \text{$B_1$} & $a\blackdiamond h\blackdiamond c\blackdiamond e\blackdiamond d$ & \text{$A$},\;\ref{8_12_6},\;\ref{8_12_4} & \\
3 & \text{$B_2$} & $a\blackdiamond c\blackdiamond h\blackdiamond d$ & \text{A},\;\ref{8_12_6} & \\
4 & \text{$C_1$} & $a\blackdiamond c\blackdiamond e\blackdiamond h\blackdiamond d$ & \text{$B_2$},\;\ref{8_12_4} & \\
5 & \text{$C_2$} & $a\blackdiamond c\blackdiamond h\blackdiamond e\blackdiamond d$ & \text{$B_2$},\;\ref{8_12_4} & \\
6 & \text{$D_1$} & $a\blackdiamond h\blackdiamond c\blackdiamond e\blackdiamond f\blackdiamond d$ & \text{$B_1$},\;\ref{8_12_7} & \\
7 & \text{$D_2$} & $a\blackdiamond h\blackdiamond c\blackdiamond e\blackdiamond d\blackdiamond f$ & \text{$B_1$},\;\ref{8_12_7} & \\
8 & \text{$E_1$} & $a\blackdiamond h\blackdiamond b\blackdiamond c\blackdiamond e\blackdiamond f\blackdiamond d$ & \text{$D_1$},\;\ref{8_12_3} & \\
9 & \text{$E_2$} & $a\blackdiamond b\blackdiamond h\blackdiamond c\blackdiamond e\blackdiamond f\blackdiamond d$ & \text{$D_1$},\;\ref{8_12_3} & \\
10 & \text{$E_3$} & $g\blackdiamond b\blackdiamond a\blackdiamond h\blackdiamond c\blackdiamond e\blackdiamond f\blackdiamond d$ & \text{$D_1$},\;\ref{8_12_3},\;\ref{8_12_1} & {\color{red} \ref{8_12_2}},\;{\color{red} \ref{8_12_8}} \\
11 & \text{$F_1$} & $a\blackdiamond h\blackdiamond b\blackdiamond g\blackdiamond c\blackdiamond e\blackdiamond f\blackdiamond d$ & \text{$E_1$},\;\ref{8_12_2} & {\color{blue} \ref{8_12_1}},\;{\color{blue} \ref{8_12_8}} \\
12 & \text{$F_2$} & $a\blackdiamond h\blackdiamond b\blackdiamond c\blackdiamond g\blackdiamond e\blackdiamond f\blackdiamond d$ & \text{$E_1$},\;\ref{8_12_2} & {\color{blue} \ref{8_12_1}},\;{\color{blue} \ref{8_12_8}} \\
13 & \text{$F_3$} & $a\blackdiamond h\blackdiamond b\blackdiamond c\blackdiamond e\blackdiamond g\blackdiamond f\blackdiamond d$ & \text{$E_1$},\;\ref{8_12_2} & {\color{blue} \ref{8_12_1}},\;{\color{blue} \ref{8_12_8}} \\
14 & \text{$G_1$} & $a\blackdiamond b\blackdiamond h\blackdiamond g\blackdiamond c\blackdiamond e\blackdiamond f\blackdiamond d$ & \text{$E_2$},\;\ref{8_12_8} & {\color{blue} \ref{8_12_1}},\;{\color{blue} \ref{8_12_2}} \\
15 & \text{$G_2$} & $a\blackdiamond b\blackdiamond h\blackdiamond c\blackdiamond g\blackdiamond e\blackdiamond f\blackdiamond d$ & \text{$E_2$},\;\ref{8_12_8} & {\color{blue} \ref{8_12_1}},\;{\color{blue} \ref{8_12_2}} \\
16 & \text{$G_3$} & $a\blackdiamond b\blackdiamond h\blackdiamond c\blackdiamond e\blackdiamond g\blackdiamond f\blackdiamond d$ & \text{$E_2$},\;\ref{8_12_8} & {\color{blue} \ref{8_12_1}},\;{\color{blue} \ref{8_12_2}} \\
17 & \text{$H_1$} & $a\blackdiamond h\blackdiamond b\blackdiamond c\blackdiamond e\blackdiamond d\blackdiamond f$ & \text{$D_2$},\;\ref{8_12_3} & \\
18 & \text{$H_2$} & $a\blackdiamond b\blackdiamond h\blackdiamond c\blackdiamond e\blackdiamond d\blackdiamond f$ & \text{$D_2$},\;\ref{8_12_3} & \\
19 & \text{$H_3$} & $g\blackdiamond b\blackdiamond a\blackdiamond h\blackdiamond c\blackdiamond e\blackdiamond d\blackdiamond f$ & \text{$D_2$},\;\ref{8_12_3},\;\ref{8_12_1} & {\color{red} \ref{8_12_2}},\;{\color{red} \ref{8_12_8}} \\
20 & \text{$I_1$} & $a\blackdiamond h\blackdiamond b\blackdiamond g\blackdiamond c\blackdiamond e\blackdiamond d\blackdiamond f$ & \text{$H_1$},\;\ref{8_12_2} & {\color{blue} \ref{8_12_1}},\;{\color{blue} \ref{8_12_8}} \\
21 & \text{$I_2$} & $a\blackdiamond h\blackdiamond b\blackdiamond c\blackdiamond g\blackdiamond e\blackdiamond d\blackdiamond f$ & \text{$H_1$},\;\ref{8_12_2} & {\color{blue} \ref{8_12_1}},\;{\color{blue} \ref{8_12_8}} \\
22 & \text{$I_3$} & $a\blackdiamond h\blackdiamond b\blackdiamond c\blackdiamond e\blackdiamond g\blackdiamond d\blackdiamond f$ & \text{$H_1$},\;\ref{8_12_2} & {\color{blue} \ref{8_12_1}},\;{\color{blue} \ref{8_12_8}} \\
23 & \text{$I_4$} & $a\blackdiamond h\blackdiamond b\blackdiamond c\blackdiamond e\blackdiamond d\blackdiamond g\blackdiamond f$ & \text{$H_1$},\;\ref{8_12_2} & {\color{blue} \ref{8_12_1}},\;{\color{blue} \ref{8_12_8}} \\
24 & \text{$J_1$} & $a\blackdiamond b\blackdiamond h\blackdiamond g\blackdiamond c\blackdiamond e\blackdiamond d\blackdiamond f$ & \text{$H_2$},\;\ref{8_12_8} & {\color{blue} \ref{8_12_1}},\;{\color{blue} \ref{8_12_2}} \\
25 & \text{$J_2$} & $a\blackdiamond b\blackdiamond h\blackdiamond c\blackdiamond g\blackdiamond e\blackdiamond d\blackdiamond f$ & \text{$H_2$},\;\ref{8_12_8} & {\color{blue} \ref{8_12_1}},\;{\color{blue} \ref{8_12_2}} \\
26 & \text{$J_3$} & $a\blackdiamond b\blackdiamond h\blackdiamond c\blackdiamond e\blackdiamond g\blackdiamond d\blackdiamond f$ & \text{$H_2$},\;\ref{8_12_8} & {\color{blue} \ref{8_12_1}},\;{\color{blue} \ref{8_12_2}} \\
27 & \text{$J_4$} & $a\blackdiamond b\blackdiamond h\blackdiamond c\blackdiamond e\blackdiamond d\blackdiamond g\blackdiamond f$ & \text{$H_2$},\;\ref{8_12_8} & {\color{blue} \ref{8_12_1}},\;{\color{blue} \ref{8_12_2}} \\
28 & \text{$K_1$} & $a\blackdiamond b\blackdiamond c\blackdiamond e\blackdiamond h\blackdiamond d$ & \text{$C_1$},\;\ref{8_12_3} & \\
29 & \text{$K_2$} & $f\blackdiamond g\blackdiamond b\blackdiamond a\blackdiamond c\blackdiamond e\blackdiamond h\blackdiamond d$ & \text{$C_1$},\;\ref{8_12_3},\;\ref{8_12_1},\;\ref{8_12_2} & {\color{blue} \ref{8_12_7}},\;{\color{blue} \ref{8_12_8}} \\
30 & \text{$L_1$} & $a\blackdiamond b\blackdiamond c\blackdiamond f\blackdiamond e\blackdiamond h\blackdiamond d$ & \text{$K_1$},\;\ref{8_12_7} &  \\
31 & \text{$L_2$} & $a\blackdiamond b\blackdiamond g\blackdiamond f\blackdiamond c\blackdiamond e\blackdiamond h\blackdiamond d$ & \text{$K_1$},\;\ref{8_12_7},\;\ref{8_12_2} & {\color{blue} \ref{8_12_1}},\;{\color{red} \ref{8_12_8}} \\
32 & \text{$L_3$} & $a\blackdiamond f\blackdiamond g\blackdiamond b\blackdiamond c\blackdiamond e\blackdiamond h\blackdiamond d$ & \text{$K_1$},\;\ref{8_12_7},\;\ref{8_12_2} & {\color{red} \ref{8_12_1}},\;{\color{blue} \ref{8_12_8}} \\
33 & \text{$L_4$} & $f\blackdiamond a\blackdiamond b\blackdiamond c\blackdiamond e\blackdiamond h\blackdiamond d$ & \text{$K_1$},\;\ref{8_12_7} &  \\
34 & \text{$M_1$} & $a\blackdiamond b\blackdiamond g\blackdiamond c\blackdiamond f\blackdiamond e\blackdiamond h\blackdiamond d$ & \text{$L_1$},\;\ref{8_12_2} & {\color{blue} \ref{8_12_1}},\;{\color{red} \ref{8_12_8}} \\
35 & \text{$M_2$} & $a\blackdiamond b\blackdiamond c\blackdiamond g\blackdiamond f\blackdiamond e\blackdiamond h\blackdiamond d$ & \text{$L_1$},\;\ref{8_12_2} & {\color{blue} \ref{8_12_1}},\;{\color{red} \ref{8_12_8}} \\
36 & \text{$N_1$} & $f\blackdiamond g\blackdiamond a\blackdiamond b\blackdiamond c\blackdiamond e\blackdiamond h\blackdiamond d$ & \text{$L_4$},\;\ref{8_12_2} & {\color{red} \ref{8_12_1}},\;{\color{blue} \ref{8_12_8}} \\
37 & \text{$N_2$} & $f\blackdiamond a\blackdiamond g\blackdiamond b\blackdiamond c\blackdiamond e\blackdiamond h\blackdiamond d$ & \text{$L_4$},\;\ref{8_12_2} & {\color{red} \ref{8_12_1}},\;{\color{blue} \ref{8_12_8}} \\
38 & \text{$O_1$} & $a\blackdiamond b\blackdiamond c\blackdiamond h\blackdiamond e\blackdiamond d$ & \text{$C_2$},\;\ref{8_12_3} &  \\
39 & \text{$O_2$} & $f\blackdiamond g\blackdiamond b\blackdiamond a\blackdiamond c\blackdiamond h\blackdiamond e\blackdiamond d$ & \text{$C_2$},\;\ref{8_12_3},\;\ref{8_12_1},\;\ref{8_12_2} & {\color{red} \ref{8_12_7}},\;{\color{blue} \ref{8_12_8}} \\
40 & \text{$P_1$} & $a\blackdiamond b\blackdiamond c\blackdiamond h\blackdiamond e\blackdiamond f\blackdiamond d$ & \text{$O_1$},\;\ref{8_12_7} &  \\
41 & \text{$P_2$} & $a\blackdiamond b\blackdiamond c\blackdiamond h\blackdiamond e\blackdiamond d\blackdiamond f$ & \text{$O_1$},\;\ref{8_12_7} &  \\
42 & \text{$Q_1$} & $a\blackdiamond b\blackdiamond c\blackdiamond h\blackdiamond g\blackdiamond e\blackdiamond f\blackdiamond d$ & \text{$P_1$},\;\ref{8_12_8} & {\color{blue} \ref{8_12_1}},\;{\color{blue} \ref{8_12_2}}  \\
43 & \text{$Q_2$} & $a\blackdiamond b\blackdiamond c\blackdiamond h\blackdiamond e\blackdiamond g\blackdiamond f\blackdiamond d$ & \text{$P_1$},\;\ref{8_12_8} & {\color{blue} \ref{8_12_1}},\;{\color{blue} \ref{8_12_2}}  \\
44 & \text{$R_1$} & $a\blackdiamond b\blackdiamond c\blackdiamond h\blackdiamond g\blackdiamond e\blackdiamond d\blackdiamond f$ & \text{$P_2$},\;\ref{8_12_8} & {\color{blue} \ref{8_12_1}},\;{\color{blue} \ref{8_12_2}}  \\
45 & \text{$R_2$} & $a\blackdiamond b\blackdiamond c\blackdiamond h\blackdiamond e\blackdiamond g\blackdiamond d\blackdiamond f$ & \text{$P_2$},\;\ref{8_12_8} & {\color{blue} \ref{8_12_1}},\;{\color{blue} \ref{8_12_2}}  \\
46 & \text{$R_3$} & $a\blackdiamond b\blackdiamond c\blackdiamond h\blackdiamond e\blackdiamond d\blackdiamond g\blackdiamond f$ & \text{$P_2$},\;\ref{8_12_8} & {\color{blue} \ref{8_12_1}},\;{\color{blue} \ref{8_12_2}}  \\

\hline
\end{longtblr}

Linear orders on the generating set of $Q(8_{12})$ as in \text{$F_1$}, \text{$F_2$}, \text{$F_3$}, \text{$G_1$}, \text{$G_2$}, \text{$G_3$}, \text{$I_1$}, \text{$I_2$}, \text{$I_3$}, \text{$I_4$}, \text{$J_1$}, \text{$J_2$}, \text{$J_3$}, \text{$J_4$}, \text{$K_2$}, \text{$Q_1$}, \text{$Q_2$}, \text{$R_1$}, \text{$R_2$} and \text{$R_3$} in Table \ref{tbl8_12} could be extendable to biorders on $Q(8_{12})$.

\subsection{Knot $8_{13}$}\label{8_13}

\begin{figure}[H]
	\centering
	\includegraphics[scale=0.3]{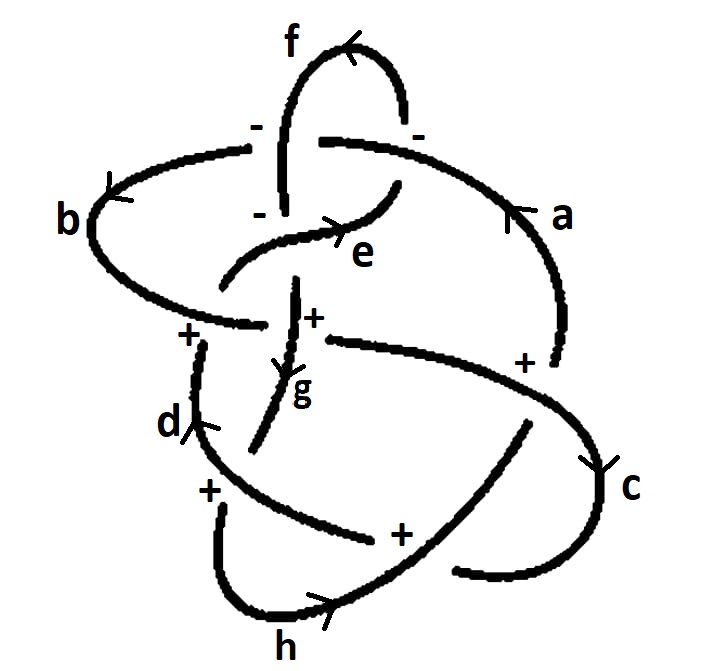}
	\caption{Knot $8_{13}$}
    \label{fig8_13}
\end{figure}

Consider the diagram of the knot $8_{13}$ as in Figure \ref{fig8_13}. The quandle $Q(8_{13})$ is generated by labelings of arcs the diagram, and the relations are given by

\begin{align}
	b*f&=a&b\blackdiamond a\blackdiamond f \label{8_13_1}\\
	f*a&=e&f\blackdiamond e\blackdiamond a \label{8_13_2}\\
	g*e&=f&g\blackdiamond f\blackdiamond e \label{8_13_3}\\
	d*b&=e&d\blackdiamond e\blackdiamond b \label{8_13_4}\\
	b*g&=c&b\blackdiamond c\blackdiamond g \label{8_13_5}\\
	h*c&=a&h\blackdiamond a\blackdiamond c \label{8_13_6}\\
	c*h&=d&c\blackdiamond d\blackdiamond h \label{8_13_7}\\
	g*d&=h&g\blackdiamond h\blackdiamond d \label{8_13_8}
\end{align}

We have the table as follows.

\begin{longtblr}[caption = {Orders that could be extendable to biorders on $Q(8_{13})$}, label = {tbl8_13}]
{colspec = {|c|c|l|l|l|}, rowhead = 2}
            \hline
			\text{Sr.} & \text{Label} & \text{Expression} & \text{Obtained by} & \text{{\color{red} Contradiction to} $\backslash$}\\
			\text{No.} &  & &  & \text{\color{blue} Compatible with}\\[2mm]
			\hline
1& \text{A}& $b\blackdiamond a\blackdiamond e\blackdiamond f\blackdiamond g$ & \ref{8_13_1},\;\ref{8_13_2},\;\ref{8_13_3} & \\
2& \text{$B_1$}& $b\blackdiamond a\blackdiamond e\blackdiamond d\blackdiamond f\blackdiamond g$ & \text{A},\;\ref{8_13_4} & \\
3& \text{$B_2$}& $b\blackdiamond a\blackdiamond e\blackdiamond f\blackdiamond d\blackdiamond h\blackdiamond g$ & \text{A},\;\ref{8_13_4},\;\ref{8_13_8} & \\
4& \text{$B_3$}& $b\blackdiamond a\blackdiamond e\blackdiamond f\blackdiamond g\blackdiamond h\blackdiamond d$ & \text{A},\;\ref{8_13_4},\;\ref{8_13_8} & \\
5& \text{$C_1$}& $b\blackdiamond a\blackdiamond e\blackdiamond d\blackdiamond h\blackdiamond f\blackdiamond g$ & \text{$B_1$},\;\ref{8_13_8} & \\
6& \text{$C_2$}& $b\blackdiamond a\blackdiamond e\blackdiamond d\blackdiamond f\blackdiamond h\blackdiamond g$ & \text{$B_1$},\;\ref{8_13_8} & \\
7& \text{$D_1$}& $b\blackdiamond c\blackdiamond a\blackdiamond e\blackdiamond d\blackdiamond h\blackdiamond f\blackdiamond g$ & \text{$C_1$},\;\ref{8_13_6} &{\color{blue} \ref{8_13_5},\;\ref{8_13_7}} \\
8& \text{$D_2$}& $c\blackdiamond b\blackdiamond a\blackdiamond e\blackdiamond d\blackdiamond h\blackdiamond f\blackdiamond g$ & \text{$C_1$},\;\ref{8_13_6} &{\color{red} \ref{8_13_5},}\;{\color{blue} \ref{8_13_7}} \\
9& \text{$E_1$}& $b\blackdiamond c\blackdiamond a\blackdiamond e\blackdiamond d\blackdiamond f\blackdiamond h\blackdiamond g$ & \text{$C_2$},\;\ref{8_13_6} &{\color{blue} \ref{8_13_5},}\;{\color{blue} \ref{8_13_7}} \\
10& \text{$E_2$}& $c\blackdiamond b\blackdiamond a\blackdiamond e\blackdiamond d\blackdiamond f\blackdiamond h\blackdiamond g$ & \text{$C_2$},\;\ref{8_13_6} &{\color{red} \ref{8_13_5},}\;{\color{blue} \ref{8_13_7}} \\
11& \text{$F_1$}& $b\blackdiamond c\blackdiamond a\blackdiamond e\blackdiamond f\blackdiamond d\blackdiamond h\blackdiamond g$ & \text{$B_2$},\;\ref{8_13_6} &{\color{blue} \ref{8_13_5},}\;{\color{blue} \ref{8_13_7}} \\
12& \text{$F_2$}& $c\blackdiamond b\blackdiamond a\blackdiamond e\blackdiamond f\blackdiamond d\blackdiamond h\blackdiamond g$ & \text{$B_2$},\;\ref{8_13_6} &{\color{red} \ref{8_13_5},}\;{\color{blue} \ref{8_13_7}} \\
13& \text{G}& $b\blackdiamond a\blackdiamond e\blackdiamond f\blackdiamond g\blackdiamond h\blackdiamond d\blackdiamond c$ & \text{$B_3$},\;\ref{8_13_7} &{\color{red} \ref{8_13_5},}\;{\color{red} \ref{8_13_6}} \\
	\hline
\end{longtblr}

Linear orders on the generating set of $Q(8_{13})$ as in \text{$D_1$}, \text{$E_1$} and \text{$F_1$} in Table \ref{tbl8_13} could be extendable to biorders on $Q(8_{13})$.

\subsection{Knot $8_{14}$}\label{8_14}

\begin{figure}[H]
	\centering
	\includegraphics[scale=0.3]{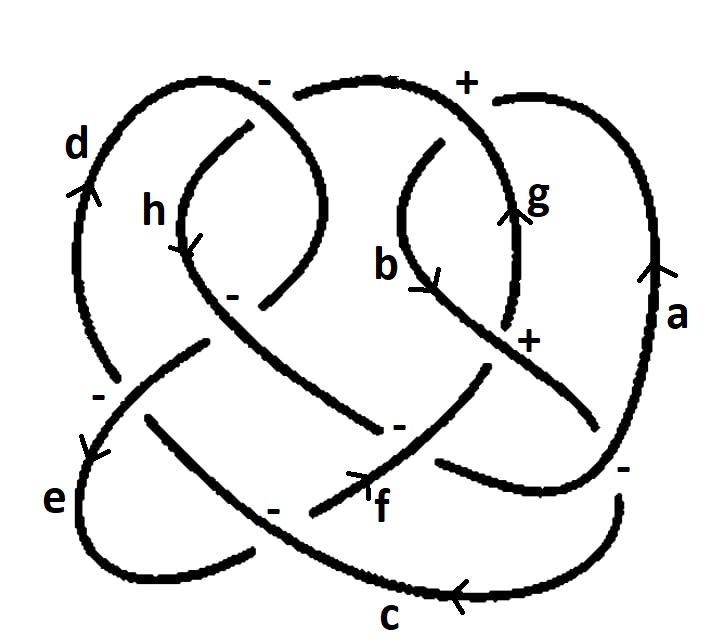}
	\caption{Knot $8_{14}$}
    \label{fig8_14}
\end{figure}

Consider the diagram of the knot $8_{14}$ as in Figure \ref{fig8_14}. The quandle $Q(8_{14})$ is generated by labelings of arcs of the diagram, and the relations are given by
\begin{align}
	a*g&=b&a\blackdiamond b\blackdiamond g \label{8_14_1}\\
	f*b&=g&f\blackdiamond g\blackdiamond b \label{8_14_2}\\
	c*a&=b&c\blackdiamond b\blackdiamond a \label{8_14_3}\\
	a*f&=h&a\blackdiamond h\blackdiamond f \label{8_14_4}\\
	f*c&=e&f\blackdiamond e\blackdiamond c \label{8_14_5}\\
	e*h&=d&e\blackdiamond d\blackdiamond h \label{8_14_6}\\
	d*e&=c&d\blackdiamond c\blackdiamond e \label{8_14_7}\\
	h*d&=g&h\blackdiamond g\blackdiamond d \label{8_14_8}
\end{align}

The determinant of knot $8_{14}$ is $31$. By \cite[Proposition 4.4]{ms}, the generators $a$, $b$, $c$, $d$, $e$, $f$, $g$ and $h$ are all pairwise distinct. Applying the relations above ((\ref{8_14_1}) to (\ref{8_14_8})), we have the table as follows.

\begin{longtblr}[caption = {Orders that could be extendable to biorders on $Q(8_{14})$}, label = {tbl8_14}]
{colspec = {|c|c|l|l|l|}, rowhead = 2}
            \hline
			\text{Sr.} & \text{Label} & \text{Expression} & \text{Obtained by} & \text{{\color{red} Contradiction to} $\backslash$}\\
			\text{No.} &  & &  & \text{\color{blue} Compatible with}\\[2mm]
			\hline
1& \text{A}& $f\blackdiamond e\blackdiamond c\blackdiamond d$ & \ref{8_14_5},\;\ref{8_14_7} & \\
2& \text{B}& $f\blackdiamond e\blackdiamond c\blackdiamond d\blackdiamond h\blackdiamond a$ & \text{A},\;\ref{8_14_6},\;\ref{8_14_4}& \\
3& \text{C}& $f\blackdiamond e\blackdiamond c\blackdiamond d\blackdiamond g\blackdiamond h\blackdiamond a$ & \text{B},\;\ref{8_14_8}& \\
4& \text{$D_1$}& $f\blackdiamond e\blackdiamond c\blackdiamond d\blackdiamond g\blackdiamond b\blackdiamond h\blackdiamond a$ & \text{C},\;\ref{8_14_1}& {\color{blue} \ref{8_14_2},\;\ref{8_14_3}} \\
5& \text{$D_2$}& $f\blackdiamond e\blackdiamond c\blackdiamond d\blackdiamond g\blackdiamond h\blackdiamond b\blackdiamond a$ & \text{C},\;\ref{8_14_1}& {\color{blue} \ref{8_14_2},\;\ref{8_14_3}} \\
\hline
\end{longtblr}

In Table \ref{tbl8_14}, linear orders \text{$D_1$} and \text{$D_2$} are compatible with all the relations. Therefore, these orders could be extendable to biorders on $Q(8_{14})$.

\subsection{Knot $8_{15}$}\label{8_15}

The knot $8_{15}$ is prime, alternating and positive (or negative). It is a pretzel knot and hence Montesinos one, see \cite[Theorem 18]{dm}. By \cite[Corollary 6.7]{rss}, the quandle of this knot is not biorderable.

\subsection{Knot $8_{16}$}\label{8_16}

\begin{figure}[H]
	\centering
	\includegraphics[scale=0.3]{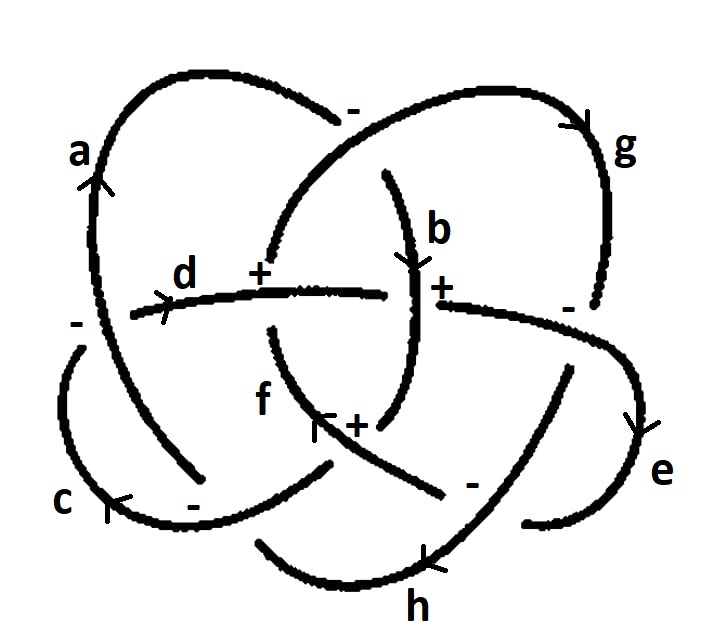}
	\caption{Knot $8_{16}$}
    \label{fig8_16}
\end{figure}

Consider the diagram of the knot $8_{16}$ as in Figure \ref{fig8_16}. The quandle $Q(8_{16})$ is generated by labelings of arcs of the diagram, and the relations are given by
\begin{align}
	b*g&=a&b\blackdiamond a\blackdiamond g \label{8_16_1}\\
	d*a&=c&d\blackdiamond c\blackdiamond a \label{8_16_2}\\
	a*c&=h&a\blackdiamond h\blackdiamond c \label{8_16_3}\\
	f*h&=e&f\blackdiamond e\blackdiamond h \label{8_16_4}\\
	h*e&=g&h\blackdiamond g\blackdiamond e \label{8_16_5}\\
	f*d&=g&f\blackdiamond g\blackdiamond d \label{8_16_6}\\
	d*b&=e&d\blackdiamond e\blackdiamond b \label{8_16_7}\\
	b*f&=c&b\blackdiamond c\blackdiamond f \label{8_16_8}
\end{align}

Since there is a quandle homomorphism  $\phi: Q(8_{16})\to\text{Core}(\mathbb{Z}_{35})$ given by $a\mapsto 1, b\mapsto 21, c\mapsto 0, d\mapsto 2, e\mapsto 5, f\mapsto28, g\mapsto11$ and $h\mapsto34$, the generators $a$, $b$, $c$, $d$, $e$, $f$, $g$ and $h$ are all pairwise distinct. Applying the relations above ((\ref{8_16_1}) to (\ref{8_16_8})), we can make the table as follows.

\begin{longtblr}[caption = {Orders that could be extendable to biorders on $Q(8_{16})$}, label = {tbl8_16}]
{colspec = {|c|c|l|l|l|}, rowhead = 2}
            \hline
			\text{Sr.} & \text{Label} & \text{Expression} & \text{Obtained by} & \text{{\color{red} Contradiction to} $\backslash$}\\
			\text{No.} &  & &  & \text{\color{blue} Compatible with}\\[2mm]
			\hline
1& \text{A} & $h\blackdiamond g\blackdiamond e\blackdiamond f$ & \ref{8_16_5},\;\ref{8_16_4}& \\
2& \text{$B_1$} & $h\blackdiamond d\blackdiamond g\blackdiamond e\blackdiamond f$	& \text{A},\;\ref{8_16_6} & \\
3& \text{$B_2$} & $d\blackdiamond h\blackdiamond g\blackdiamond e\blackdiamond f$	& \text{A},\;\ref{8_16_6} & \\
4& \text{$C_1$} & $a\blackdiamond h\blackdiamond d\blackdiamond g\blackdiamond e\blackdiamond b\blackdiamond c\blackdiamond f$ & \text{$B_1$},\;\ref{8_16_7},\;\ref{8_16_8},\;\ref{8_16_3} & {\color{red} \ref{8_16_1},\;\ref{8_16_2}}\\
5& \text{$C_2$} & $a\blackdiamond h\blackdiamond d\blackdiamond g\blackdiamond e\blackdiamond f\blackdiamond c\blackdiamond b$ & \text{$B_1$},\;\ref{8_16_7},\;\ref{8_16_8},\;\ref{8_16_3} & {\color{red} \ref{8_16_1},\;\ref{8_16_2}}\\
6& \text{$D_1$} & $d\blackdiamond h\blackdiamond g\blackdiamond e\blackdiamond b\blackdiamond c\blackdiamond f$ & \text{$B_2$},\;\ref{8_16_7},\;\ref{8_16_8}&\\
7& \text{$D_2$} & $d\blackdiamond h\blackdiamond g\blackdiamond e\blackdiamond f\blackdiamond c\blackdiamond b$ & \text{$B_2$},\;\ref{8_16_7},\;\ref{8_16_8}&\\
8& \text{$E_1$} & $d\blackdiamond a\blackdiamond h\blackdiamond g\blackdiamond e\blackdiamond b\blackdiamond c\blackdiamond f$ & \text{$D_1$},\;\ref{8_16_3} & {\color{red} \ref{8_16_1},\;\ref{8_16_2}}\\
9& \text{$E_2$} & $a\blackdiamond d\blackdiamond h\blackdiamond g\blackdiamond e\blackdiamond b\blackdiamond c\blackdiamond f$ & \text{$D_1$},\;\ref{8_16_3} & {\color{red} \ref{8_16_1},\;\ref{8_16_2}}\\
10& \text{$F_1$} & $d\blackdiamond a\blackdiamond h\blackdiamond g\blackdiamond e\blackdiamond f\blackdiamond c\blackdiamond b$ & \text{$D_2$},\;\ref{8_16_3} & {\color{red} \ref{8_16_1},\;\ref{8_16_2}}\\
11& \text{$F_2$} & $a\blackdiamond d\blackdiamond h\blackdiamond g\blackdiamond e\blackdiamond f\blackdiamond c\blackdiamond b$ & \text{$D_2$},\;\ref{8_16_3} & {\color{red} \ref{8_16_1},\;\ref{8_16_2}}\\
\hline
\end{longtblr}

In Table \ref{tbl8_16}, linear orders \text{$C_1$}, \text{$C_2$}, \text{$E_1$}, \text{$E_2$}, \text{$F_1$} and \text{$F_2$} contradict (\ref{8_16_1}) and (\ref{8_16_2}). Therefore, we conclude that $Q(8_{16})$ can not be biorderable.

\subsection{Knot $8_{17}$}\label{8_17}

\begin{figure}[H]
	\centering
	\includegraphics[scale=0.3]{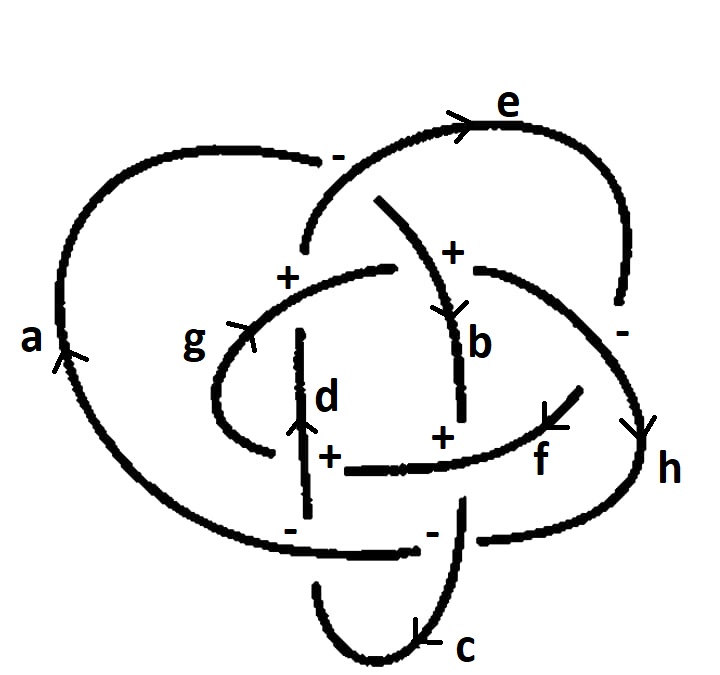}
	\caption{Knot $8_{17}$}
    \label{fig8_17}
\end{figure}

Consider the diagram of the knot $8_{17}$ as in Figure \ref{fig8_17}. The quandle $Q(8_{17})$ is generated by labelings of arcs of the diagram, and the relations are given by
\begin{align}
	b*e&=a&b\blackdiamond a\blackdiamond e \label{8_17_1}\\
	f*h&=e&f\blackdiamond e\blackdiamond h \label{8_17_2}\\
	a*c&=h&a\blackdiamond h\blackdiamond c \label{8_17_3}\\
	d*a&=c&d\blackdiamond c\blackdiamond a \label{8_17_4}\\
	f*d&=g&f\blackdiamond g\blackdiamond d \label{8_17_5}\\
	d*g&=e&d\blackdiamond e\blackdiamond g \label{8_17_6}\\
	g*b&=h&g\blackdiamond h\blackdiamond b \label{8_17_7}\\
	b*f&=c&b\blackdiamond c\blackdiamond f \label{8_17_8}
\end{align}

The determinant of knot $8_{17}$ is $37$. By \cite[Proposition 4.4]{ms}, the generators $a$, $b$, $c$, $d$, $e$, $f$, $g$ and $h$ are all pairwise distinct. Applying the relations above ((\ref{8_17_1}) to (\ref{8_17_8})), we have the table as follows.

\begin{longtblr}[caption = {Orders that could be extendable to biorders on $Q(8_{17})$}, label = {tbl8_17}]
{colspec = {|c|c|l|l|l|}, rowhead = 2}
            \hline
			\text{Sr.} & \text{Label} & \text{Expression} & \text{Obtained by} & \text{{\color{red} Contradiction to} $\backslash$}\\
			\text{No.} &  & &  & \text{\color{blue} Compatible with}\\[2mm]
			\hline
1& \text{A} & $f\blackdiamond g\blackdiamond e\blackdiamond d$ & \ref{8_17_5},\;\ref{8_17_6}& \\
2& \text{$B_1$} & $f\blackdiamond g\blackdiamond e\blackdiamond h\blackdiamond d$ & \text{A},\;\ref{8_17_2}&\\
3& \text{$B_2$} & $f\blackdiamond g\blackdiamond e\blackdiamond d\blackdiamond h\blackdiamond b$ & \text{A},\;\ref{8_17_2},\;\ref{8_17_7}&\\
4& \text{$C_1$} & $f\blackdiamond g\blackdiamond e\blackdiamond h\blackdiamond b\blackdiamond d$ & \text{$B_1$},\;\ref{8_17_7}&\\
5& \text{$C_2$} & $f\blackdiamond g\blackdiamond e\blackdiamond h\blackdiamond d\blackdiamond b$ & \text{$B_1$},\;\ref{8_17_7}&\\
6& \text{$D_1$} & $f\blackdiamond g\blackdiamond e\blackdiamond a\blackdiamond h\blackdiamond b\blackdiamond d$ & \text{$C_1$},\;\ref{8_17_1}&\\
7& \text{$D_2$} & $f\blackdiamond g\blackdiamond e\blackdiamond h\blackdiamond a\blackdiamond b\blackdiamond d$ & \text{$C_1$},\;\ref{8_17_1}&\\
8& \text{$E_1$} & $f\blackdiamond g\blackdiamond e\blackdiamond a\blackdiamond c\blackdiamond h\blackdiamond b\blackdiamond d$ & \text{$D_1$},\;\ref{8_17_4}& {\color{red} \ref{8_17_3},}\;{\color{blue} \ref{8_17_8}}\\
9& \text{$E_2$} & $f\blackdiamond g\blackdiamond e\blackdiamond a\blackdiamond h\blackdiamond c\blackdiamond b\blackdiamond d$ & \text{$D_1$},\;\ref{8_17_4}& {\color{blue} \ref{8_17_3},}\;{\color{blue} \ref{8_17_8}}\\
10& \text{$E_3$} & $f\blackdiamond g\blackdiamond e\blackdiamond a\blackdiamond h\blackdiamond b\blackdiamond c\blackdiamond d$ & \text{$D_1$},\;\ref{8_17_4}& {\color{blue} \ref{8_17_3},}\;{\color{red} \ref{8_17_8}}\\
11& \text{$F_1$} & $f\blackdiamond g\blackdiamond e\blackdiamond h\blackdiamond a\blackdiamond c\blackdiamond b\blackdiamond d$ & \text{$D_2$},\;\ref{8_17_4}& {\color{red} \ref{8_17_3},}\;{\color{blue} \ref{8_17_8}}\\
12& \text{$F_2$} & $f\blackdiamond g\blackdiamond e\blackdiamond h\blackdiamond a\blackdiamond b\blackdiamond c\blackdiamond d$ & \text{$D_2$},\;\ref{8_17_4}& {\color{red} \ref{8_17_3},}\;{\color{red} \ref{8_17_8}}\\
13& \text{$G_1$} & $f\blackdiamond g\blackdiamond e\blackdiamond a\blackdiamond h\blackdiamond d\blackdiamond b$ & \text{$C_2$},\;\ref{8_17_1}&\\
14& \text{$G_2$} & $f\blackdiamond g\blackdiamond e\blackdiamond h\blackdiamond a\blackdiamond c\blackdiamond d\blackdiamond b$ & \text{$C_2$},\;\ref{8_17_1},\;\ref{8_17_4}& {\color{red} \ref{8_17_3},}\;{\color{blue} \ref{8_17_8}}\\
15& \text{$G_3$} & $f\blackdiamond g\blackdiamond e\blackdiamond h\blackdiamond d\blackdiamond c\blackdiamond a\blackdiamond b$ & \text{$C_2$},\;\ref{8_17_1},\;\ref{8_17_4}& {\color{red} \ref{8_17_3},}\;{\color{blue} \ref{8_17_8}}\\
16& \text{$H_1$} & $f\blackdiamond g\blackdiamond e\blackdiamond a\blackdiamond c\blackdiamond h\blackdiamond d\blackdiamond b$ & \text{$G_1$},\;\ref{8_17_4}& {\color{red} \ref{8_17_3},}\;{\color{blue} \ref{8_17_8}}\\
17& \text{$H_2$} & $f\blackdiamond g\blackdiamond e\blackdiamond a\blackdiamond h\blackdiamond c\blackdiamond d\blackdiamond b$ & \text{$G_1$},\;\ref{8_17_4}& {\color{blue} \ref{8_17_3},}\;{\color{blue} \ref{8_17_8}}\\
18& \text{$I_1$} & $f\blackdiamond g\blackdiamond e\blackdiamond a\blackdiamond c\blackdiamond d\blackdiamond h\blackdiamond b$ & \text{$B_2$},\;\ref{8_17_1},\;\ref{8_17_4}& {\color{red} \ref{8_17_3},}\;{\color{blue} \ref{8_17_8}}\\
19& \text{$I_2$} & $f\blackdiamond g\blackdiamond e\blackdiamond d\blackdiamond c\blackdiamond a\blackdiamond h\blackdiamond b$ & \text{$B_2$},\;\ref{8_17_1},\;\ref{8_17_4}& {\color{red} \ref{8_17_3},}\;{\color{blue} \ref{8_17_8}}\\
20& \text{$I_3$} & $f\blackdiamond g\blackdiamond e\blackdiamond d\blackdiamond h\blackdiamond a\blackdiamond b$ & \text{$B_2$},\;\ref{8_17_1}& \\
21& \text{$J_1$} & $f\blackdiamond g\blackdiamond e\blackdiamond d\blackdiamond c\blackdiamond h\blackdiamond a\blackdiamond b$ & \text{$I_3$},\;\ref{8_17_4} & {\color{blue} \ref{8_17_3},\;\ref{8_17_8}}\\
22& \text{$J_2$} & $f\blackdiamond g\blackdiamond e\blackdiamond d\blackdiamond h\blackdiamond c\blackdiamond a\blackdiamond b$ & \text{$I_3$},\;\ref{8_17_4} & {\color{red} \ref{8_17_3},}\;{\color{blue} \ref{8_17_8}}\\

\hline
\end{longtblr}

Linear orders on the generating set of $Q(8_{17})$ as in \text{$E_2$}, \text{$H_2$} and \text{$J_1$} in Table \ref{tbl8_17} could be extendable to biorders on $Q(8_{17})$.

\subsection{Knot $8_{18}$}\label{8_18}

\begin{figure}[H]
	\centering
	\includegraphics[scale=0.3]{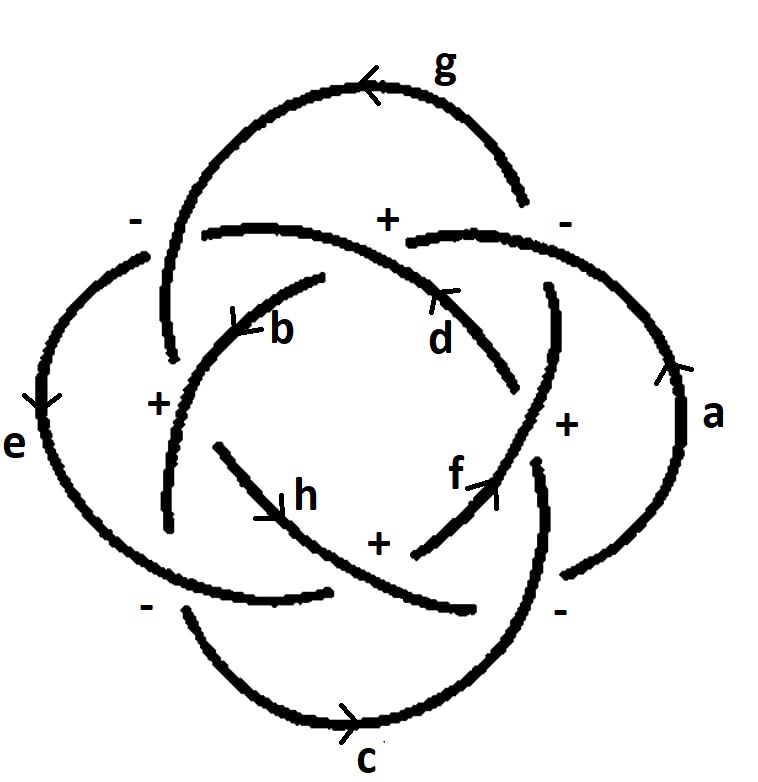}
	\caption{Knot $8_{18}$}
	\label{fig8_18}
\end{figure}

Consider the diagram of the knot $8_{18}$ as in Figure \ref{fig8_18}. The quandle $Q(8_{18})$ is generated by labelings of arcs of the diagram, and the relations are given by
\begin{align}
	a*c&=h&a\blackdiamond h\blackdiamond c \label{8_18_1}\\
	g*a&=f&g\blackdiamond f\blackdiamond a \label{8_18_2}\\
	e*g&=d&e\blackdiamond d\blackdiamond g \label{8_18_3}\\
	c*e&=b&c\blackdiamond b\blackdiamond e \label{8_18_4}\\
	e*h&=f&e\blackdiamond f\blackdiamond h \label{8_18_5}\\
	c*f&=d&c\blackdiamond d\blackdiamond f \label{8_18_6}\\
	a*d&=b&a\blackdiamond b\blackdiamond d \label{8_18_7}\\
	g*b&=h&g\blackdiamond h\blackdiamond b \label{8_18_8}
\end{align}

There is a quandle homomorphism  $\phi: Q(8_{18})\to\text{Core}(\mathbb{Z}_3)$ given by  $a\mapsto 0, b\mapsto 0, c\mapsto 1, d\mapsto 0, e\mapsto 2, f\mapsto2, g\mapsto1$, $h\mapsto2$, and another one $\psi: Q(8_{18})\to\text{Core}(\mathbb{Z}_3)$  given by $a\mapsto 1, b\mapsto 2, c\mapsto 0, d\mapsto 0, e\mapsto 1, f\mapsto0, g\mapsto2$, $h\mapsto2$. Both homomorphisms together gives us that the generators $a$, $b$, $c$, $d$, $e$, $f$, $g$ and $h$ are all pairwise distinct. Applying the relations above ((\ref{8_18_1}) to (\ref{8_18_8})), we have the table as follows.

\begin{longtblr}[caption = {Orders that could be extendable to biorders on $Q(8_{18})$}, label = {tbl8_18}]
	{colspec = {|c|c|l|l|l|}, rowhead = 2}
	\hline
	\text{Sr.} & \text{Label} & \text{Expression} & \text{Obtained by} & \text{{\color{red} Contradiction to} $\backslash$}\\
	\text{No.} &  & &  & \text{\color{blue} Compatible with}\\[2mm]
	\hline
	1 & A &$ e\blackdiamond f\blackdiamond h$ & \ref{8_18_5} & \\
	2 & \text{$B_1$} &$ g\blackdiamond e\blackdiamond f\blackdiamond h $& A & \\
	3 & \text{$B_2$} &$ e\blackdiamond g\blackdiamond f\blackdiamond h$ & A & \\
    4 & \text{$B_3$} &$ e\blackdiamond f\blackdiamond g\blackdiamond h$ & A & \\
    5 & \text{$B_4$} &$ e\blackdiamond f\blackdiamond h\blackdiamond g$ & A & \\
    6 & \text{C} &$ g\blackdiamond d\blackdiamond e\blackdiamond f\blackdiamond h\blackdiamond b\blackdiamond a$ & \text{$B_1$},\;\ref{8_18_3},\;\ref{8_18_8},\;\ref{8_18_7}& \\
   7 & \text{$D_1$} &$ g\blackdiamond d\blackdiamond e\blackdiamond f\blackdiamond h\blackdiamond b\blackdiamond c\blackdiamond a$ & \text{C},\;\ref{8_18_4} & {\color{red} \ref{8_18_6}} \\
   8 & \text{$D_2$} &$ g\blackdiamond d\blackdiamond e\blackdiamond f\blackdiamond h\blackdiamond b\blackdiamond a\blackdiamond c$ & \text{C},\;\ref{8_18_4} & {\color{red} \ref{8_18_6}} \\
   9 & \text{E} &$ e\blackdiamond d\blackdiamond g\blackdiamond f\blackdiamond h\blackdiamond b\blackdiamond a$ & \text{$B_2$},\;\ref{8_18_3},\;\ref{8_18_8},\;\ref{8_18_7}\\
  10 & \text{$F_1$} &$ e\blackdiamond d\blackdiamond g\blackdiamond f\blackdiamond h\blackdiamond b\blackdiamond c\blackdiamond a $& \text{E},\;\ref{8_18_4} & {\color{red} \ref{8_18_1}} \\
 11 & \text{$F_2$} &$ e\blackdiamond d\blackdiamond g\blackdiamond f\blackdiamond h\blackdiamond b\blackdiamond a\blackdiamond c$ & \text{E},\;\ref{8_18_4} & {\color{red} \ref{8_18_1}} \\
12 & \text{G} &$ e\blackdiamond f\blackdiamond g\blackdiamond h\blackdiamond b\blackdiamond c$ & \text{$B_3$},\;\ref{8_18_8},\;\ref{8_18_4}& \\
13 & \text{H} &$ e\blackdiamond f\blackdiamond d\blackdiamond g\blackdiamond h\blackdiamond b\blackdiamond c$ & \text{G},\;\ref{8_18_3},\;\ref{8_18_6}& \\
14 & \text{$I_1$} &$ a\blackdiamond e\blackdiamond f\blackdiamond d\blackdiamond g\blackdiamond h\blackdiamond b\blackdiamond c$ & \text{H},\;\ref{8_18_2} & {\color{red} \ref{8_18_7}} \\
15 & \text{$I_2$} &$ e\blackdiamond a\blackdiamond f\blackdiamond d\blackdiamond g\blackdiamond h\blackdiamond b\blackdiamond c$ & \text{H},\;\ref{8_18_2} & {\color{red} \ref{8_18_7}} \\
16 & \text{$J_1$} &$ a\blackdiamond e\blackdiamond f\blackdiamond h\blackdiamond g$ & \text{$B_4$},\;\ref{8_18_2} & \\
17 & \text{$J_2$} &$ e\blackdiamond a\blackdiamond f\blackdiamond h\blackdiamond g$ & \text{$B_4$},\;\ref{8_18_2} & \\
18 & \text{$K_1$} & $ a\blackdiamond e\blackdiamond f\blackdiamond h\blackdiamond c\blackdiamond g$ & \text{$J_1$},\;\ref{8_18_1} & \\
19 & \text{$K_2$} & $ a\blackdiamond e\blackdiamond f\blackdiamond h\blackdiamond g\blackdiamond c $& \text{$J_1$},\;\ref{8_18_1} & \\
20 & \text{$L_1$} & $ e\blackdiamond a\blackdiamond f\blackdiamond h\blackdiamond c\blackdiamond g$ & \text{$J_2$},\;\ref{8_18_1} & \\
21 & \text{$L_2$} &$ e\blackdiamond a\blackdiamond f\blackdiamond h\blackdiamond g\blackdiamond c$ & \text{$J_2$},\;\ref{8_18_1} & \\
22 & \text{$M_1$} & $ a\blackdiamond e\blackdiamond f\blackdiamond d\blackdiamond h\blackdiamond c\blackdiamond g$ & \text{$K_1$},\;\ref{8_18_6} & \\
23 & \text{$M_2$} & $ a\blackdiamond e\blackdiamond f\blackdiamond h\blackdiamond d\blackdiamond c\blackdiamond g$ & \text{$K_1$},\;\ref{8_18_6} & \\
24 & \text{$N_1$} & $ a\blackdiamond e\blackdiamond b\blackdiamond f\blackdiamond d\blackdiamond h\blackdiamond c\blackdiamond g$ & \text{$M_1$},\;\ref{8_18_4},\;\ref{8_18_7} & {\color{blue} \ref{8_18_3},\;\ref{8_18_8}} \\
25 & \text{$N_2$} & $ a\blackdiamond e\blackdiamond f\blackdiamond b\blackdiamond d\blackdiamond h\blackdiamond c\blackdiamond g$ & \text{$M_1$},\;\ref{8_18_4},\;\ref{8_18_7} & {\color{blue} \ref{8_18_3},\;\ref{8_18_8}} \\
26 & \text{$O_1$} &$ a\blackdiamond e\blackdiamond b\blackdiamond f\blackdiamond h\blackdiamond d\blackdiamond c\blackdiamond g $& \text{$M_2$},\;\ref{8_18_4},\;\ref{8_18_7},\;\ref{8_18_8} & {\color{blue} \ref{8_18_3}} \\
27 & \text{$O_2$} & $ a\blackdiamond e\blackdiamond f\blackdiamond b\blackdiamond h\blackdiamond d\blackdiamond c\blackdiamond g$ & \text{$M_2$},\;\ref{8_18_4},\;\ref{8_18_7},\;\ref{8_18_8} & {\color{blue} \ref{8_18_3}} \\
28 & \text{$P_1$} & $ a\blackdiamond e\blackdiamond f\blackdiamond d\blackdiamond h\blackdiamond g\blackdiamond c$ & \text{$K_2$},\;\ref{8_18_3},\;\ref{8_18_6} & \\
29 & \text{$P_2$} &$  a\blackdiamond e\blackdiamond f\blackdiamond h\blackdiamond d\blackdiamond g\blackdiamond c$ & \text{$K_2$},\;\ref{8_18_3},\;\ref{8_18_6} & \\
30 & \text{$Q_1$} & $ a\blackdiamond e\blackdiamond b\blackdiamond f\blackdiamond d\blackdiamond h\blackdiamond g\blackdiamond c$&  \text{$P_1$},\;\ref{8_18_4},\;\ref{8_18_7} & {\color{blue} \ref{8_18_8}} \\
31 & \text{$Q_2$} & $ a\blackdiamond e\blackdiamond f\blackdiamond b\blackdiamond d\blackdiamond h\blackdiamond g\blackdiamond c$ &  \text{$P_1$},\;\ref{8_18_4},\;\ref{8_18_7} & {\color{blue} \ref{8_18_8}} \\
32 & \text{$R_1$} & $ a\blackdiamond e\blackdiamond b\blackdiamond f\blackdiamond h\blackdiamond d\blackdiamond g\blackdiamond c$ & \text{$P_2$},\;\ref{8_18_4},\;\ref{8_18_7},\;\ref{8_18_8} & {\color{blue} used all relations} \\
33 & \text{$R_2$} &$  a\blackdiamond e\blackdiamond f\blackdiamond b\blackdiamond h\blackdiamond d\blackdiamond g\blackdiamond c $& \text{$P_2$},\;\ref{8_18_4},\;\ref{8_18_7},\;\ref{8_18_8} & {\color{blue} used all relations} \\
34 & \text{$S_1$} & $ e\blackdiamond a\blackdiamond f\blackdiamond d\blackdiamond h\blackdiamond c\blackdiamond g$ & \text{$L_1$},\;\ref{8_18_6} & \\
35 & \text{$S_2$} & $ e\blackdiamond a\blackdiamond f\blackdiamond h\blackdiamond d\blackdiamond c\blackdiamond g$ & \text{$L_1$},\;\ref{8_18_6} & \\
36 & \text{$T_1$} & $ e\blackdiamond a\blackdiamond b\blackdiamond f\blackdiamond d\blackdiamond h\blackdiamond c\blackdiamond g$ & \text{$S_1$},\;\ref{8_18_7}& {\color{blue} \ref{8_18_3},\;\ref{8_18_4},\;\ref{8_18_8}} \\
37 & \text{$T_2$} & $ e\blackdiamond a\blackdiamond f\blackdiamond b\blackdiamond d\blackdiamond h\blackdiamond c\blackdiamond g$ & \text{$S_1$},\;\ref{8_18_7}& {\color{blue} \ref{8_18_3},\;\ref{8_18_4},\;\ref{8_18_8}} \\
38 & \text{$U_1$} &$  e\blackdiamond a\blackdiamond b\blackdiamond f\blackdiamond h\blackdiamond d\blackdiamond c\blackdiamond g$ & \text{$S_2$},\;\ref{8_18_7},\;\ref{8_18_8} & {\color{blue} \ref{8_18_3},\;\ref{8_18_4}} \\
39 & \text{$U_2$} &$  e\blackdiamond a\blackdiamond f\blackdiamond b\blackdiamond h\blackdiamond d\blackdiamond c\blackdiamond g$ & \text{$S_2$},\;\ref{8_18_7},\;\ref{8_18_8} & {\color{blue} \ref{8_18_3},\;\ref{8_18_4}} \\
40 & \text{$V_1$} &$ e\blackdiamond a\blackdiamond f\blackdiamond d\blackdiamond h\blackdiamond g\blackdiamond c$ & \text{$L_2$},\;\ref{8_18_3},\;\ref{8_18_6} & \\
41 & \text{$V_2$} & $e\blackdiamond a\blackdiamond f\blackdiamond h\blackdiamond d\blackdiamond g\blackdiamond c$ & \text{$L_2$},\;\ref{8_18_3},\;\ref{8_18_6} & \\
42 & \text{$W_1$} & $ e\blackdiamond a\blackdiamond b\blackdiamond f\blackdiamond d\blackdiamond h\blackdiamond g\blackdiamond c$ & \text{$V_1$},\;\ref{8_18_7} & {\color{blue} \ref{8_18_4},\;\ref{8_18_8}} \\
43 & \text{$W_2$} & $ e\blackdiamond a\blackdiamond f\blackdiamond b\blackdiamond d\blackdiamond h\blackdiamond g\blackdiamond c$ & \text{$V_1$},\;\ref{8_18_7} & {\color{blue} \ref{8_18_4},\;\ref{8_18_8}} \\
44 & \text{$X_1$} & $ e\blackdiamond a\blackdiamond b\blackdiamond f\blackdiamond h\blackdiamond d\blackdiamond g\blackdiamond c$ & \text{$V_2$},\;\ref{8_18_7},\;\ref{8_18_8} & {\color{blue} \ref{8_18_4}} \\
45 & \text{$X_2$} & $ e\blackdiamond a\blackdiamond f\blackdiamond b\blackdiamond h\blackdiamond d\blackdiamond g\blackdiamond c$ & \text{$V_2$},\;\ref{8_18_7},\;\ref{8_18_8} & {\color{blue} \ref{8_18_4}} \\
	\hline
\end{longtblr}
Linear orders on the generating set of $Q(8_{18})$ as in \text{$N_1$}, \text{$N_2$}, \text{$O_1$}, \text{$O_2$}, \text{$Q_1$}, \text{$Q_2$}, \text{$R_1$}, \text{$R_2$}, \text{$T_1$}, \text{$T_2$}, \text{$U_1$}, \text{$U_2$}, \text{$W_1$}, \text{$W_2$}, \text{$X_1$} and \text{$X_2$} in Table \ref{tbl8_18} could be extendable to biorders on $Q(8_{18})$.

\subsection{Knot $8_{19}$}\label{8_19}

Since the knot $8_{19}$ is the torus knot \text{$T(4,3)$}, by \cite[Theorem 7.2]{rss}, the quandle $Q(8_{19})$ is not right-orderable, and thus it is not biorderable.

\subsection{Knot $8_{20}$}\label{8_20}

\begin{figure}[H]
	\centering
	\includegraphics[scale=0.3]{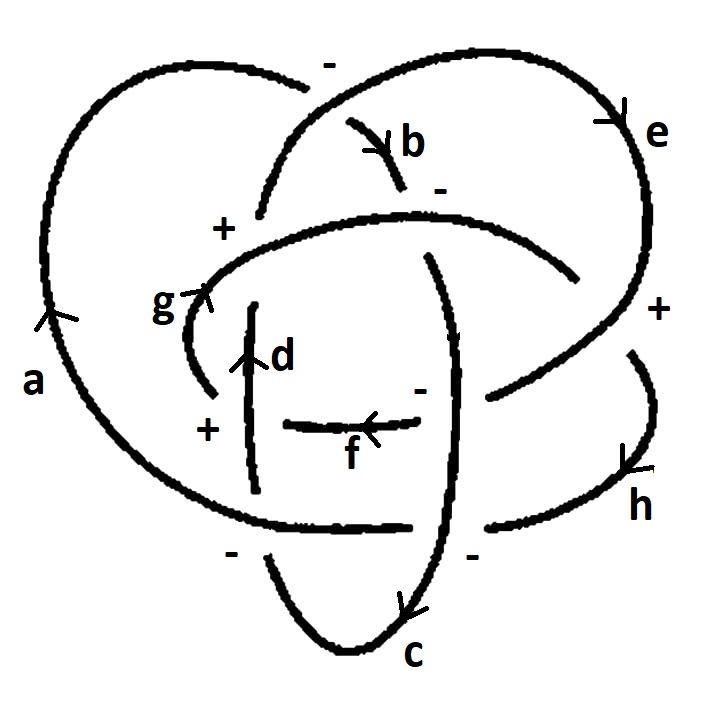}
	\caption{Knot $8_{20}$}
	\label{fig8_20}
\end{figure}

Consider the diagram of the knot $8_{20}$ as in Figure \ref{fig8_20}. The quandle $Q(8_{20})$ is generated by labelings of arcs of the diagram, and the relations are given by
\begin{align}
	b*e&=a&b\blackdiamond a\blackdiamond e \label{8_20_1}\\
	c*g&=b&c\blackdiamond b\blackdiamond g \label{8_20_2}\\
	d*g&=e&d\blackdiamond e\blackdiamond g \label{8_20_3}\\
	g*e&=h&g\blackdiamond h\blackdiamond e \label{8_20_4}\\
	f*d&=g&f\blackdiamond g\blackdiamond d \label{8_20_5}\\
	d*a&=c&d\blackdiamond c\blackdiamond a \label{8_20_6}\\
	a*c&=h&a\blackdiamond h\blackdiamond c \label{8_20_7}\\
	f*c&=e&f\blackdiamond e\blackdiamond c \label{8_20_8}
\end{align}

Since there is a quandle homomorphism  $\phi: Q(8_{20})\to\text{Core}(\mathbb{Z}_9)$ given by $a\mapsto 7, b\mapsto 6, c\mapsto 5, d\mapsto 0, e\mapsto 2, f\mapsto8, g\mapsto1$ and $h\mapsto3$, the generators $a$, $b$, $c$, $d$, $e$, $f$, $g$ and $h$ are all pairwise distinct. Applying the relations above ((\ref{8_20_1}) to (\ref{8_20_8})), we have the table as follows.

\begin{longtblr}[caption = {Orders that could be extendable to biorders on $Q(8_{20})$}, label = {tbl8_20}]
{colspec = {|c|c|l|l|l|}, rowhead = 2}
    \hline
	\text{Sr.} & \text{Label} & \text{Expression} & \text{Obtained by} & \text{{\color{red} Contradiction to} $\backslash$}\\
	\text{No.} &  & &  & \text{\color{blue} Compatible with}\\[2mm]
	\hline
1 & \text{A} & $d\blackdiamond e\blackdiamond h\blackdiamond g\blackdiamond f$ & \ref{8_20_3},\;\ref{8_20_4},\;\ref{8_20_5}& \\
2 & \text{$B_1$} &	$d\blackdiamond c\blackdiamond e\blackdiamond h\blackdiamond g\blackdiamond f$ & \text{A},\;\ref{8_20_8} & \\
3 & \text{$B_2$} &	$c\blackdiamond d\blackdiamond e\blackdiamond h\blackdiamond g\blackdiamond f$ & \text{A},\;\ref{8_20_8} & \\
4 & \text{$C_1$} &  $d\blackdiamond c\blackdiamond b\blackdiamond a\blackdiamond e\blackdiamond h\blackdiamond g\blackdiamond f$ & \text{$B_1$},\;\ref{8_20_2},\;\ref{8_20_1} & {\color{blue} \ref{8_20_6}},\;{\color{red} \ref{8_20_7}}	\\
5 & \text{$C_2$} &  $d\blackdiamond c\blackdiamond e\blackdiamond a\blackdiamond b\blackdiamond h\blackdiamond g\blackdiamond f$ & \text{$B_1$},\;\ref{8_20_2},\;\ref{8_20_1} & {\color{blue} \ref{8_20_6}},\;{\color{red} \ref{8_20_7}}	\\
6 & \text{$C_3$} &  $d\blackdiamond c\blackdiamond e\blackdiamond h\blackdiamond b\blackdiamond g\blackdiamond f$ & \text{$B_1$},\;\ref{8_20_2} & \\
7 & \text{$D_1$} &  $d\blackdiamond c\blackdiamond e\blackdiamond a\blackdiamond h\blackdiamond b\blackdiamond g\blackdiamond f$ & \text{$C_3$},\;\ref{8_20_1} & {\color{blue} \ref{8_20_6}},\;{\color{red} \ref{8_20_7}}	\\
8 & \text{$D_2$} &  $d\blackdiamond c\blackdiamond e\blackdiamond h\blackdiamond a\blackdiamond b\blackdiamond g\blackdiamond f$ & \text{$C_3$},\;\ref{8_20_1} & {\color{blue} \ref{8_20_6}},\;{\color{blue} \ref{8_20_7}}	\\
9 & \text{E} &  $b\blackdiamond a\blackdiamond c\blackdiamond d\blackdiamond e\blackdiamond h\blackdiamond g\blackdiamond f$ & \text{$B_2$},\;\ref{8_20_6},\;\ref{8_20_1} & {\color{red} \ref{8_20_2}},\;{\color{red} \ref{8_20_7}}	\\
	\hline
	\end{longtblr}

Linear order on the generating set of $Q(8_{20})$ as in \text{$D_2$} in Table \ref{tbl8_20} could be extendable to a biorder on $Q(8_{20})$.

\subsection{Knot $8_{21}$}\label{8_21}

\begin{figure}[H]
	\centering
	\includegraphics[scale=0.3]{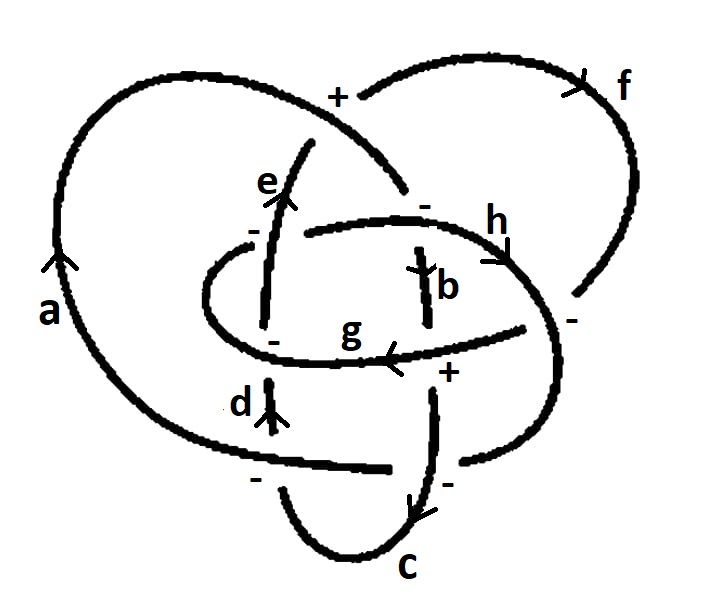}
	\caption{Knot $8_{21}$}
	\label{fig8_21}
\end{figure}

Consider the diagram of the knot $8_{21}$ as in Figure \ref{fig8_21}. The quandle $Q(8_{21})$ is generated by labelings of arcs of the diagram, and the relations are given by
\begin{align}
	e*a&=f&e\blackdiamond f\blackdiamond a \label{8_21_1}\\
	h*e&=g&h\blackdiamond g\blackdiamond e \label{8_21_2}\\
	b*h&=a&b\blackdiamond a\blackdiamond h \label{8_21_3}\\
	e*g&=d&e\blackdiamond d\blackdiamond g \label{8_21_4}\\
	b*g&=c&b\blackdiamond c\blackdiamond g \label{8_21_5}\\
	g*h&=f&g\blackdiamond f\blackdiamond h \label{8_21_6}\\
	a*c&=h&a\blackdiamond h\blackdiamond c \label{8_21_7}\\
	d*a&=c&d\blackdiamond c\blackdiamond a \label{8_21_8}
\end{align}

Since there is a quandle homomorphism  $\phi: Q(8_{21})\to\text{Core}(\mathbb{Z}_{15})$ given by $a\mapsto 1, b\mapsto 12, c\mapsto 0, d\mapsto 2, e\mapsto 10, f\mapsto7, g\mapsto6$ and $h\mapsto14$, the generators $a$, $b$, $c$, $d$, $e$, $f$, $g$ and $h$ are all pairwise distinct. Applying the relations above ((\ref{8_21_1}) to (\ref{8_21_8})), we have the table as follows.

\begin{longtblr}[caption = {Orders that could be extendable to biorders on $Q(8_{21})$}, label = {tbl8_21}]
{colspec = {|c|c|l|l|l|}, rowhead = 2}
            \hline
			\text{Sr.} & \text{Label} & \text{Expression} & \text{Obtained by} & \text{{\color{red} Contradiction to} $\backslash$}\\
			\text{No.} &  & &  & \text{\color{blue} Compatible with}\\[2mm]
			\hline
			1 & \text{A} & $b\blackdiamond a\blackdiamond h\blackdiamond c\blackdiamond d$ & \ref{8_21_7},\;\ref{8_21_8},\;\ref{8_21_3} & \\
			2 & \text{$B_1$} & $b\blackdiamond a\blackdiamond h\blackdiamond c\blackdiamond g\blackdiamond d\blackdiamond e$ & \text{A},\;\ref{8_21_5},\;\ref{8_21_4} & \\
			3 & \text{$B_2$} & $b\blackdiamond a\blackdiamond h\blackdiamond c\blackdiamond d\blackdiamond g\blackdiamond e$ & \text{A},\;\ref{8_21_5},\;\ref{8_21_2} & \\
			4 & \text{$C_1$} & $b\blackdiamond a\blackdiamond h\blackdiamond f\blackdiamond c\blackdiamond g\blackdiamond d\blackdiamond e$ & \text{$B_1$},\;\ref{8_21_6} & {\color{blue} \ref{8_21_1},\;\ref{8_21_2}} \\
			5 & \text{$C_2$} & $b\blackdiamond a\blackdiamond h\blackdiamond c\blackdiamond f\blackdiamond g\blackdiamond d\blackdiamond e$ & \text{$B_1$},\;\ref{8_21_6} & {\color{blue} \ref{8_21_1},\;\ref{8_21_2}} \\
			6 & \text{$D_1$} & $b\blackdiamond a\blackdiamond h\blackdiamond f\blackdiamond c\blackdiamond d\blackdiamond g\blackdiamond e$ & \text{$B_2$},\;\ref{8_21_6} & {\color{blue} \ref{8_21_1}},\;{\color{red}\ref{8_21_4}} \\
			7 & \text{$D_2$} & $b\blackdiamond a\blackdiamond h\blackdiamond c\blackdiamond f\blackdiamond d\blackdiamond g\blackdiamond e$ & \text{$B_2$},\;\ref{8_21_6} & {\color{blue} \ref{8_21_1}},\;{\color{red}\ref{8_21_4}} \\
			8 & \text{$D_3$} & $b\blackdiamond a\blackdiamond h\blackdiamond c\blackdiamond d\blackdiamond f\blackdiamond g\blackdiamond e$ & \text{$B_2$},\;\ref{8_21_6} & {\color{blue} \ref{8_21_1}},\;{\color{red}\ref{8_21_4}} \\
			\hline
\end{longtblr}

Linear orders on the generating set of $Q(8_{21})$ as in \text{$C_1$} and \text{$C_2$} in Table \ref{tbl8_21} could be extendable to biorders on $Q(8_{21})$.

\medskip

\begin{ack}
Vaishnavi Gupta is supported by the Junior Research Fellowship (JRF), University Grants Commission, India.
\end{ack}

\end{document}